\documentclass[11pt]{article}
\usepackage{amsfonts}
\usepackage{color}
\newtheorem{theo}{Theorem}[section]
\newtheorem{lem}[theo]{Lemma}
\newtheorem{cor}[theo]{Corollary}

\newtheorem{defi}{Definition}[section]

\newcommand{\mysection}[1]{\section{#1} \setcounter{equation}{0}}
\newcommand{\proof}{{\sc Proof.} \quad}
\newcommand{\proofc}{{\sc Proof} \ }
\newcommand{\be}{\begin{equation} \label}
\newcommand{\ee}{\end{equation}}
\newcommand{\bea}{\begin{eqnarray}\label}
\newcommand{\eea}{\end{eqnarray}}
\newcommand{\bas}{\begin{eqnarray*}}
\newcommand{\eas}{\end{eqnarray*}}
\newcommand{\bit}{\begin{itemize}}
\newcommand{\eit}{\end{itemize}}
\newcommand{\qed}{\hfill$\Box$ \vskip.2cm}
\newcommand{\nn}{\nonumber}
\newcommand{\R}{\mathbb{R}}
\newcommand{\N}{\mathbb{N}}
\newcommand{\pO}{\partial\Omega}

\newcommand{\eps}{\varepsilon}

\newcommand{\wto}{\rightharpoonup}

\newcommand{\hra}{\hookrightarrow}

\newcommand{\abs}{\\[5pt]}

\newcommand{\tme}{T_{max, \eps}}
\newcommand{\hatt}{\widehat{T}_\eps}
\newcommand{\io}{\int_\Omega}
\newcommand{\bom}{\overline{\Omega}}
\newcommand{\ue}{u_\eps}
\newcommand{\uex}{u_{\eps x}}
\newcommand{\uexx}{u_{\eps xx}}
\newcommand{\uexxx}{u_{\eps xxx}}
\newcommand{\uexxxx}{u_{\eps xxxx}}
\newcommand{\uet}{u_{\eps t}}
\newcommand{\ve}{v_\eps}
\newcommand{\vex}{v_{\eps x}}
\newcommand{\vexx}{v_{\eps xx}}
\newcommand{\vexxx}{v_{\eps xxx}}
\newcommand{\vexxxx}{v_{\eps xxxx}}
\newcommand{\vet}{v_{\eps t}}
\newcommand{\epss}{\eps_\star}

\newcommand{\F}{{\mathcal{F}}_\eps}
\newcommand{\D}{{\mathcal{D}}_\eps}

\newcommand{\Eo}{{\mathcal{E}}_{1,\eps}}
\newcommand{\Et}{{\mathcal{E}}_{2,\eps}}
\newcommand{\Do}{{\mathcal{D}}_{1,\eps}}
\newcommand{\Dt}{{\mathcal{D}}_{2,\eps}}

\newcommand{\ye}{y_\eps}
\newcommand{\ts}{T_\star}
\newcommand{\chis}{\chi_\star}

\newcommand{\chiss}{\chi_{\star\star}}
\newcommand{\oy}{\overline{y}}
\newcommand{\us}{u_{\star}}
\newcommand{\vs}{v_{\star}}
\newcommand{\xis}{\xi_{\star}}
\begin{document}
\title{Existence theory and qualitative analysis\\ for a fully cross-diffusive predator-prey system}
%
\author
{
Youshan Tao\footnote{taoys@sjtu.edu.cn}\\
{\small School of Mathematical Sciences, Shanghai Jiao Tong University,}\\
{\small Shanghai 200240, P.R.~China}
 \and
Michael Winkler\footnote{michael.winkler@math.uni-paderborn.de}\\
{\small Institut f\"ur Mathematik, Universit\"at Paderborn,}\\
{\small 33098 Paderborn, Germany} }
\date{}
\maketitle
\begin{abstract}
\noindent
This manuscript considers a Neumann initial-boundary value
problem for the         
predator-prey system
 \bas
        \left\{ \begin{array}{l}
    u_t = D_1 u_{xx} - \chi_1 (uv_x)_x + u(\lambda_1-u+a_1 v), \\[1mm]
    v_t = D_2 v_{xx} + \chi_2 (vu_x)_x + v(\lambda_2-v-a_2 u),
        \end{array} \right.
    \qquad \qquad (\star)
\eas
in an open bounded interval $\Omega$ as the spatial domain, where for
$i\in\{1,2\}$ the parameters $D_i, a_i, \lambda_i$ and $\chi_i$ are
positive.\abs
Due to the simultaneous appearance of two mutually interacting taxis-type cross-diffusive mechanisms,
one of which even being attractive, it seems unclear how far a solution theory can be built upon classical results
on parabolic evolution problems.
In order to nevertheless create an analytical setup capable of providing global existence results as well as detailed
information on qualitative behavior, this work pursues a strategy via parabolic regularization, in the course of which
($\star$) is approximated by means of certain fourth-order problems involving degenerate diffusion operators of thin film type.\abs
During the design thereof, a major challenge is related to the ambition to retain consistency with some fundamental
entropy-like structures formally associated with ($\star$); in particular, this will motivate the construction of
an approximation scheme including two free parameters which will finally be fixed in different ways, depending on
the size of $\lambda_2$ relative to $a_2 \lambda_1$.\abs        
Adequately coping with this will firstly yield a result on global existence of weak solutions for arbitrary choices
of the parameters in ($\star$) and arbitrarily large positive initial data from $H^1$,
and secondly allow for the conclusion that in both cases $\lambda_2 > a_2\lambda_1$ and $\lambda_2 \le a_2\lambda_1$,
the respectively obtained spatially homogeneous coexistence and prey-extinction states uphold their
global asymptotic stability properties well-known to be present in the corresponding ODE setting,
provided that both tactic sensitivities $\chi_1$ and $\chi_2$ are suitably small.\abs
 {\bf Key words:} cross-diffusion, thin-film equation, large time behavior\\
 {\bf MSC (2010):} 35K51, 35K52 (primary); 35K65, 35B40, 35Q92, 92C17 (secondary)
\end{abstract}
\newpage
\mysection{Introduction}
Taxis-type migration processes have been identified as impellent
mechanisms of crucial importance for the emergence of multifarious
dynamics in numerous biological systems at various levels of
complexity. Especially the destabilizing potential of attractive
chemotaxis, and its responsibility for striking experimental
findings, e.g.~on paradigmatic bacterial aggregation phenomena but
also on self-organization in more intricate frameworks, has received
considerable interest (\cite{KS}, \cite{hillen_painter2009}).
Accordingly, a meanwhile abundant mathematical literature has been
focusing on issues related to the ability of such taxis mechanisms
to enforce the formation of structures, not only in contexts of
simple paradigmatic Keller-Segel type chemotaxis models
(\cite{herrero_velazquez}, \cite{perthame_EJDE}, \cite{nagai},
\cite{nagai_senba_yoshida}, \cite{win_JMPA}) but also in some more
general triangular cross-diffusion systems embedding taxis processes
into more complex settings
(\cite{li_yan_blowup_two_species}, \cite{taowin_JEMS},
\cite{DLM}, \cite{win_ct_fluid_ARMA}, \cite{taowin_261}).\abs
In comparison to this, the knowledge seems much less developed in cases in which several taxis mechanisms are mutually coupled.
A prototypical situation of this type is addressed in \cite{tsyganov_prl}, where systems of the form
\be{00}
    \left\{ \begin{array}{l}
    u_t=D_1 \Delta u-\chi_1 \nabla \cdot (u\nabla v) +f(u, v), \\[1mm]
    v_t=D_2 \Delta v +\chi_2 \nabla \cdot (v \nabla u)+g(u, v),
      \end{array} \right.
\ee are proposed as refinements of classical reaction-diffusion
models for predator-prey interaction, supplemented by attractive
taxis of predators toward regions of increasing prey population
densities, and by repulsive cross-diffusive migration of prey
individuals downward population gradients of the predators (cf.~also
\cite{junping_shi_M3AS}, \cite{junping_shi_JDE}, \cite{tyutyunov_mmnp},
\cite{kareiva_odell} and \cite{hillen_lewis} for further modeling
aspects related to predator-taxis and prey-taxis). As documented in
\cite{tsyganov_prl}  and \cite{tyutyunov_mmnp}, formal linear
analysis indicates that the indeed introduction of such taxis
mechanisms enables {\em pursuit-evasion systems} of the form
(\ref{00}) to exhibit quite colorful wave-like solution behavior.
From a perspective of rigorous analysis, however, passing on to such
fully cross-diffusive systems, in which hence the collection of
migration mechanisms can no longer be arranged in the form of a
triangular diffusion operator, apparently amounts to entering quite
uncharted territories. In fact, studies on cross-diffusion systems
involving fully occupied diffusion matrices seem to essentially
concentrate on systems of Shigesada-Kawasaki-Teramoto type, for
which indeed a considerably comprehensive theory at least with
regard to questions of global solvability, but partially even going
beyond, could be developed (\cite{chen_jungel1},
\cite{chen_jungel2}, \cite{jungel_ARMA},
\cite{desvillettes_lepoutre_moussa_trescases},
\cite{lou_ni_jde1996}, \cite{taowin_PLMS}, \cite{lou_win}). However,
since the migration mechanisms therein do not only exhibit
structures evidently different from those in (\ref{00}),
but since they moreover, and yet more drastically, are exclusively of repulsive character,
such model classes can only be viewed as far relatives of (\ref{00}), with hence quite limited potential for accessibility
to similar techniques.
Accordingly, the analytical literature concerned with a doubly tactic system of type (\ref{00}) apparently reduces
to the single precedent \cite{taowin_double_crossdiff},
in which a method is designed that in the absence of sources, that is, in the case
when $f=g=0$, can be used to establish results on global existence and stabilization toward     
spatial averages in a corresponding one-dimensional boundary value poblem.\abs
{\bf Objectives and challenges.} \quad
The intention of this work is to address a system of type (\ref{00}) in a more realistic setting of
predator-prey evolution in which doubly tactic pursuit-evasion interplay is coupled to appropriate kinetics.
Concentrating on classical Lotka-Volterra interaction as a paradigmatic framework therefor,
we shall henceforth consider the apparently prototypical version of a fully cross-diffusive predator-prey system given by
\be{01}
        \left\{ \begin{array}{l}
    u_t = D_1 u_{xx} - \chi_1 (uv_x)_x + u(\lambda_1-u+a_1 v), \\[1mm]
    v_t = D_2 v_{xx} + \chi_2 (vu_x)_x + v(\lambda_2-v-a_2 u),
        \end{array} \right.
\ee
in an open bounded interval as the spatial domain, where for $i\in\{1,2\}$ the parameters
$D_i, a_i, \lambda_i$ and $\chi_i$ are positive.\abs
Our particular purpose consists firstly
in establishing a result on global existence of solutions
within a natural weak framework, and secondly in attempting to
undertake a basic step toward understanding qualitative effects of
the doubly cross-diffusive mechanisms in (\ref{01}) when accompanied by zero-order predator-prey interaction.
Indeed, standard results on local existence of smooth solutions to
Keller-Segel-type systems accounting for cross-diffusion exclusively
in one quantity (\cite{horstmann_win}, \cite{BBTW}) confirm that the
diffusion-induced relaxing behavior known from the situation when
$\chi_1=\chi_2=0$ persists at least during suitably small initial
time intervals when $\chi_1 \ne 0$ as long as $\chi_2$ yet vanishes.
On the other hand, the simultaneous appearance of taxis-type
cross-diffusion in both equations from (\ref{01}) may quite
drastically affect this picture in that when both taxis mechanisms
are attractive in the sense that $\chi_1>0>\chi_2$, then in general
not even local solutions can be expected to exist, even in the case
when the kinetic terms therein are completely disregarded and the
initial data belong to $(C^\infty(\bom))^2$
(\cite{taowin_double_crossdiff}).\abs
As already suggested by the outcome of the precedent work \cite{taowin_double_crossdiff},
such an instantaneously and tho\-roughly destabilizing role may not be played by doubly taxis-like interaction
in cases when only one of the cross-diffusive migration processes is attractive, with the other one being repulsive.
In particular,  
addressing the simplified variant of (\ref{01}) given by
\be{02}
        \left\{ \begin{array}{l}
    u_t = D_1 u_{xx} - \chi_1 (uv_x)_x, \\[1mm]
    v_t = D_2 v_{xx} + \chi_2 (vu_x)_x,
        \end{array} \right.
\ee with $\chi_1$ and $\chi_2$ both assumed positive, the main
results in \cite{taowin_double_crossdiff} assert that in an open
bounded domain $\Omega\subset\R$ and for all reasonably regular
nonnegative initial data, a corresponding Neumann-type
initial-boundary value problem possesses a globally defined
nonnegative weak solution, inter alia belonging to the space
$(L^\infty((0,\infty);L\log L(\Omega)))^2$. Moreover, this solution
is shown to become bounded in the pointwise sense at least
eventually, and that both of its components uniformly approach their
respective temporally constant spatial mean in the large time limit.
This inter alia indicates that the pattern-supporting potential of
the attractive taxis mechanism therein, well-known e.g.~as
generating nonconstant steady states in associated Keller-Segel
systems (\cite{schaaf}, \cite{lin_ni_takagi}), is insufficient to
enforce any structure formation on large time scales when
accompanied by repulsion as in (\ref{02}).\abs
As a fundamental technical obstacle, any analytical approach to
questions concerning solvability in (\ref{01}) needs to face the
circumstance that no general theory seems available which might at
least warrant local existence of some solutions. Hence forced to
construct solutions either from a very basic starting point
e.g.~within a suitable fixed point framework, or via approximation,
in this work we choose the latter type of ansatz by means of a
fourth-order regularization reminiscent of the well-studied thin
film equation (\cite{bernis_friedman},
\cite{beretta_bertsch_dalpasso},
\cite{giacomelli_grun_waiting_time}, \cite{otto}). Indeed, we shall
see that when carefully designed, beyond providing accessibility to
well-established local solution theory (\cite{amann}), this approach
will in quite a natural manner bring about the important advantage
of paving the way for a qualitative analysis through the derivation
of a priori estimates at the level of approximate solutions, and
might thereby turn out to be more appropriate than
e.g.~discretization-based methods (\cite{chen_jungel1}, \cite{chen_jungel2}, \cite{jungel_ARMA}) which in the
present context, beyond apparently unsolved problems already at the
stage of solvability, seem to offer somewhat less flexibility with
regard to testing procedures.\abs
{\bf Thin-film-type approximation.} \quad
To make our concrete strategy more precise, let us first specify the full
setting in which (\ref{01}) will be studied, and subsequently
consider
\be{0}
        \left\{ \begin{array}{lcll}
    u_t &=& D_1 u_{xx} - \chi_1 (uv_x)_x + u(\lambda_1-u+a_1 v),
    & x\in\Omega, \ t>0, \\[1mm]
    v_t &=& D_2 v_{xx} + \chi_2 (vu_x)_x + v(\lambda_2-v-a_2 u),
    & x\in\Omega, \ t>0, \\[1mm]
    & & \hspace*{-13mm}
    u_x=v_x=0,
    & x\in\pO, \ t>0, \\[1mm]
    & & \hspace*{-13mm}
    u(x,0)=u_0(x),
    \quad v(x,0)=v_0(x),
    & x\in\Omega,
        \end{array} \right.
\ee
in an open bounded interval $\Omega\subset\R$,
where the initial data $u_0$ and $v_0$, along with approximating families
$(u_{0\eps})_{\eps\in (0,1)}$ and $(v_{0\eps})_{\eps\in (0,1)}$,
will be assumed to be such that
\be{ie}
    \left\{ \begin{array}{l}
    u_0 \in W^{1,2}(\Omega)
    \mbox{ and }
    v_0\in W^{1,2}(\Omega)
    \quad \mbox{satisfy $u_0>0$ and $v_0>0$ in $\bom$, that} \\[1mm]
    u_{0\eps}\in C^5(\bom)
    \mbox{ and }
    v_{0\eps} \in C^5(\bom)
    \quad \mbox{for all $\eps\in (0,1)$ \quad with} \\[1mm]
    u_{0\eps x}=u_{0\eps xxx}= v_{0\eps x}=v_{0\eps xxx}=0
    \mbox{ on $\pO$ \quad for all $\eps\in (0,1)$ \quad and}\\[1mm]
    \frac{1}{2} \inf_\Omega u_0 \le u_{0\eps} \le u_0+1
    \mbox{ in $\Omega$ \quad and \quad}
    \io u_{0\eps x}^2 \le \io u_{0x}^2+1
    \quad \mbox{for all $\eps\in (0,1)$, \quad that} \\[1mm]
    \frac{1}{2} \inf_\Omega v_0 \le v_{0\eps} \le v_0+1
    \mbox{ in $\Omega$ \quad and \quad}
    \io v_{0\eps x}^2 \le \io v_{0x}^2+1
    \quad \mbox{for all $\eps\in (0,1)$, \quad and that} \\[1mm]
    u_{0\eps} \to u_0
    \mbox{ \ and \ }
    v_{0\eps} \to v_0
    \mbox{\quad a.e.~in $\Omega$ \quad as } \eps\searrow 0.
    \end{array} \right.
\ee
We henceforth fix any $\alpha\in (0,\frac{1}{2}]$, and with free parameters $n_1>0$ and $n_2>0$
to be specified below we consider the regularized parabolic system
\be{0eps}
        \left\{ \begin{array}{l}
    \uet = - \eps \Big( \frac{\ue^4}{\ue^{4-n_1}+\eps} \uexxx\Big)_x
    + \eps^\frac{\alpha}{2} (\ue^{-\alpha} \uex)_x
    + D_1 \uexx
    - \chi_1 \Big(\frac{\ue^{5-n_1}}{\ue^{4-n_1}+\eps} \vex \Big)_x
    + \frac{3\ue^3}{3\ue^2+\eps} \cdot (\lambda_1-\ue+a_1 \ve), \\[2mm]
    \hspace*{140mm}
    x\in\Omega, \ t>0, \\[2mm]
    \vet = - \eps \Big( \frac{\ve^4}{\ve^{4-n_2}+\eps} \vexxx\Big)_x
    + \eps^\frac{\alpha}{2} (\ve^{-\alpha} \vex)_x
    + D_2 \vexx
    + \chi_2 \Big( \frac{ \ve^{5-n_2}}{\ve^{4-n_2}+\eps} \uex \Big)_x
    + \frac{3\ve^3}{3\ve^2+\eps} \cdot (\lambda_2-\ve-a_2 \ue), \\[2mm]
    \hspace*{140mm}
    x\in\Omega, \ t>0,      
        \end{array} \right.
\ee
along with homogeneous Neumann-type boundary data and under the initial conditions given by
\be{0epsb}
    \left\{ \begin{array}{ll}
    \uex=\vex=\uexxx=\vexxx=0,
    &
    x\in\pO, \ t>0, \\[1mm]
    \ue(x,0)=u_{0\eps}(x),
    \quad \ve(x,0)=v_{0\eps}(x),
    \qquad
    &
    x\in\Omega,
        \end{array} \right.
\ee
for $\eps\in (0,1)$.
An introduction both of similar thin-film-type fourth-order diffusion operators and of corresponding
second-order fast diffusion corrections has already been underlying the analysis in \cite{taowin_double_crossdiff};
a difference of crucial importance in the present approach, however, consists in the circumstance that the
parameters $n_1$ and $n_2$ , which may be viewed as measuring a certain intermediate-scale degeneracy of the considered
fourth-order diffusion mechanisms, will be allowed to attain different values here:
While the development of our existence theory will merely require that $n_i\in [1,2]$ for $i\in\{1,2\}$,
our subsequent qualitative analysis will rely on the specific choices $(n_1,n_2)=(2,2)$ and $(n_1,n_2)=(2,1)$, respectively,
depending on the parameter setting dictated by the Lotka-Volterra interaction in (\ref{0}).\abs
{\bf Main results.} \quad
In fact, we shall firstly see that if we
restrict $n_1$ and $n_2$ so as to satisfy $n_i\in [1,2]$ for
$i\in\{1,2\}$, then the artificial terms introduced herein
conveniently cooperate with a quasi-entropy property formally enjoyed
by the functional
\be{entropy}
    \io u\ln u - \io u
    + \frac{\chi_1}{\chi_2} \io v\ln v
    - \frac{\chi_1}{\chi_2} \io v
\ee
when evaluated along suitably smooth trajectories of (\ref{0}).
This notice, generalizing an observation already made in the context of a
slightly different regularization for (\ref{02}) in
\cite{taowin_double_crossdiff}, will lead us to the first of our
main results, by asserting global solvability
without any restrictions on the system parameters in (\ref{0}), and for widely arbitrary initial data:
%
%
\begin{theo}\label{theo555}
  Let $\Omega\subset\R$ be a bounded open interval, for $i\in\{1,2\}$ let the constants
  $D_i, a_i, \lambda_i$ and $\chi_i$ be positive, and suppose that
  $u_0\in W^{1,2}(\Omega)$ and $v_0\in W^{1,2}(\Omega)$  are such that $u_0>0$ and $v_0>0$
  in $\bom$.
  Then there exist nonnegative functions $u$ and $v$ defined a.e.~in $\Omega\times (0,\infty)$ and satisfying
  \be{555.1}
    \{u,v\} \subset
    C^0_w([0,\infty);L^1(\Omega))
    \ \cap \ L^3_{loc}(\bom\times [0,\infty))
    \ \cap \ L^\frac{3}{2}_{loc}([0,\infty);W^{1,\frac{3}{2}}(\Omega))
    \ \cap \ L^\infty((0,\infty);L\log L(\Omega)),
  \ee
  which are such that $(u,v)$ forms a global weak solution of (\ref{0}) in the sense of Definition \ref{defi_weak} below.\\
  This solution can be obtained as the limit of solutions to (\ref{0eps})-(\ref{0epsb}) in that
  whenever $n_1\in [1,2]$, $n_2\in [1,2]$ and $\alpha\in (0,\frac{1}{2}]$ and $(u_{0\eps})_{\eps\in (0,1)}$
  and $(v_{0\eps})_{\eps\in (0,1)}$ satisfy (\ref{ie}),
  there exists
  $(\eps_j)_{j\in \N} \subset (0,1)$
  such that $\eps_j\searrow 0$ as $j\to\infty$, and that
  \be{555.2}
    \ue\to u
    \quad \mbox{as well as}  \quad
    \ve\to v
    \quad \mbox{a.e.~in $\Omega\times (0,\infty)$ \qquad as $\eps=\eps_j\searrow 0$.}
  \ee
\end{theo}
Beyond establishing the above existence result which with regard to its outcome essentially
parallels Theorem 1.1 in \cite{taowin_double_crossdiff}, our
analysis related to (\ref{entropy}) will moreover be organized in
such a way that a beneficent dependence of correspondingly gained
estimates on the sensitivities $\chi_1$ and $\chi_2$ can be cropped
out. Together with a rather straightforward observation identifying
a $(\chi_1,\chi_2)$-independent absorbing set in $(L^1(\Omega))^2$
for (\ref{0eps})-(\ref{0epsb}) (Lemma \ref{lem3}), the regularity information thus
provided by a quasi-entropy inequality associated with (\ref{entropy})
will imply eventual bounds for solutions to (\ref{0eps})-(\ref{0epsb}) in somewhat
stronger topologies, independent of the initial data and depending
on $\chi_1$ and $\chi_2$ in a favorably controllable manner (Lemma
\ref{lem19}). In a key step to be achieved in Lemma \ref{lem62},
these ultimate bounds will further be improved so as to become
manifest even in $(W^{1,2}(\Omega))^2$, provided that $\chi_1$ and
$\chi_2$ satisfy a smallness condition which, importantly, involves
essentially no knowledge on the initial data.\abs
This will bear fruit by implying global asymptotic stability of the nontrivial spatially homogeneous steady states
of (\ref{0}), in the case $\lambda_2>a_2 \lambda_1$ given by $(\us,\vs)$ with
\be{usvs}
    \us:=\frac{\lambda_1 + a_1\lambda_2}{1+a_1 a_2}
    \qquad \mbox{and} \qquad
    \vs:=\frac{\lambda_2 - a_2 \lambda_1}{1+a_1 a_2}
\ee
and by $(\lambda_1,0)$ if $\lambda_2 \le a_2 \lambda_1$, whenever $\chi_1$ and $\chi_2$ are suitably small:
In Section \ref{sect_asy1} addressing the former case, we shall see that then for small $\chi_1$ and $\chi_2$ and
any choice of initial data compatible with (\ref{ie}), a regularized variant of
\be{entropy2}
    \io \Big( u-\us - \us \ln \frac{u}{\us})
       + \frac{a_1}{a_2}    \io \Big( v-\vs - \vs \ln \frac{v}{\vs})
\ee
eventually plays the role of a genuine entropy functional for
(\ref{0eps})-(\ref{0epsb}), up to a regularization error, under a crucial
additional assumption on structural consistency of our
approximation, here expressed in the hypothesis that $n_1=n_2=2$
which is fortunately in compliance with Theorem \ref{theo555} (Lemma
\ref{lem25}). In conjunction with the compactness properties of
trajectories entailed by Lemma \ref{lem62}, this further dissipative
structure will lead us to the following second of our main results:
\begin{theo}\label{theo33}
  Let $\Omega\subset \R$ be an open bounded interval.
  Then given $D_1>0, D_2>0, a_1>0, a_2>0,\lambda_1>0$ and $\lambda_2>0$ fulfilling
  \bas
    \lambda_2> a_2\lambda_1,
  \eas
  one can find $\chiss>0$ such that if $\chi_1\in (0,\chiss)$ and $\chi_2\in (0,\chiss)$ as well as
  \bas
    n_1=2
    \quad \mbox{and} \quad
    n_2=2,
  \eas
  then for any choice of initial data fulfilling (\ref{ie}) the global weak solution $(u,v)$ of (\ref{0}) from
  Theorem \ref{theo555} has the properties that
  \be{reg}
    \{u,v\} \subset \in C^0(\bom\times (T,\infty)) \cap L^\infty(\Omega\times (0,\infty))
    \qquad \mbox{for some } T>0,
  \ee
  and that as $t\to\infty$,
  \bas
    u(\cdot,t) \to \us
    \quad \mbox{in } L^\infty(\Omega)
    \qquad \mbox{as well as} \qquad
    v(\cdot,t) \to \vs
    \quad \mbox{in } L^\infty(\Omega),
  \eas
  where $\us>0$ and $\vs>0$ are as in (\ref{usvs}).
\end{theo}
In the opposite parameter range where $\lambda_2 \le a_2\lambda_1$,
an accordingly modified analysis identifying a corresponding entropy property of a functional approximating
\be{entropy3}
    \io \Big( u-\lambda_1 - \lambda_1 \ln \frac{u}{\lambda_1}\Big)
    + \frac{a_1}{a_2} \io v
    +\frac{a_1}{2a_2 \lambda_2}\io v^2,
\ee
will finally reveal in Section \ref{sect_asy2}
the following analogue of the above statement under the assumption that yet $n_1=2$ but that
now $n_2$ satisfies the alternative consistency condition $n_2=1$, both still admissible in Theorem \ref{theo555}:
\begin{theo}\label{theo333}
  Whenever $\Omega\subset \R$ is an open bounded interval and
  $D_1>0, D_2>0, a_1>0, a_2>0,\lambda_1>0$ and $\lambda_2>0$ are such that
  \bas
    \lambda_2 \le a_2\lambda_1,
  \eas
  there exists $\chiss>0$ such that if $\chi_1\in (0,\chiss)$, $\chi_2\in (0,\chiss)$,
  \bas
    n_1=2
    \quad \mbox{and} \quad
    n_2=1
  \eas
  and if (\ref{ie}) holds, then the global weak solution $(u,v)$ of (\ref{0}) from
  Theorem \ref{theo555} satisfies (\ref{reg}) as well as
  \bas
    u(\cdot,t) \to \lambda_1
    \quad \mbox{in } L^\infty(\Omega)
    \qquad \mbox{and} \qquad
    v(\cdot,t) \to 0
    \quad \mbox{in } L^\infty(\Omega)
    \qquad \mbox{as } t\to\infty.
  \eas
\end{theo}
\vspace{0.3cm}
\mysection{Basic properties of the approximate problems}
To begin with, we recall that Amann's theory provides a basic statement on local existence and extensibility
of solutions to (\ref{0eps})-(\ref{0epsb}) in the following sense.
\begin{lem}\label{lem_loc}
  For $i\in\{1,2\}$ let $D_i>0, a_i>0, \lambda_i>0, \chi_i>0$ and $n_i\in (0,4)$,
  and suppose that (\ref{ie}) holds.
  Then for all $\eps\in (0,1)$
  there exist $\tme \in (0,\infty]$ and a pair $(\ue,\ve)$ of functions
  \bas
    \left\{ \begin{array}{l}
    \ue \in \bigcap_{s\in \frac{3}{2},2)} C^0([0,\tme);W^{s,2}(\Omega)) \cap C^{4,1} (\bom\times (0,\tme)) \qquad \mbox{and} \\[1mm]
        \ve \in \bigcap_{s\in \frac{3}{2},2)} C^0([0,\tme);W^{s,2}(\Omega)) \cap C^{4,1} (\bom\times (0,\tme)), 
    \end{array} \right.
  \eas
  satisfying $\ue>0$ and $\ve>0$ in $\bom\times [0,\tme)$, which are such that $(\ue,\ve)$ solves (\ref{0eps})-(\ref{0epsb})
  classically in $\Omega\times (0,\tme)$, and that
  \bea{ext}
    & & \hspace*{-15mm}
    \mbox{either $\tme=\infty$, \quad or \quad}  \nn\\
    & & \limsup_{t\nearrow \tme} \bigg\{
    \|\ue(\cdot,t)\|_{W^{2,2}(\Omega)} + \Big\| \frac{1}{\ue(\cdot,t)} \Big\|_{L^\infty(\Omega)}
    + \|\ve(\cdot,t)\|_{W^{2,2}(\Omega)} + \Big\| \frac{1}{\ve(\cdot,t)} \Big\|_{L^\infty(\Omega)} \bigg\}
    = \infty.
  \eea
\end{lem}
\proof
  This can be seen by adapting the argument from \cite[Lemma 2.1]{taowin_double_crossdiff} in a straightforward manner.
\qed
\subsection{Mass evolution. Absorbing sets in $L^1(\Omega)$}
Simple integration in (\ref{0eps})-(\ref{0epsb}) yields a basic information on the evolution of the total mass functionals
related to both solution components.
\begin{lem}\label{lem5}
  Let $D_i>0, a_i>0, \lambda_i>0, \chi_i>0$ and $n_i\in (0,4)$ for $i\in\{1,2\}$,
  and assume (\ref{ie}).
  Then for all $\eps\in (0,1)$, we have
  \be{5.1}
    \frac{d}{dt} \io \ue
    \le \Big( \lambda_1 + \frac{\sqrt{\eps}}{2\sqrt{3}}\Big) \cdot \io \ue
    - \io \ue^2
    + a_1 \io \ue\ve
    \qquad \mbox{for all } t\in (0,\tme)
  \ee
  and
  \be{5.2}
    \frac{d}{dt} \io \ve
    \le \Big(\lambda_2 + \frac{\sqrt{\eps}}{2\sqrt{3}}\Big) \io \ve
    - \io \ve^2
    - a_2 \io \ue \ve
    + \frac{a_2 \sqrt{\eps}}{2\sqrt{3}} \io \ue
    \qquad \mbox{for all } t\in (0,\tme).
  \ee
\end{lem}
\proof
  Writing $g_\eps(s):=\frac{3s^3}{3s^2+\eps}$ for $s\ge 0$ and $\eps\in (0,1)$, on the basis of (\ref{0eps})-(\ref{0epsb}) we compute
  \bea{5.3}
    \frac{d}{dt} \io \ue
    &=& \lambda_1 \io g_\eps(\ue)
    - \io g_\eps(\ue) \ue
    + a_1 \io g_\eps(\ue) \ve \nn\\
    &=& \lambda_1 \io \ue
    - \io \ue^2
    + a_1 \io \ue \ve \nn\\
    & & - \lambda_1 \io \Big(\ue-g_\eps(\ue)\Big)
    + \io \Big(\ue-g_\eps(\ue)\Big) \ue
    - a_1 \io \Big(\ue-g_\eps(\ue)\Big) \ve
  \eea
  for $t\in (0,\tme)$.
  Using that $[0,\infty) \ni s\mapsto s-g_\eps(s)$ is nonnegative and
  attains its maximum at $s=\sqrt{\frac{\eps}{3}}$ with extremal value
  $\frac{\sqrt{\eps}}{2\sqrt{3}}$, from (\ref{5.3}) we directly obtain (\ref{5.1}), whereas (\ref{5.2}) can be derived
  similarly.
\qed
By taking a suitable linear combination of the latter inequalities, we obtain some genuine $L^1$ estimates
which will later on play a fundamental role in our derivation of eventual bounds which do not depend on the
size of the initial data.
\begin{lem}\label{lem3}
  For $i\in\{1,2\}$ let $D_i>0, a_i>0, \lambda_i>0, \chi_i>0$ and $n_i\in (0,4)$,
  and assume that (\ref{ie}) holds.
  Then there exists a bounded function $m: (0,\infty) \to \R$ such that for all $\eps\in (0,1)$,
  \be{3.1}
    \io \ue(\cdot,t) + \io \ve(\cdot,t) \le m(t)
    \qquad \mbox{for all } t\in (0,\tme),
  \ee
  and such that
  \be{3.2}
    \limsup_{t\to\infty} m(t)
    \le m_\infty
    := \frac{|\Omega|}{2} \cdot \Big(\lambda_1 + \frac{1}{2\sqrt{3}} + 1
    + \max \Big\{ a_1^2,  1\Big\}\cdot \frac{a_2}{2\sqrt{3}} \Big)^2
    + \frac{|\Omega|}{2} \cdot \Big( \lambda_2 + \frac{1}{2\sqrt{3}}+1\Big)^2 \cdot
      \max \Big\{ a_1^2 \, , \, 1\Big\}
  \ee
\end{lem}
\proof
  We let
  \be{3.5}
    \beta:=\max \Big\{ a_1^2 \, , \, 1\Big\}
  \ee
  and combine (\ref{5.1}) with (\ref{5.2}) in estimating
  \bea{3.6}
    & & \hspace*{-20mm}
    \frac{d}{dt} \bigg\{ \io \ue + \beta \io \ve \bigg\}
    + \bigg\{ \io \ue + \beta\io \ve \bigg\} \nn\\
    &\le& \Big(\lambda_1 + \frac{1}{2\sqrt{3}} + 1 + \frac{\beta a_2}{2\sqrt{3}}\Big) \io \ue
    - \io \ue^2 \nn\\
    & & + \beta \Big(\lambda_2 + \frac{1}{2\sqrt{3}}+1\Big) \io \ve
    - \beta \io \ve^2 \nn\\
    & & + (a_1 - \beta a_2) \io \ue\ve
    \qquad \mbox{for all } t\in (0,\tme).
  \eea
  Here three applications of Young's inequality show that
  \bas
    \Big(\lambda_1 + \frac{1}{2\sqrt{3}} + 1 + \frac{\beta a_2}{2\sqrt{3}}\Big) \io \ue
    \le \frac{1}{2} \io \ue^2
    + \frac{|\Omega|}{2} \Big(\lambda_1 + \frac{1}{2\sqrt{3}} + 1 + \frac{\beta a_2}{2\sqrt{3}}\Big)^2
    \qquad \mbox{for all } t\in (0,\tme)
  \eas
  and
  \bas
    \beta\Big(\lambda_2 + \frac{1}{2\sqrt{3}}+1\Big) \io \ve
    \le \frac{\beta}{2} \io \ve^2
    + \frac{\beta |\Omega|}{2} \Big(\lambda_2 + \frac{1}{2\sqrt{3}}+1\Big)^2
    \qquad \mbox{for all } t\in (0,\tme)
  \eas
  as well as
  \bas
    (a_1 - \beta a_2) \io \ue\ve
    &\le& a_1 \io \ue\ve \\
    &\le& \frac{1}{2} \io \ue^2
    + \frac{a_1^2}{2} \io \ve^2 \\
    &\le& \frac{1}{2} \io \ue^2
    + \frac{\beta}{2} \io \ve^2
    \qquad \mbox{for all } t\in (0,\tme),
  \eas
  the latter relying on (\ref{3.5}).
  From (\ref{3.6}) we therefore obtain that
  \bas
    & & \hspace*{-10mm}
    \frac{d}{dt} \bigg\{ \io \ue + \beta \io \ve \bigg\}
    + \bigg\{ \io \ue + \beta\io \ve \bigg\} \nn\\
    &\le& c_4:=\frac{|\Omega|}{2} \Big(\lambda_1 + \frac{1}{2\sqrt{3}} + 1
    + \frac{\beta a_2}{2\sqrt{3}}\Big)^2
    + \frac{\beta |\Omega|}{2} \Big( \lambda_2 + \frac{1}{2\sqrt{3}}+1\Big)^2
    \qquad \mbox{for all } t\in (0,\tme)
  \eas
  and that accordingly, again due to (\ref{ie}),
  \bas
    \io \ue + \beta \io \ve
    &\le& \bigg\{ \io u_0 + \beta \io v_0 + (\beta+1)|\Omega| \bigg\} \cdot e^{-t}
    + c_4 \int_0^t e^{-(t-s)} ds \\
    &\le& \bigg\{ \io u_0 + \beta \io v_0 + (\beta+1)|\Omega| \bigg\} \cdot e^{-t}
    + c_4
    \qquad \mbox{for all } t\in (0,\tme).
  \eas
  As $\beta\ge 1$ by (\ref{3.5}), upon an abvious choice of $m$ this asserts both (\ref{3.1}) and (\ref{3.2})
  in this case.
\qed
\subsection{Global extensibility in the regularized problems}
In order to show that the solutions from Lemma \ref{lem_loc} are actually global in time, in view
of (\ref{ext}) our goal will be to establish a priori bounds, throughout this part possibly depending on $\eps$,
for $\ue$ and $\ve$, and for $\frac{1}{\ue}$ and $\frac{1}{\ve}$, in $W^{2,2}(\Omega)$
and $L^\infty(\Omega)$, respectively.
In a first step toward this, but moreover also later on in the context of our qualitative analysis (see Lemma \ref{lem62}),
we will make use of the following interpolation inequality that has extensively been exploited
in the analysis of the thin film equation, albeit mostly for slightly different purposes there
(\cite{bernis} and \cite{beretta_bertsch_dalpasso}).
\begin{lem}\label{lem63}
  Let $\beta\in\R$ be such that $\beta\ne 1$. Then
  \be{63.1}
    \io \varphi^{\beta-2} \varphi_x^4  \le \frac{9}{(\beta-1)^2} \io \varphi^\beta \varphi_{xx}^2
  \ee
  for all $\varphi\in C^2(\bom)$ which are such that $\varphi>0$ in $\bom$ and $\varphi_x=0$ on $\pO$.
\end{lem}
\proof
  We integrate by parts and use the Cauchy-Schwarz inequality to see that
  \bas
    \io \varphi^{\beta-2} \varphi_x^4
    = - \frac{3}{\beta-1} \io \varphi^{\beta-1} \varphi_x^2 \varphi_{xx}
    \le \frac{3}{|\beta-1|} \bigg\{ \io \varphi^{\beta-2} \varphi_x^4 \bigg\}^\frac{1}{2} \cdot
    \bigg\{ \io \varphi^\beta \varphi_{xx}^2 \bigg\}^\frac{1}{2},
  \eas
  from which (\ref{63.1}) immediately follows.
\qed
In the present context, our first application thereof will guarantee a certain consistency of the artificial
second-order fast diffusion term in (\ref{0eps})-(\ref{0epsb}) with regard to the evolution of the $H^1$ norms of solutions.
Already in the following basic statement concerning this, to be immediately utilized in Lemma \ref{lem2}
but recalled also later on in the crucial Lemma \ref{lem62},
our standing assumption that $\alpha \le \frac{1}{2}$ plays a major role.
\begin{lem}\label{lem1}
  For $i\in\{1,2\}$ let $D_i>0, a_i>0, \lambda_i>0, \chi_i>0$ and $n_i\in (0,4)$,
  and assume (\ref{ie}).
  The for all $\eps\in (0,1)$ and any $t\in (0,\tme)$,
  \bea{1.1}
    & & \hspace*{-20mm}
    \frac{1}{2} \frac{d}{dt} \io \uex^2
    + \eps \io \frac{\ue^4}{\ue^{4-n_1}+\eps} \uexxx^2
    + \frac{D_1}{2} \io \uexx^2 \nn\\
    &\le& \chi_1 \io \Big( \frac{\ue^{5-n_1}}{\ue^{4-n_1}+\eps} \vex\Big)_x \uexx
    + 3\lambda_1 \io \uex^2
    + \frac{a_1^2}{2D_1} \io \ue^2 \ve^2
  \eea
  and
  \bea{1.2}
    & & \hspace*{-20mm}
    \frac{1}{2} \frac{d}{dt} \io \vex^2
    + \eps \io \frac{\ve^4}{\ve^{4-n_2}+\eps} \vexxx^2
    + \frac{D_2}{2} \io \vexx^2 \nn\\
    &\le& -\chi_2 \io \Big( \frac{\ve^{5-n_2}}{\ve^{4-n_2}+\eps} \uex\Big)_x \vexx
    + 3\lambda_2 \io \vex^2
    + \frac{a_2^2}{2D_2} \io \ue^2 \ve^2.
  \eea
\end{lem}
\proof
  We multiply the first equation in (\ref{0eps}) by $-\uexx$ and integrate by parts to see that
  \bea{1.4}
    & & \hspace*{-20mm}
    \frac{1}{2} \frac{d}{dt} \io \uex^2
    + \eps \io \frac{\ue^4}{\ue^{4-n_1}+\eps} \uexxx^2
    + D_1 \io \uexx^2 \nn\\
    &=& - \eps^\frac{\alpha}{2} \io (\ue^{-\alpha} \uex)_x \uexx
    + \chi_1 \io \Big( \frac{\ue^{5-n_1}}{\ue^{4-n_1}+\eps} \vex\Big)_x \uexx \nn\\
    & & - \io \frac{3\ue^3}{3\ue^2+\eps} (\lambda_1-\ue+a_1 \ve) \uexx
    \qquad \mbox{for all } t\in (0,\tme),
  \eea
  where another integration by parts, followed by an application of Lemma \ref{lem63}, shows that
  \bea{1.5}
    - \io (\ue^{-\alpha} \uex)_x \uexx
    &=& - \io \ue^{-\alpha} \uexx^2
    + \alpha \io \ue^{-\alpha-1} \uex^2 \uexx \nn\\
    &=& - \io \ue^{-\alpha} \uexx^2
    + \frac{\alpha(\alpha+1)}{3} \io \ue^{-\alpha-2} \uex^4 \nn\\[2mm]
    &\le& 0
    \qquad \mbox{for all } t\in (0,\tme),
  \eea
  because $\frac{\alpha(\alpha+1)}{3} \cdot \frac{9}{(\alpha+1)^2} = \frac{3\alpha}{\alpha+1} \le 1$ thanks to our assumption
  that $\alpha\le \frac{1}{2}$.
  Apart from that, observing that for all $\eps\in (0,1)$ and any $s\ge 0$ we have
  \bas
    \frac{d}{ds} \Big( \frac{3s^3}{3s^2 + \eps}\Big)
    = \frac{9s^4 + 9\eps s^2}{(3s^2+\eps)^2}
    \le \frac{9s^2}{3s^2+\eps}
    \le 3
  \eas
  and
  \bas
    \frac{d}{ds} \Big( \frac{3s^4}{3s^2 + \eps}\Big)
    = \frac{18s^5 + 12\eps s^3}{(3s^2+\eps)^2}
    \ge 0,
  \eas
  in the rightmost expression in (\ref{1.4}) we may once more rely on integration by parts to estimate
  \bea{1.6}
    - \io \frac{3\ue^3}{3\ue^2+\eps} (\lambda_1-\ue) \uexx
    &=& \lambda_1 \io \Big( \frac{3\ue^3}{3\ue^2+\eps}\Big)_x \uex
    - \io \Big( \frac{3\ue^4}{3\ue^2+\eps}\Big)_x \uex \nn\\
    &\le& 3 \lambda_1 \io \uex^2
    \qquad \mbox{for all } t\in (0,\tme).
  \eea
  As finally
  \bas
    - a_1 \io \frac{3\ue^3}{3\ue^2+\eps} \ve \uexx
    &\le& \frac{D_1}{2} \io \uexx^2
    + \frac{a_1^2}{2D_1} \io \Big( \frac{3\ue^3}{3\ue^2+\eps} \Big)^2 \ve^2 \\
    &\le& \frac{D_1}{2} \io \uexx^2
    + \frac{a_1^2}{2D_1} \io \ue^2 \ve^2
    \qquad \mbox{for all } t\in (0,\tme)
  \eas
  by Young's inequality, in view of (\ref{1.5}) and (\ref{1.6}) we readily infer (\ref{1.1}) from (\ref{1.4}).
  The inequality (\ref{1.2}) can be derived similarly.
\qed
Under the additional requirements that $1\le n_1\le 2$ and $1\le
n_2\le 2$, for each fixed $\eps\in (0,1)$ the integrals in
(\ref{1.1}) and (\ref{1.2}) originating from the cross-diffusive
interaction in (\ref{0eps}) can conveniently be estimated in terms
of the respective higher-order dissipative contributions, thus
leading to the following $\eps$-dependent bound with respect to the
norm in $W^{1,2}(\Omega)$.
\begin{lem}\label{lem2}
  Let $D_i>0, a_i>0, \lambda_i>0, \chi_i>0$ and $n_i\in [1,2]$ for $i\in\{1,2\}$,
  and let (\ref{ie}) be valid.
  Then for all $\eps\in (0,1)$ and $T>0$ there exists $C(\eps,T)>0$ such that
  \be{2.1}
    \io \uex^2(x,t) dx
    + \io \vex^2(x,t) dx
    \le C(\eps,T)
    \qquad \mbox{for all } t\in (0,\hatt)
  \ee
  and
  \be{7.1}
    \|\ue(\cdot,t)\|_{L^\infty(\Omega)}
    + \|\ve(\cdot,t)\|_{L^\infty(\Omega)} \le C(\eps,T)
    \qquad \mbox{for all } t\in (0,\hatt),
  \ee

  where $\hatt:=\min\{T,\tme\}$.
\end{lem}
\proof
  In the cross-diffusive term on the right of (\ref{1.1}), we integrate by parts and use Young's inequality to find that
  for all $t\in (0,\tme)$,
  \bea{2.3}
    \chi_1 \io \Big( \frac{\ue^{5-n_1}}{\ue^{4-n_1}+\eps} \vex\Big)_x \uexx
    &=& - \chi_1 \io \frac{\ue^{5-n_1}}{\ue^{4-n_1}+ \eps} \vex \uexxx \nn\\
    &\le& \eps \io \frac{\ue^4}{\ue^{4-n_1} +\eps} \uexxx^2
    + \frac{\chi_1^2}{4\eps} \io \frac{\ue^{6-2n_1}}{\ue^{4-n_1}+\eps} \vex^2,
  \eea
  where since $1\le n_1 \le 2$, again due to Young's inequality we obtain that
  \bea{2.4}
    \io \frac{\ue^{6-2n_1}}{\ue^{4-n_1}+\eps} \vex^2
    \le \io \ue^{2-n_1} \vex^2
    &\le& \io (\ue+1) \vex^2\nn\\
    &\le& \|\vex\|_{L^\infty(\Omega)}^2 \cdot \io (\ue+1)
    \qquad \mbox{for all } t\in (0,\tme).
  \eea
  Here we use that
  \be{2.44}
    \sup_{t\in (0,\hatt)} \bigg\{\io \ue + \io \ve \bigg\} <\infty
  \ee
  by Lemma \ref{lem3},
  which in conjunction with the Gagliardo-Nirenberg inequality and Young's
  inequality shows that (\ref{2.4}) implies that with some $c_1>0$ and $c_2>0$ we have
  \bea{2.5}
    \frac{\chi_1^2}{4\eps} \io \frac{\ue^{6-2n_1}}{\ue^{4-n_1}+\eps} \vex^2
    &\le& c_1 \|\vex\|_{L^\infty(\Omega)}^2 \nn\\
    &\le& c_2 \|\vexx\|_{L^2(\Omega)} \|\vex\|_{L^2(\Omega)} \nn\\
    &\le& \frac{D_2}{2} \io \vexx^2
    + \frac{c_2^2}{2 D_2} \io \vex^2
    \qquad \mbox{for all } t\in (0,\hatt).
  \eea
  Also the last summand in (\ref{1.1}) can be controlled by means of Young's inequality and the Gagliardo-Nirenberg
  inequality, which thanks to (\ref{2.44}) namely provide $c_3>0$ and $c_4>0$ fulfilling
  \bas
    \frac{a_1^2}{2D_1} \io \ue^2 \ve^2
    &\le& \frac{a_1^2}{4D_1} \io \ue^4
    + \frac{a_1^2}{4D_1} \io \ve^4 \nn\\
    &\le& c_3 \bigg\{ \io \uex^2 \bigg\} \cdot \bigg\{ \io \ue \bigg\}^2
    + c_3 \bigg\{ \io \ue \bigg\}^4
    + c_3 \bigg\{ \io \vex^2 \bigg\} \cdot \bigg\{ \io \ve \bigg\}^2
    + c_3 \bigg\{ \io \ve \bigg\}^4 \nn\\
    &\le& c_4 \io \uex^2 + c_4 \io \vex^2
    \qquad \mbox{for all } t\in (0,\hatt).
  \eas
  Combined with (\ref{2.3}) and (\ref{2.5}), this shows that (\ref{1.1}) entails the inequality
  \bas
    \frac{1}{2} \frac{d}{dt} \io \uex^2
    + \frac{D_1}{2} \io \uexx^2
    \le \frac{D_2}{2} \io \vexx^2
    + (3\lambda_1 + c_4) \io \uex^2
    + \Big( \frac{c_2^2}{2D_2} + c_4 \Big) \io \vex^2
    \qquad \mbox{for all } t\in (0,\hatt),
  \eas
  so that applying quite a similar procedure to (\ref{1.2}) we infer the existence of $c_5>0$ such that
  \bas
    \frac{d}{dt} \bigg\{ \io \uex^2 + \io \vex^2 \bigg\}
    \le c_5 \cdot \bigg\{ \io \uex^2 + \io \vex^2 \bigg\}
    \qquad \mbox{for all } t\in (0,\hatt),
  \eas
  which upon integration yields (\ref{2.1}) and thereby also implies (\ref{7.1}) due to the continuity
  of the embedding $W^{1,2}(\Omega) \hra L^\infty(\Omega)$.
\qed
The following conclusion of Lemma \ref{lem2} parallels that of Lemma 2.3 from \cite{taowin_double_crossdiff}
in its statement, but its derivation proceeds in a slightly different manner.
\begin{lem}\label{lem6}
  Let $D_i>0, a_i>0, \lambda_i>0, \chi_i>0$ and $n_i\in [1,2]$ for $i\in\{1,2\}$,
  and assume (\ref{ie}).
  Then given any $\eps\in (0,1)$, for all $T>0$ one can find $C(\eps,T)>0$ such that again writing
  $\hatt:=\min\{T,\tme\}$ we have
  \be{6.1}
    \io \frac{1}{\ue^2(x,t)} dx
    + \io \frac{1}{\ve^2(x,t)} dx \le C(\eps,T)
    \qquad \mbox{for all } t\in (0,\hatt).
  \ee
\end{lem}
\proof
  On the basis of (\ref{0eps})-(\ref{0epsb}) and several integrations by parts, we compute
  \bea{6.3}
    \frac{1}{2} \frac{d}{dt} \io \frac{1}{\ue^2}
    &=& \eps \io \frac{1}{\ue^3} \cdot \Big( \frac{\ue^4}{\ue^{4-n_1}+\eps} \uexxx\Big)_x
    -\eps^\frac{\alpha}{2} \io \frac{1}{\ue^3} \cdot (\ue^{-\alpha}
    \uex)_x\nn\\
    & &   - D_1 \io \frac{1}{\ue^3} \uexx
    + \chi_1 \io \frac{1}{\ue^3} \cdot \Big( \frac{\ue^{5-n_1}}{\ue^{4-n_1}+\eps} \vex \Big)_x \nn\\
    & & -3\io \frac{1}{3\ue^2+\eps} \cdot (\lambda_1-\ue+a_1\ve) \nn\\
    &=& 3\eps \io \frac{1}{\ue^{4-n_1}+\eps} \uex \uexxx
       -3\eps^\frac{\alpha}{2} \io \frac{\uex^2}{\ue^{4+\alpha}}\nn\\
    & &- 3D_1 \io \frac{\uex^2}{\ue^4}
    + 3\chi_1 \io \frac{1}{\ue^{n_1-1} (\ue^{4-n_1}+\eps)} \uex \vex \nn\\
    & & -3\io \frac{1}{3\ue^2+\eps} \cdot (\lambda_1-\ue+a_1\ve)
    \qquad \mbox{for all } t\in (0,\tme),
  \eea
  where one more integration by parts, followed by an application of Young's inequality, shows that
  for all $t\in (0,\tme)$,
  \bea{6.4}
    3\eps \io \frac{1}{\ue^{4-n_1}+\eps} \uex \uexxx
    &=& -3\eps \io \frac{1}{\ue^{4-n_1}+\eps} \uexx^2
    + 3(4-n_1)\eps \io \frac{\ue^{3-n_1}}{(\ue^{4-n_1}+\eps)^2} \uex^2  \uexx \nn\\
    &\le& -\frac{3\eps}{2} \io \frac{1}{\ue^{4-n_1}+\eps} \uexx^2
    + \frac{3(4-n_1)^2\eps}{2} \io \frac{\ue^{6-2n_1}}{(\ue^{4-n_1}+\eps)^3} \uex^4.
  \eea
  In order to estimate the rightmost summand herein appropriately,
  we recall that according to Lemma \ref{lem2} there exist
  $c_1=c_1(\eps,T)>0$ and $c_2=c_2(\eps,T)>0$ such that
  \be{6.44}
    \io \uex^2 + \io \vex^2 \le c_1
    \qquad \mbox{for all } t\in (0,\hatt)
  \ee
  and
  \be{6.45}
    \ue(x,t) + \ve(x,t) \le c_2
    \qquad \mbox{for all $x\in\Omega$ and } t\in (0,\hatt).
  \ee
  The latter, namely ensures that
  \be{6.5}
    \io \frac{1}{\ue^{4-n_1}+\eps} \uexx^2
    \ge \frac{1}{c_2^{4-n_1} + \eps} \io \uexx^2
    \qquad \mbox{for all } t\in (0,\hatt),
  \ee
  while trivially
  \bea{6.6}
    \io \frac{\ue^{6-2n_1}}{(\ue^{4-n_1}+\eps)^3} \uex^4
    &\le& \io \frac{(\ue^{4-n_1}+\eps)^\frac{6-2n_1}{4-n_1}}{(\ue^{4-n_1}+\eps)^3} \uex^4 \nn\\
    &=& \io (\ue^{4-n_1}+\eps)^{-\frac{6-n_1}{4-n_1}} \uex^4 \nn\\
    &\le& \eps^{-\frac{6-n_1}{4-n_1}} \io \uex^4
    \qquad \mbox{for all } t\in (0,\tme)
  \eea
  by nonnegativity of $\ue$.
  Since the Gagliardo-Nirenberg inequality provides $c_3>0$ such that
  \bas
    \io \varphi^4 \le c_3 \|\varphi_{x}\|_{L^2(\Omega)} \|\varphi\|_{L^2(\Omega)}^3
    \qquad \mbox{for all } \varphi\in W_0^{1,2}(\Omega),
  \eas
  along with (\ref{6.44}) this enables us to invoke Young's inequality in making sure that
  \bea{6.7}
    \frac{3(4-n_1)^2\eps}{2} \io \frac{\ue^{6-2n_1}}{(\ue^{4-n_1}+\eps)^3} \uex^4
    &\le& \frac{3(4-n_1)^2}{2} \eps^{-\frac{2}{4-n_1}} \io \uex^4 \nn\\
    &\le& \frac{3(4-n_1)^2 c_1^\frac{3}{2} c_3}{2} \eps^{-\frac{2}{4-n_1}} \cdot \bigg\{ \io \uexx^2 \bigg\}^\frac{1}{2}
        \nn\\
    &\le& \frac{3\eps}{2(c_2^{4-n_1}+\eps)} \io \uexx^2
    + c_4
    \qquad \mbox{for all } t\in (0,\hatt)
  \eea
  with $c_4=c_4(\eps,T):=\frac{3}{8} (4-n_1)^4 c_1^3 c_3^2 (c_2^{4-n_1}+\eps) \cdot \eps^{-\frac{8-n_1}{4-n_1}}$.\\
  Next, in the second last summand in (\ref{6.3}) we use Young's inequality to see that again thanks to (\ref{6.44}),
  \bea{6.8}
    3\chi_1 \io \frac{1}{\ue^{n_1-1} (\ue^{4-n_1}+\eps)} \uex \vex
    &\le& 3D_1 \io \frac{\uex^2}{\ue^4}
    + \frac{3\chi_1^2}{4D_1} \io \frac{\ue^{6-2n_1}}{(\ue^{4-n_1}+\eps)^2} \vex^2 \nn\\
    &=& 3D_1 \io \frac{\uex^2}{\ue^4}
    + \frac{3\chi_1^2}{4D_1}
    \io \frac{\ue^{6-2n_1}}{(\ue^{4-n_1}+\eps)^\frac{6-2n_1}{4-n_1} (\ue^{4-n_1}+\eps)^\frac{2}{4-n_1}} \vex^2 \nn\\
    &\le& 3D_1 \io \frac{\uex^2}{\ue^4}
    + \frac{3\chi_1^2}{4D_1 \eps^\frac{2}{4-n_1}} \io \vex^2 \nn\\
    &\le& 3D_1 \io \frac{\uex^2}{\ue^4}
    + c_5
    \qquad \mbox{for all } t\in (0,\hatt)
  \eea
  if we let $c_5=c_5(\eps,T):= \frac{3\chi_1^2  c_1}{4D_1} \eps^{-\frac{2}{4-n_1}}$.\\
  Finally, in the last integral in (\ref{6.3}) we may use (\ref{6.45}) to achieve the rough estimate
  \bas
    |\lambda_1 - \ue + a_1\ve|
    \le \lambda_1 + (a_1+1) c_2
    \qquad \mbox{in } \Omega \times (0,\hatt),
  \eas
  which ensures that
  \bas
    -3\io \frac{1}{3\ue^2+\eps} \cdot (\lambda_1-\ue+a_1\ve)
    \le c_6=c_6(\eps,T):=\frac{3(\lambda_1+(a_1+1)c_1)|\Omega|}{\eps}
    \qquad \mbox{for all } t\in (0,\hatt).
  \eas
  In combination with (\ref{6.4}), (\ref{6.5}), (\ref{6.7}) and (\ref{6.8}), this shows that (\ref{6.3}) entails the
  inequality
  \bas
    \frac{1}{2} \frac{d}{dt} \io \frac{1}{\ue^2} \le c_4 + c_5 + c_6
    \qquad \mbox{for all } t\in (0,\hatt)
  \eas
  and thereby yields the claimed estimate for $\ue$, while that for $\ve$ can be derived similarly.
\qed
As well-known from \cite{taowin_double_crossdiff}, the estimates from \ref{lem2} and Lemma \ref{lem6} in combination yield
pointwise lower bounds for both solution components:
\begin{lem}\label{lem55}
  Let $D_i>0, a_i>0, \lambda_i>0, \chi_i>0$ and $n_i\in [1,2]$ for $i\in\{1,2\}$,
  and let (\ref{ie}) hold.
  Then for all $\eps\in (0,1)$ and $T>0$ there exists $C(\eps,T)>0$ such that
  \be{55.1}
    \ue(x,t) \ge C(\eps,T)
    \quad \mbox{and} \quad
    \ve(x,t) \ge C(\eps,T)
    \qquad \mbox{for all $x\in\Omega$ and } t\in (0,\hatt),
  \ee
  where again $\hatt:=\min\{T,\tme\}$.
\end{lem}
\proof
  Based on Lemma \ref{lem2} with Lemma \ref{lem6}, this can be obtained by verbatim copying the argument from
  \cite[Lemma 2.4]{taowin_double_crossdiff}.
\qed
Again referring to a corresponding argument from \cite{taowin_double_crossdiff}, we may refrain from giving details
concerning the derivation of $H^2$ estimates of the following flavor.
\begin{lem}\label{lem54}
  Let $D_i>0, a_i>0, \lambda_i>0, \chi_i>0$ and $n_i\in [1,2]$ for $i\in\{1,2\}$,
  and let (\ref{ie}) hold.
  Then for all $\eps\in (0,1)$ and $T>0$ one can find $C(\eps,T)>0$ fulfilling
  \be{54.1}
    \io \uexx^2(x,t)dx
    + \io \vexx^2(x,t)dx
    \le C(\eps,T)
    \qquad \mbox{for all } t\in (0,\hatt),
  \ee
  where again $\hatt:=\min\{T,\tme\}$.
\end{lem}
\proof
  This estimate can be derived by testing the first two equations in (\ref{0eps}) against $\uexxxx$ and $\vexxxx$
  and estimating all appearing ill-signed lower-order integrals in terms of the
  respective fourth-order dissipative summands, using thanks to Lemma \ref{lem55} and Lemma \ref{lem2}
  the latter are bounded
  from below by suitable positive multiples of $\io \uexxxx^2$ and $\io \vexxxx^2$.
  For a corresponding argument in a specialized setting without the kinetic terms from (\ref{0eps}) we refer to
  \cite[Lemma 2.5]{taowin_double_crossdiff}, and we may omit detailing the minor adaptations necessary in the present framework.
\qed
In view of (\ref{ext}), collecting our $\eps$-dependent regularity information we can now make sure that indeed all our approximate
solutions actually exist globally in time.
\begin{lem}\label{lem56}
  Let $D_i>0, a_i>0, \lambda_i>0, \chi_i>0$ and $n_i\in [1,2]$ for $i\in\{1,2\}$,
  and suppose that (\ref{ie}) is satisfied.
  Then for all $\eps\in (0,1)$ we have $\tme=\infty$; that is
  the solution $(\ue,\ve)$ of (\ref{0eps})-(\ref{0epsb}) from Lemma \ref{lem_loc} is global in time.
\end{lem}
\proof
  This directly results from a combination of Lemma \ref{lem54},
  Lemma \ref{lem55}, Lemma \ref{lem3} and (\ref{ext}).
\qed
\mysection{A quasi-entropy structure reminiscent of (\ref{entropy})}
We next intend to analyze a regularization-adapted modification of the functional in (\ref{entropy}).
Here in contrast to the simplified setup in \cite{taowin_double_crossdiff} in which the absence of kinetic terms implies a genuine
Lyapunov property, in the present context the appearance of the zero-order terms in (\ref{0eps}) requires
substantial additional efforts.
These will be prepared by the following basic observation.
\begin{lem}\label{lem9}
  Let $\nu\in [0,2]$. Then
  \be{9.1}
    \frac{s^\nu}{3s^2+\eps} \le \frac{1}{2} \cdot \Big(\frac{\nu}{3}\Big)^\frac{\nu}{2}
    \cdot (2-\nu)^\frac{2-\nu}{2} \eps^{-\frac{2-\nu}{2}}
    \qquad \mbox{for all $s\ge 0$ and any } \eps\in (0,1).
  \ee
\end{lem}
\proof
  In the case $\nu=2$, (\ref{9.1}) follows from the simple estimate
  \bas
    \frac{s^\nu}{3s^2+\eps} \le \frac{s^\nu}{3s^2}=\frac{1}{3}
    \qquad \mbox{for all $s\ge 0$ and any } \eps\in (0,1).
  \eas
  When $\nu\in [0,2)$, the function $\varphi(s):=\frac{s^\nu}{3s^2+\eps}, \ s\ge 0$, satisfies
  $\varphi'(s)=\frac{-3(2-\nu) s^{\nu+1} + \nu\eps s^{\nu-1}}{(3s^2+\eps)^2}$
  for all $s>0$, whence $\varphi$ attains its maximum at $s_0:=\sqrt{\frac{\nu\eps}{3(2-\nu)}}$ with
  \bas
    \varphi(s_0)
    = \frac{\big(\frac{\nu}{3(2-\nu)}\big)^\frac{\nu}{2} \eps^\frac{\nu}{2}}{\frac{\nu\eps}{2-\nu} + \eps}
    = \frac{1}{2} \cdot \Big(\frac{\nu}{3}\Big)^\frac{\nu}{2} \cdot (2-\nu)^\frac{2-\nu}{2} \eps^{-\frac{2-\nu}{2}},
  \eas
  which verifies (\ref{9.1}) also in this case.
\qed
We can now make sure that if
in addition to the above we assume $n_1\ge 1$ and $n_2\ge 1$, then our approximation in (\ref{0eps}) indeed
cooperates with a fundamental structural property of (\ref{0}) in the following sense.
\begin{lem}\label{lem8}
  Assume that $D_i>0, a_i>0, \lambda_i>0$ and $n_i\in [1,2]$ for $i\in\{1,2\}$.
  Then there exists $C>0$ such that whenever $\chi_1>0$ and $\chi_2>0$ and (\ref{ie}) holds,
  for all $\eps\in (0,1)$ the functions $\F$ and $\D$ defined by
  \bea{F}
    \hspace*{-12mm}
    \F(t)
    \!\!\!
    &:=&
    \!\!\!
    \io \ue(\cdot,t) \ln \ue(\cdot,t)
    - \io \ue(\cdot,t)
    + \frac{\eps}{(3-n_1)(4-n_1)} \io \frac{1}{\ue^{3-n_1}(\cdot,t)} \nn\\
    & & +
    \frac{\chi_1}{\chi_2} \io \ve(\cdot,t) \ln \ve(\cdot,t)
    - \frac{\chi_1}{\chi_2} \io \ve(\cdot,t)
    + \frac{\chi_1 \eps}{(3-n_2)(4-n_2)\chi_2} \io \frac{1}{\ve^{3-n_2}(\cdot,t)},
    \quad t\ge 0,
  \eea
  and
  \bea{D}
    \hspace*{-5mm}
    \D(t)
    \!\!\!
    &:=&
    \!\!\!
    \frac{D_1}{2} \io \frac{\uex^2(\cdot,t)}{\ue(\cdot,t)}
    + \eps \io \ue^{n_1-1}(\cdot,t) \uexx^2(\cdot,t)
    + D_1 \eps \io \frac{\uex^2(\cdot,t)}{\ue^{5-n_1}(\cdot,t)} \nn\\
    & & + \frac{\chi_1 D_2}{2\chi_2} \io \frac{\vex^2(\cdot,t)}{\ve(\cdot,t)}
    + \frac{\chi_1 \eps}{\chi_2} \io \ve^{n_2-1}(\cdot,t) \vexx^2(\cdot,t)
    + \frac{\chi_1 D_2 \eps}{\chi_2} \io \frac{\vex^2(\cdot,t)}{\ve^{5-n_2}(\cdot,t)},
    \qquad t>0,
  \eea
  satisfy
  \bea{8.1}
    \frac{d}{dt} \F(t)
    + \D(t)
    \le C \cdot \Big( 1 + \frac{\chi_1}{\chi_2} \Big) \cdot
    \Bigg\{ 1 + \bigg\{ \io \ue(\cdot,t) \bigg\}^7 + \bigg\{ \io \ve(\cdot,t) \bigg\}^7 \Bigg\}
    \qquad \mbox{for all } t>0.
  \eea
\end{lem}
\proof
  For $i\in\{1,2\}$ abbreviating
  \bas
    L_i(s):= s\ln s - s + \frac{\eps}{(4-n_i)(3-n_i)} \cdot \frac{1}{s^{3-n_i}},
    \qquad s>0,
  \eas
  we see that
  \be{8.01}
    L_i'(s) = \ln s - \frac{\eps}{4-n_i} \cdot \frac{1}{s^{4-n_i}}
    \qquad \mbox{for all } s>0
  \ee
  and
  \be{8.02}
    L_i''(s) = \frac{1}{s} + \frac{\eps}{s^{5-n_i}}
    \qquad \mbox{for all } s>0,
  \ee
  and we observe that thus, in particular,
  \be{8.03}
    \frac{s^4}{s^{4-n_i}+\eps} \cdot L_i''(s)
    = s^{n_i-1}
    \qquad \mbox{for all } s>0
  \ee
  and hence
  \be{8.2}
    \frac{d^2}{ds^2} \Big( \frac{s^4}{s^{4-n_i}+\eps} L_i''(s)\Big)
    = (n_i-1)(n_i-2) s^{n_i-3} \le 0
    \qquad \mbox{for all } s>0
  \ee
  according to our assumption that $n_i\in [1,2]$.
  In order to make appropriate use of this, we now go back to (\ref{0eps})-(\ref{0epsb}) and integrate by parts in computing
  \bea{8.3}
    \frac{d}{dt} \io L_1(\ue)
    &=& \io L_1'(\ue) \uet \nn\\
    &=& - \eps \io L_1'(\ue) \cdot \Big(\frac{\ue^4}{\ue^{4-n_1}+\eps} \uexxx\Big)_x
    + \eps^\frac{\alpha}{2} \io L_1'(\ue) \cdot (\ue^{-\alpha} \uex)_x \nn\\
    & & + D_1 \io L_1'(\ue) \uexx \nn\\
    & & - \chi_1 \io L_1'(\ue) \cdot \Big( \frac{\ue^{5-n_1}}{\ue^{4-n_1}+\eps} \vex \Big)_x \nn\\
    & & + \io L_1'(\ue) \cdot \frac{3\ue^3}{3\ue^2+ \eps} \cdot (\lambda_1 - \ue + a_1\ve) \nn\\[1mm]
    &=& \eps \io \frac{\ue^4}{\ue^{4-n_1}+\eps} \cdot L_\eps''(\ue) \uex\uexxx
    - \eps^\frac{\alpha}{2} \io \ue^{-\alpha} L_1''(\ue) \uex^2 \nn\\
    & & - D_1 \io L_1''(\ue) \uex^2 \nn\\
    & & + \chi_1 \io \frac{\ue^{5-n_1}}{\ue^{4-n_1}+\eps} L_1''(\ue) \uex\vex \nn\\
    & & + \io L_1'(\ue) \cdot \frac{3\ue^3}{3\ue^2+\eps} (\lambda_1 - \ue + a_1 \ve)
    \qquad \mbox{for all } t>0,
  \eea
  where (\ref{8.02}) and (\ref{8.01}) directly yield
  \be{8.33}
    -\eps^\frac{\alpha}{2} \io \ue^{-\alpha} L_1''(\ue) \uex^2 \le 0
    \qquad \mbox{for all } t>0
  \ee
  and
  \be{8.4}
    -D_1 \io L_1''(\ue) \uex^2
    = - D_1 \io \frac {\uex^2}{\ue}
    - D_1 \eps \io \frac{\uex^2}{\ue^{5-n_1}}
    \qquad \mbox{for all } t>0
  \ee
  as well as
  \be{8.5}
    \chi_1 \io \frac{\ue^{5-n_1}}{\ue^{4-n_1}+\eps} L_1''(\ue) \uex\vex
    = \chi_1 \io \uex \vex
    \qquad \mbox{for all } t>0
  \ee
  and
  \bea{8.6}
    \hspace*{-10mm}
    \io L_1'(\ue) \cdot \frac{3\ue^3}{3\ue^2+\eps} (\lambda_1 - \ue + a_1 \ve)
    &=& 3\io \frac{\ue^3 \ln \ue}{3\ue^2+\eps} (\lambda_1 - \ue + a_1\ve) \nn\\
    & & - \frac{3\eps}{4-n_1} \io \frac{\ue^{n_1-1}}{3\ue^2+\eps} (\lambda_1-\ue+a_1\ve)
    \qquad \mbox{for all } t>0.
  \eea
  In the first summand on the right of (\ref{8.3}), we integrate by parts two more times to find that thanks to
  (\ref{8.2}) and (\ref{8.03}),
  \bea{8.7}
    \eps \io \frac{\ue^4}{\ue^{4-n_1}+\eps} \cdot L_1''(\ue) \uex\uexxx
    &=& - \eps  \io \frac{\ue^4}{\ue^{4-n_1}+\eps} \cdot L_1''(\ue) \uexx^2 \nn\\
    & & - \eps \io \frac{d}{ds} \Big(\frac{s^4}{s^{4-n_1}+\eps} L_1''(s)\Big) \bigg|_{s=\ue} \cdot \uex^2 \uexx \nn\\
    &=& - \eps  \io \frac{\ue^4}{\ue^{4-n_1}+\eps} \cdot L_1''(\ue) \uexx^2 \nn\\
    & & + \frac{\eps}{3} \io \frac{d^2}{ds^2} \Big(\frac{s^4}{s^{4-n_1}+\eps} L_1''(s)\Big) \bigg|_{s=\ue}
        \cdot \uex^4 \nn\\
    &\le& - \eps  \io \frac{\ue^4}{\ue^{4-n_1}+\eps} \cdot L_1''(\ue) \uexx^2 \nn\\
    &=& - \eps \io \ue^{n_1-1} \uexx^2
    \qquad \mbox{for all } t>0,
  \eea
  so that collecting (\ref{8.3})-(\ref{8.7}) we obtain
  \bas
    & & \hspace*{-20mm}
    \frac{d}{dt} \io L_1(\ue)
    + \eps \io \ue^{n_1-1} \uexx^2
    + D_1 \io \frac{\uex^2}{\ue}
    + D_1 \eps \io \frac{\uex^2}{\ue^{5-n_1}} \\
    &\le& \chi_1 \io \uex\vex  \\
    & & + 3\io \frac{\ue^3 \ln \ue}{3\ue^2 + \eps} (\lambda_1-\ue+a_1\ve)
    - \frac{3\eps}{4-n_1} \io \frac{\ue^{n_1-1}}{3\ue^2+\eps} (\lambda_1-\ue+a_1\ve)
    \qquad \mbox{for all } t>0.
  \eas
  As similarly
  \bas
    & & \hspace*{-20mm}
    \frac{d}{dt} \io L_2(\ve)
    + \eps \io \ve^{n_2-1} \vexx^2
    + D_2 \io \frac{\vex^2}{\ve}
    + D_2 \eps \io \frac{\vex^2}{\ve^{5-n_2}} \\
    &\le& - \chi_2 \io \uex\vex \\
    & & + 3\io \frac{\ve^3 \ln \ve}{3\ve^2 + \eps} (\lambda_2-\ve-a_2 \ue)
    - \frac{3\eps}{4-n_2} \io \frac{\ve^{n_2-1}}{3\ve^2+\eps} (\lambda_2-\ve-a_2 \ue)
    \qquad \mbox{for all } t>0,
  \eas
  on taking a suitable linear combination of the latter two relations we can achieve a cancellation of the respective
  cross-diffusive contributions and thereby obtain the inequality
  \bea{8.8}
    \hspace*{-10mm}
    \frac{d}{dt} \F(t) + \D(t)
    &=& \frac{d}{dt} \bigg\{ \io L_1(\ue) + \frac{\chi_1}{\chi_2} \io L_2(\ve)\bigg\}
    + \D(t) \nn\\
    &\le& \!\!\! - \frac{D_1}{2} \io \frac{\uex^2}{\ue}
    - \frac{\chi_1 D_2}{2\chi_2} \io \frac{\vex^2}{\ve} \nn\\
    & & \!\!\! + 3\io \frac{\ue^3 \ln \ue}{3\ue^2+\eps} (\lambda_1-\ue+a_1\ve)
    - \frac{3\eps}{4-n_1} \io \frac{\ue^{n_1-1}}{3\ue^2+\eps} (\lambda_1-\ue+a_1\ve) \nn\\
    & & \!\!\! + \frac{3\chi_1}{\chi_2} \io \frac{\ve^3 \ln \ve}{3\ve^2+\eps} (\lambda_2-\ve-a_2\ue)
    - \frac{3\chi_1 \eps}{(4-n_2)\chi_2} \io \frac{\ve^{n_2-1}}{3\ve^2+\eps} (\lambda_2-\ve-a_2\ue)
  \eea
  for all $t>0$.
  Here clearly
  \bas
    3\io \frac{\ue^3 \ln \ue}{3\ue^2+\eps} (\lambda_1-\ue+a_1\ve)
    \le \lambda_1 \int_{\{\ue\ge 1\}} \ue\ln\ue
    - \int_{\{\ue\le 1\}} \ue^2 \ln \ue
    + a_1 \int_{\{\ue\ge 1\}} \ue\ve \ln\ue
     \eas
for all $t>0$,  so that using the estimates
  \be{8.999}
    \ln\xi \le 2\sqrt{\xi}
    \quad \mbox{for $\xi\ge 1$}
    \qquad \mbox{and} \qquad
    \xi^2\ln\xi \ge -\frac{1}{2e}
    \quad \mbox{for } \xi \in (0,1],
  \ee
  by means of Young's inequality we infer that
  \bea{8.9}
  \hspace*{-5mm}
    3\io \frac{\ue^3 \ln \ue}{3\ue^2+\eps} (\lambda_1-\ue+a_1\ve)
    &\le& 2\lambda_1 \io \ue^\frac{3}{2}
    + \frac{|\Omega|}{2e}
    + 2a_1 \io \ue^\frac{3}{2} \ve \nn\\
    &\le& 2\bigg\{ \lambda_1 |\Omega|+ a_1 \io \ve \bigg\} \cdot \|\ue\|_{L^\infty(\Omega)}^\frac{3}{2}
    + \frac{|\Omega|}{2e}
    \qquad \mbox{for all } t>0.
  \eea
  Since the Gagliardo-Nirenberg inequality provides $c_1>0$ fulfilling
  \bas
    \|\varphi\|_{L^\infty(\Omega)}^3
    \le c_1 \|\varphi_x\|_{L^2(\Omega)}^\frac{3}{2} \|\varphi\|_{L^2(\Omega)}^\frac{3}{2}
    + c_1 \|\varphi\|_{L^2(\Omega)}^3
    \qquad \mbox{for all } \varphi\in W^{1,2}(\Omega).
  \eas
  by several applications of Young's inequality we can herein estimate
  \bea{8.99}
    2\bigg\{ \lambda_1 |\Omega|+ a_1 \io \ve \bigg\} \cdot \|\ue\|_{L^\infty(\Omega)}^\frac{3}{2}
    &\le& 2c_1 \bigg\{ \lambda_1 |\Omega|+ a_1 \io \ve \bigg\}
    \cdot \|(\sqrt{\ue})_x\|_{L^2(\Omega)}^\frac{3}{2} \cdot \bigg\{ \io \ue \bigg\}^\frac{3}{4} \nn\\
    & & + 2c_1 \bigg\{ \lambda_1 |\Omega|+ a_1 \io \ve \bigg\} \cdot \bigg\{ \io \ue \bigg\}^\frac{3}{2} \nn\\
    &\le& 2D_1 \|(\sqrt{\ue})_x\|_{L^2(\Omega)}^2
    + \frac{2c_1^4}{D_1^3} \cdot \bigg\{ \lambda_1 |\Omega|+ a_1 \io \ve \bigg\}^4 \cdot \bigg\{ \io \ue \bigg\}^3  \nn\\
    & & + 2c_1 \bigg\{ \lambda_1 |\Omega|+ a_1 \io \ve \bigg\} \cdot \bigg\{ \io \ue \bigg\}^\frac{3}{2} \nn\\
    &=& \frac{D_1}{2} \io \frac{\uex^2}{\ue}
    + \frac{2c_1^4}{D_1^3} \cdot \bigg\{ \lambda_1 |\Omega| + a_1 \io \ve \bigg\}^4 \cdot \bigg\{ \io \ue \bigg\}^3  \nn\\
    & & + 2c_1 \bigg\{ \lambda_1 |\Omega|+ a_1 \io \ve \bigg\} \cdot \bigg\{ \io \ue \bigg\}^\frac{3}{2}  \nn\\
    &\le&  \frac{D_1}{2} \io \frac{\uex^2}{\ue} \nn\\
    & & + \Big(\frac{2c_1^4}{D_1^3} + 2c_1\Big) \cdot \bigg\{ \lambda_1 |\Omega| + a_1\io \ve \bigg\}^7 \nn\\
    & & + \Big(\frac{2c_1^4}{D_1^3} + 2c_1\Big) \cdot \bigg\{ \io \ue\bigg\}^7 + 2c_1 \nn\\
    &\le&  \frac{D_1}{2} \io \frac{\uex^2}{\ue} \nn\\
    & & + 2^6 \cdot \Big(\frac{2c_1^4}{D_1^3} + 2c_1\Big) \cdot \Bigg\{ \lambda_1^7 |\Omega|^7+ a_1^7 \bigg\{ \io \ve \bigg\}^7 \Bigg\}
    \nn\\
    & & + \Big(\frac{2c_1^4}{D_1^3} + 2c_1\Big) \cdot \bigg\{ \io \ue\bigg\}^7 + 2c_1
  \eea
  for all $t>0$.
  Next, noting that again thanks to (\ref{8.999}),
  \bas
   - \frac{\ve^3 \ln \ve}{3\ve^2+\eps} \cdot \ve
    \le -\frac{1}{3} \ve^2 \ln \ve \le \frac{1}{6e}
    \qquad \mbox{if } \ve \le 1,
  \eas
  that
  \bas
    \frac{\ve^3 \ln \ve}{3\ve^2+\eps}  (\lambda_2-\ve)\le 0
        \qquad \mbox{if } \ve \ge \max\{1,\lambda_2\}
  \eas
  and that
  \bas
    \frac{\ve^3 \ln \ve}{3\ve^2+\eps}
        \le \frac{1}{3} \ve\ln \ve \le \frac{2}{3}\ve^\frac{3}{2}\le
        \frac{2}{3} \lambda_2^\frac{3}{2}
        \qquad \mbox{if } 1<\ve<\lambda_2,
  \eas
  writing $c_2:=\max\{\frac{1}{6e} \, , \, \frac{2}{3}\lambda_2^\frac{5}{2}\}$ we see that
  \bea{8.10}
    \hspace*{-10mm}
    \frac{3\chi_1}{\chi_2} \io \frac{\ve^3 \ln \ve}{3\ve^2+\eps} (\lambda_2-\ve-a_2\ue)
    &\le& \frac{3\chi_1  c_2 |\Omega|}{\chi_2}
    - \frac{\chi_1 a_2}{\chi_2} \int_{\{\ve\le 1\}} \ue\ve\ln\ve \nn\\
    &\le& \frac{3\chi_1 c_2 |\Omega|}{\chi_2}
    + \frac{\chi_1 a_2}{e \chi_2} \io \ue \nn\\
    &\le& \frac{3\chi_1  c_2 |\Omega|}{\chi_2}
    + \frac{\chi_1 a_2}{e \chi_2} \cdot \bigg\{ \io \ue \bigg\}^7
    + \frac{\chi_1 a_2}{e \chi_2}
    \qquad \mbox{for all } t>0
  \eea
   due to the inequality $-\xi\ln \xi\le \frac{1}{e}$ for $\xi \in (0,1]$.
  In the last and the third last summand in (\ref{8.8}) we make use of Lemma \ref{lem9}, which when applied to
  $\nu=n_1$ entails that
  \bea{8.11}
    - \frac{3\eps}{4-n_1} \io \frac{\ue^{n_1-1}}{3\ue^2+\eps} (\lambda_1-\ue+a_1\ve)
    &\le& \frac{3\eps}{4-n_1} \io \frac{\ue^{n_1}}{3\ue^2+\eps} \nn\\
    &\le& \frac{3\eps}{4-n_1} \cdot \frac{1}{2} \Big(\frac{n_1}{3}\Big)^\frac{n_1}{2} (2-n_1)^\frac{2-n_1}{2}
    \eps^{-\frac{2-n_1}{2}} |\Omega| \nn\\
    &=& \frac{[3(2-n_1)]^\frac{2-n_1}{2} n_1^\frac{n_1}{2} |\Omega|}{2(4-n_1)}  \cdot \eps^\frac{n_1}{2} \nn\\
    &\le& \frac{[3(2-n_1)]^\frac{2-n_1}{2} n_1^\frac{n_1}{2} |\Omega|}{2(4-n_1)}
    \qquad \mbox{for all } t>0
  \eea
  because of our restriction that $\eps<1$.
  Twice more employing Lemma \ref{lem9}, with $\nu=n_2$ and $\nu=n_2-1$, respectively, furthermore shows that
  again due to Young's inequality,
  \bea{8.12}
    & & \hspace*{-40mm}
    - \frac{3\chi_1 \eps}{(4-n_2)\chi_2} \io \frac{\ve^{n_2-1}}{3\ve^2+\eps} (\lambda_2-\ve-a_2\ue) \nn\\
    &\le& \frac{3\chi_1 \eps}{(4-n_2)\chi_2} \io \frac{\ve^{n_2}}{3\ve^2+\eps}
    + \frac{3a_2 \chi_1 \eps}{(4-n_2)\chi_2} \io \frac{\ve^{n_2-1}}{3\ve^2+\eps} \cdot \ue \nn\\
    &\le& \frac{3\chi_1 \eps}{(4-n_2)\chi_2} \cdot \frac{1}{2}\Big(\frac{n_2}{3}\Big)^\frac{n_2}{2}
    (2-n_2)^\frac{2-n_2}{2} \eps^{-\frac{2-n_2}{2}} |\Omega| \nn\\
    & & + \frac{3a_2 \chi_1 \eps}{(4-n_2)\chi_2} \cdot \frac{1}{2}\Big(\frac{n_2-1}{3}\Big)^\frac{n_2-1}{2}
    (3-n_2)^\frac{3-n_2}{2} \eps^{-\frac{3-n_2}{2}} \io \ue \nn\\
    &=& \frac{[3(2-n_2)]^\frac{2-n_2}{2} n_2^\frac{n_2}{2} \chi_1 |\Omega|}{2(4-n_2)\chi_2} \eps^\frac{n_2}{2} \nn\\
    & & + \frac{[3(3-n_2)]^\frac{3-n_2}{2} (n_2-1)^\frac{n_2-1}{2} a_2 \chi_1}{2(4-n_2)\chi_2} \eps^\frac{n_2-1}{2}
        \io \ue \nn\\
    &\le& \frac{[3(2-n_2)]^\frac{2-n_2}{2} n_2^\frac{n_2}{2} \chi_1 |\Omega|}{2(4-n_2)\chi_2} \nn\\
    & & + \frac{[3(3-n_2)]^\frac{3-n_2}{2} (n_2-1)^\frac{n_2-1}{2} a_2 \chi_1}{2(4-n_2)\chi_2}
    \cdot \Bigg\{ 1 + \bigg\{ \io \ue \bigg\}^7 \Bigg\}
  \eea
  for all $t>0$,
  where we have once more used that $\eps<1$, and that $n_2\ge 1$.\abs
  In summary, we only need to insert (\ref{8.99}) into (\ref{8.9}), combine the latter with (\ref{8.10}), (\ref{8.11})
  and (\ref{8.12}), and to particularly observe the dependence of the obtained estimate on $\chi_1$ and $\chi_2$ exclusively
  through the factor $1+\frac{\chi_1}{\chi_2}$, to infer from (\ref{8.8}) that indeed (\ref{8.1}) holds with an appropriate
  choice of $C$.
\qed
In order to prepare appropriate integration of (\ref{8.1}) for the purpose of deriving estimates essential for our
existence proof (Lemma \ref{lem16}), but also later in our qualitative analysis (Lemma \ref{lem19}),
let us note the following relationship between $\F$ and the nonnegative quantities
$\io \ue \ln(\ue+e)$ and $\io \ve \ln (\ve+e)$.
\begin{lem}\label{lem87}
  Let $D_i>0, a_i>0, \lambda_i>0$ and $n_i\in [1,2]$ for $i\in\{1,2\}$.
  Then there exists $C>0$ such that given any $\chi_1>0$ and $\chi_2>0$ and initial data such that
  (\ref{ie}) holds, with $\F$ as in (\ref{F}) we have
  \bea{87.1}
    & & \hspace*{-20mm}
    \io \ue(\cdot,t) \ln (\ue(\cdot,t)+e)
    + \frac{\chi_1}{\chi_2} \io \ve(\cdot,t) \ln (\ve(\cdot,t)+e) \nn\\
    &\le& \F(t) + C\cdot \Big(1+\frac{\chi_1}{\chi_2}\Big) \cdot
    \bigg\{ 1 + \io \ue(\cdot,t) + \io \ve(\cdot,t) \bigg\}
    \qquad \mbox{for all $t>0$ and } \eps\in (0,1).
  \eea
  Moreover,
  \be{87.2}
    \sup_{\eps\in (0,1)} \F(0) <\infty.
  \ee
\end{lem}
\proof
  Again relying on the inequality $\xi\ln\xi \ge -\frac{1}{e}$ for $\xi\in (0,e]$, for all $t>0$ and $\eps\in (0,1)$
  we can estimate
  \bas
    \io \ue \ln (\ue+e)
    &=& \int_{\{\ue\ge e\}} \ue \ln (\ue+e)
    + \int_{\{\ue<e\}} \ue \ln (\ue+e) \nn\\
    &\le& \int_{\{\ue\ge e\}} \ue \cdot (\ln \ue + \ln 2)
    + \ln (2e) \int_{\{\ue<e\}} \ue \\
    &\le& \io \ue\ln \ue
    + (\ln 2 + \ln (2e)) \io \ue + \frac{|\Omega|}{e},
  \eas
  and performing the same operations on $\ve$ we obtain that by (\ref{F}),
  \bas
    & & \hspace*{-10mm}
    \io \ue \ln (\ue+e)
    + \frac{\chi_1}{\chi_2} \io \ve \ln (\ve+e) \\
    &\le& \io \ue\ln \ue
    + \frac{\chi_1}{\chi_2} \io \ve \ln \ve \\
    & & + (\ln 2 + \ln (2e)) \cdot \bigg\{ \io \ue + \frac{\chi_1}{\chi_2} \io \ve \bigg\}
    + \Big(1+\frac{\chi_1}{\chi_2}\Big) \cdot \frac{|\Omega|}{e} \\
    &\le& \F(t)
    + (\ln 2 + \ln (2e)+1) \cdot \bigg\{ \io \ue + \frac{\chi_1}{\chi_2} \io \ve \bigg\}
    + \Big(1+\frac{\chi_1}{\chi_2}\Big) \cdot \frac{|\Omega|}{e}
    \quad \mbox{for all $t>0$ and } \eps\in (0,1),
  \eas
  which establishes (\ref{87.1}).\abs
  Apart from that, since $u_0$ and $v_0$ are positive, and since (\ref{ie}) requires that $u_{0\eps} \to u_0$ and
  $v_{0\eps} \to v_0$ in $L^\infty(\Omega)$ as $\eps\searrow 0$, it is evident from (\ref{F}) that
  $\F(0) \to \io u_0\ln u_0 - \io u_0 + \frac{\chi_1}{\chi_2} \io v_0\ln v_0 - \frac{\chi_1}{\chi_2} \io v_0$
  as $\eps\searrow 0$, which clearly entails (\ref{87.2}).
\qed
\mysection{Global existence. Proof of Theorem \ref{theo555}}\label{sect_conv}
Before proceeding, let us now sharpen our objective by specifying our solution concept in the following natural sense.
\begin{defi}\label{defi_weak}
  For $i\in \{1,2\}$, let $D_i>0, a_i>0, \lambda_i>0$ and $\chi_i>0$,
  and let $u_0\in L^1(\Omega)$ and $v_0\in L^1(\Omega)$ be nonnegative.
  Then a pair $(u,v)$ of nonnegative measurable functions defined in $\Omega\times (0,\infty)$ and satisfying
  \be{w1}
    \{ u^2,v^2,u_x,v_x,uv_x, vu_x\} \subset L^1_{loc}(\bom\times [0,\infty))
  \ee
  will be called a {\em global weak solution} of (\ref{0}) if
  \be{w2}
    - \int_0^\infty \io u\varphi_t - \io u_0 \varphi(\cdot,0)
    = - D_1 \int_0^\infty \io u_x\varphi_x
    + \chi_1 \int_0^\infty \io uv_x\varphi_x
    + \int_0^\infty \io u(\lambda_1-u+a_1 v) \varphi
  \ee
  and
  \be{w3}
    - \int_0^\infty \io v\varphi_t - \io v_0 \varphi(\cdot,0)
    = -D_2 \int_0^\infty \io v_x \varphi_x
    - \chi_2 \int_0^\infty \io vu_x\varphi_x
    + \int_0^\infty \io v(\lambda_2-v-a_1 u) \varphi
  \ee
  for all $\varphi\in C_0^\infty(\bom\times [0,\infty))$.
\end{defi}
Based on Lemma \ref{lem8}, in this section we shall derive estimates
for solutions to (\ref{0eps})-(\ref{0epsb}) which allow for an appropriate
subsequence extraction process, in the limit yielding a solution in
the above sense. As our approach parallels that in
\cite{taowin_double_crossdiff} in some parts, we may concentrate on
the essential and novel aspects here, beyond this restricting
ourselves to outlining the main steps.
\subsection{A first integration of (\ref{8.1}): Global regularity properties for fixed $\chi_i$}
Using the $L^1$ bounds provided by Lemma \ref{lem3}, for arbitrary but fixed cross-diffusion coefficients
we obtain the following from a first integration of (\ref{8.1}) in a straightforward manner, emphasizing that throughout this
section all appearing constants may depend on $\chi_1$ and $\chi_2$.
\begin{lem}\label{lem16}
  For $i\in \{1.2\}$, let $D_i>0, a_i>0, \lambda_i>0, \chi_i>0$ and $n_i\in [1,2]$,
  and assume (\ref{ie}).
  Then for all $T>0$ there exists $C(T)>0$ such that for any $\eps\in (0,1)$,
  \be{16.1}
    \io \ue(\cdot,t) \ln (\ue(\cdot,t)+e)
    + \io \ve(\cdot,t) \ln (\ve(\cdot,t)+e)
    \le C(T)
    \qquad \mbox{for all } t\in (0,T),
  \ee
  that with $\F$ as in (\ref{F}) we have
  \be{16.11}
    \F(t) \le C(T)
    \qquad \mbox{for all } t\in (0,T),
  \ee
  and that
  \be{16.2}
    \int_0^T \io \frac{\uex^2}{\ue}
    + \int_0^T \io \frac{\vex^2}{\ve}
    \le C(T)
  \ee
  as well as
  \be{16.3}
    \eps \int_0^T \io \ue^{n_1-1} \uexx^2
    + \eps\int_0^T \io \ve^{n_2-1} \vexx^2
    \le C(T).
  \ee
\end{lem}
\proof
  We employ Lemma \ref{lem3} to find $c_1>0$ such that
  \be{16.33}
    \io \ue + \io \ve \le c_1
    \qquad \mbox{for all $t>0$ and } \eps\in (0,1),
  \ee
  and observe that therefore Lemma \ref{lem8} implies that with some $c_2>0$, the functions $\F$ and $\D$ in (\ref{F})
  and (\ref{D}) satisfy
  \bas
    \frac{d}{dt} \F(t) + \D(t) \le c_2
    \qquad \mbox{for all $t>0$ and any } \eps\in (0,1)
  \eas
  and hence, upon integration,
  \be{16.4}
    \F(t) + \int_0^t \D(s) ds
    \le \F(0) + c_2 t
    \qquad \mbox{for all $t>0$ and each } \eps\in (0,1).
  \ee
  Thus, if in accordance with Lemma \ref{lem87} and (\ref{16.33}) we pick $c_3>0$ and $c_4>0$ large enough such that
  \bas
    \F(0) \le c_3
    \qquad \mbox{for all } \eps\in (0,1)
  \eas
  and
  \bas
    \io \ue \ln (\ue+e) + \frac{\chi_1}{\chi_2}\io \ve \ln (\ve+e) \le \F(t) + c_4
    \qquad \mbox{for all $t>0$ and } \eps\in (0,1),
  \eas
  then from (\ref{16.4}) we infer that
  \bas
    \io \ue \ln (\ue+e) + \frac{\chi_1}{\chi_2}\io \ve \ln (\ve+e) + \int_0^t \D(s) ds
    &\le& \F(t) + c_4 + \int_0^t \D(s) ds \\
    &\le& c_2 t + c_3 + c_4
    \qquad \mbox{for all $t>0$ and } \eps\in (0,1),
  \eas
  from which in view of (\ref{D}) both (\ref{16.1}) and (\ref{16.11}) as well as (\ref{16.2}) and (\ref{16.3}) directly
  follow.
\qed
Further consequences can be drawn by means of the following
interpolation inequality which can be found in \cite[Corollary
7.6]{taowin_double_crossdiff}.
\begin{lem}\label{lem13}
  There exists $C>0$ such that
  for all $\varphi\in W^{1,2}(\Omega)$ satisfying $\varphi>0$ in $\bom$,
  \bas
    \io \varphi^3 \ln (\varphi+e)
    \le C \cdot \bigg\{ \io \frac{\varphi_x^2}{\varphi} \bigg\} \cdot
    \bigg\{ \io \varphi \ln (\varphi+e) \bigg\}^2
    + C \cdot \bigg\{ \io \varphi \ln (\varphi+e) \bigg\}^3.
  \eas
\end{lem}
Inter alia by utilizing the latter, from Lemma \ref{lem16} we can derive three more spatio-temporal estimates,
the first and the third among which will yield precompactness of the cross-diffusive fluxes from (\ref{0eps})
in spatio-temporal $L^1$ spaces, and the second of which will play an important role in our asymptotic analysis by
implying positive lower bounds for the mass functionals in Lemma \ref{lem18}.
\begin{cor}\label{cor17}
  Let $D_i>0, a_i>0, \lambda_i>0, \chi_i>0$ and $n_i\in [1,2]$ for $i\in\{1,2\}$,
  and let (\ref{ie}) be fulfilled.
  Then for all $T>0$ there exists $C(T)>0$ such that
  \be{17.1}
    \int_0^T \io \ue^3 \ln (\ue +e)
    + \int_0^T \io \ve^3 \ln (\ve +e)
    \le C(T)
    \qquad \mbox{for all } \eps\in (0,1)
  \ee
  and
  \be{17.2}
    \int_0^T \|\ue(\cdot,t)\|_{L^\infty(\Omega)}^2 dt
    + \int_0^T \|\ve(\cdot,t)\|_{L^\infty(\Omega)}^2 dt
    \le C(T)
    \qquad \mbox{for all } \eps\in (0,1)
  \ee
  as well as
  \be{17.3}
    \int_0^T \io |\uex|^\frac{3}{2}
    + \int_0^T \io |\vex|^\frac{3}{2}
    \le C(T)
    \qquad \mbox{for all } \eps\in (0,1).
  \ee
\end{cor}
\proof
  Given $T>0$, we first apply Lemma \ref{lem16} to find $c_1(T)>0$ and $c_2(T)>0$ such that
  \be{17.4}
    \io \ue \ln (\ue+e) + \io \ve \ln (\ve+e)
    \le c_1(T)
    \qquad \mbox{for all $t\in (0,T)$ and } \eps\in (0,1)
  \ee
  and
  \be{17.5}
    \int_0^T \io \frac{\uex^2}{\ue}
    + \int_0^T \io \frac{\vex^2}{\ve} \le c_2(T)
    \qquad \mbox{for all } \eps\in (0,1).
  \ee
  As a consequence of Lemma \ref{lem13}, combining these inequality immediately yields (\ref{17.1}).\\
  Since by the H\"older inequality,
  \bas
    \int_0^T \io |\uex|^\frac{3}{2}
    \le \bigg\{ \int_0^T \io \frac{\uex^2}{\ue} \bigg\}^\frac{3}{4} \cdot
    \bigg\{ \int_0^T \io \ue^3 \bigg\}^\frac{1}{4}
    \qquad \mbox{for all } \eps\in (0,1),
  \eas
  using that $\ln (\ue+e)\ge 1$ and arguing similarly for $\ve$ we thereafter obtain (\ref{17.3}) from
  (\ref{17.1}) and (\ref{17.5}).\\
  Finally, as the Gagliardo-Nirenberg inequality provides $c_3>0$ fulfilling
  \bas
    \|\varphi\|_{L^\infty(\Omega)}^4
    \le c_3\|\varphi_x\|_{L^2(\Omega)}^2 \|\varphi\|_{L^2(\Omega)}^2
    + c_3\|\varphi\|_{L^2(\Omega)}^4
    \qquad \mbox{for all } \varphi\in W^{1,2}(\Omega),
  \eas
  once more combining (\ref{17.4}) with (\ref{17.5}) we infer that
  \bas
    \int_0^T \|\ue\|_{L^\infty(\Omega)}^2
    &\le& c_3 \int_0^T \|(\sqrt{\ue})_x\|_{L^2(\Omega)}^2 \|\sqrt{\ue}\|_{L^2(\Omega)}^2
    + c_3 \int_0^T  \|\sqrt{\ue}\|_{L^2(\Omega)}^4 \\
    &\le& \frac{c_1(T)c_2(T) c_3}{4} + c_1^2(T) c_3 \cdot T
    \qquad \mbox{for all } \eps\in (0,1)
  \eas
  and thereby conclude that also (\ref{17.2}) holds.
\qed
A last preparation for our limit procedure asserts $\eps$-independent estimates for the time derivatives
of $\ue$ and $\ve$ with respect to the norm in some suitably large spaces.
Beyond this, these bounds will later on once more be recalled in the course of our qualitative analysis
(cf.~Corollary \ref{cor24}).
\begin{lem}\label{lem52}
  Let $D_i>0, a_i>0, \lambda_i>0, \chi_i>0$ and $n_i\in [1,2]$ for $i\in\{1,2\}$,
  and assume (\ref{ie}).
  Then for all $T>0$ one can find $C(T)>0$ such that
  \bas
    \int_0^T \|\uet(\cdot,t)\|_{(W_0^{3,2}(\Omega))^\star} dt
    + \int_0^T \|\vet(\cdot,t)\|_{(W_0^{3,2}(\Omega))^\star} dt
    \le C(T)
    \qquad \mbox{for all } \eps\in (0,\epss).
  \eas
\end{lem}

\proof
  This can be obtained by minor modification of the reasoning in \cite[Lemma 3.4]{taowin_double_crossdiff}.
\qed
Collecting the above, we arrive at the following.
\begin{lem}\label{lem_conv}
  Let $D_i>0, a_i>0, \lambda_i>0, \chi_i>0$ and $n_i\in [1,2]$ for $i\in\{1,2\}$,
  and assume (\ref{ie}).
  Then one can find $(\eps_j)_{j\in \N} \subset (0,1)$ and nonnegative functions
  $u$ and $v$ defined in $\Omega\times (0,\infty)$ such that (\ref{555.1}) holds,
  that $\eps_j\searrow 0$ as $j\to\infty$, and that with some null set $N\subset (0,\infty)$ we have
  \begin{eqnarray}
    & & \ue\to u
    \quad \mbox{and} \quad
    \ve\to v
    \qquad \mbox{a.e.~in } \Omega\times (0,\infty), \label{53.2} \\
    & & \ue(\cdot,t)\to u(\cdot,t)
    \quad \mbox{and} \quad
    \ve(\cdot,t)\to v(\cdot,t)
    \qquad \mbox{a.e.~in $\Omega$ for all } t\in (0,\infty) \setminus N, \label{53.22} \\
    & & \ue\to u
    \quad \mbox{and} \quad
    \ve\to v
    \qquad \mbox{in } L^3_{loc}(\bom\times [0,\infty))
    \qquad \mbox{and} \label{53.3} \\
    & & \uex\wto u_x
    \quad \mbox{and} \quad
    \vex\wto v_x
    \qquad \mbox{in } L^\frac{3}{2}_{loc}(\bom\times [0,\infty)) \label{53.33}
  \eea
  as $\eps=\eps_j\searrow 0$. Moreover,     
  $(u,v)$ is a global weak solution of (\ref{0}) in the sense of Definition \ref{defi_weak}.
\end{lem}
\proof
  Form Corollary \ref{cor17} we know that
  \be{53.4}
    (\ue)_{\eps\in (0,1)}
    \mbox{ and }
    (\ve)_{\eps\in (0,1)}
    \quad \mbox{are bounded in }
    L^\frac{3}{2}_{loc}([0,\infty);W^{1,\frac{3}{2}}(\Omega)),
  \ee
  and that
  \be{53.44}
    \Big(\ue^3 \ln (\ue+e)\Big)_{\eps\in (0,1)}
    \mbox{ and }
    \Big(\ve^3 \ln (\ve+e)\Big)_{\eps\in (0,1)}
    \quad \mbox{are bounded in  }
    L^1_{loc}(\bom\times [0,\infty))
  \ee
  and that hence, by the de la Vall\'ee-Poussin theorem,
  \bas
    (\ue^3)_{\eps\in (0,1)}
    \mbox{ and }
    (\ve^3)_{\eps\in (0,1)}
    \quad \mbox{are uniformly integrable over } \Omega\times (0,T)
    \mbox{ for all } T>0.
  \eas
  As Lemma \ref{lem52} warrants that
  \bas
    (\uet)_{\eps\in (0,1)}
    \mbox{ and }
    (\vet)_{\eps\in (0,1)}
    \quad \mbox{are bounded in  }
    L^1_{loc}([0,\infty);(W_0^{3,2}(\Omega))^\star),
  \eas
  the existence of $(\eps_j)_{j\in\N} \subset (0,1)$, nonnegative functions $u$ and $v$
  on $\Omega\times (0,\infty)$ and a null set $N\subset (0,\infty)$, as well as the
  convergence properties (\ref{53.2})-(\ref{53.33}), result from straightforward arguments involving the
  Aubin-Lions lemma (\cite{simon}), the Vitali convergence theorem and the Dunford-Pettis theorem.
  The verification of (\ref{w1})-(\ref{w3}) can thereafter be accompliches in much the same manner
  as demonstrated in \cite[Lemma 4.1]{taowin_double_crossdiff}, so that we may refrain from giving details here.
\qed
This in fact already contains our main result on global existence in (\ref{0}):\abs
\proofc of Theorem \ref{theo555}.\quad
  In view of Lemma \ref{lem_conv}, all statements are now obvious.
\qed
%
%
%
%
%
%
%
%
%
%
%
%
%
%
%
%
%
%
%
\mysection{Eventual bounds I. A second integration of (\ref{8.1})}
Next addressing the large time behavior of the solutions gained above in the case when $\chi_1$ and $\chi_2$ are
suitably small, since in view of Theorem \ref{theo33} and Theorem \ref{theo333} we intend to admit initial data of
arbitrary size our first step will consist in making sure that any such solution can eventually be controlled
by quantities independent of the initial data, similarly to the second conclusion from Lemma \ref{lem3} but in
more useful topological frameworks.\abs
For this purpose, we shall once more return to the quasi-entropy inequality from Lemma \ref{lem8}, but
now with the ambition to make more efficient use of the dissipation rate appearing therein.
Here an obstacle toward straightforward approches seems to be linked to the observation that
according to (\ref{F}), for each fixed $\eps\in (0,1)$ the functional $\F$ contains summands proportional
to the integrals $\io \frac{1}{\ue^{3-n_1}}$ and $\io \frac{1}{\ve^{3-n_2}}$
for which our previous estimates do not provide any meaningful information.\abs
An elementary but important interpolation lemma, formulated in a way general enough to allow
for two further applications later in Lemma \ref{lem27} and Lemma \ref{lem277},
will be of decisive support in overcoming this difficulty:
\begin{lem}\label{lem14}
  Let $\varphi\in C^1(\bom)$ be such that $\varphi>0$ in $\bom$.\abs
  i) \ For all $p>0$ an $q>0$,
  \be{14.1}
    \io \varphi^{-p}
    \le q^\frac{2p}{q} |\Omega|^\frac{p+q}{q} \cdot
    \bigg\{ \io \varphi^{-q-2} \varphi_x^2 \bigg\}^\frac{p}{q}
    + 2^\frac{2p}{q} |\Omega|^{p+1} \cdot
    \bigg\{ \io \varphi \bigg\}^{-p}.
  \ee
  ii) \ The inequality
  \be{14.2}
    - \io \ln \varphi
    \le |\Omega|^\frac{3}{2} \cdot \bigg\{ \io \frac{\varphi_x^2}{\varphi^2} \bigg\}^\frac{1}{2}
    - |\Omega| \cdot \ln \bigg\{ \io \varphi \bigg\}
    + |\Omega| \cdot \ln |\Omega|
  \ee
  holds.
\end{lem}
\proof
  i) \ We fix $x_0\in\bom$ such that $\varphi(x_0) \ge \frac{1}{|\Omega|} \io \varphi$ and then use the Cauchy-Schwarz
 inequality to see that for all $x\in\Omega$,
  \bas
    \varphi^{-\frac{q}{2}}(x)
    = \varphi^{-\frac{q}{2}}(x_0)
    + \int_{x_0}^x (\varphi^{-\frac{q}{2}})_x(y) dy
    \le |\Omega|^\frac{q}{2} \cdot \bigg\{ \io \varphi \bigg\}^{-\frac{q}{2}}
    + |\Omega|^\frac{1}{2} \cdot \bigg\{ \io (\varphi^{-\frac{q}{2}})_x^2 \bigg\}^\frac{1}{2}.
  \eas
  Thanks to Young's inequality, this implies that
  \bas
    \varphi^{-p}(x)
    &=& \Big\{ \varphi^{-\frac{q}{2}}(x) \Big\}^\frac{2p}{q} \\
    &\le& 2^\frac{2p}{q} |\Omega|^p \cdot \bigg\{ \io \varphi\bigg\}^{-p}
    + 2^\frac{2p}{q} |\Omega|^\frac{p}{q} \cdot \bigg\{ \io (\varphi^{-\frac{q}{2}})_x^2 \bigg\}^\frac{p}{q} \\
    &=& 2^\frac{2p}{q} |\Omega|^p \cdot \bigg\{ \io \varphi\bigg\}^{-p}
    + q^\frac{2p}{q} |\Omega|^\frac{p}{q} \cdot \bigg\{ \io \varphi^{-q-2} \varphi_x^2 \bigg\}^\frac{p}{q}
    \qquad \mbox{for all } x\in\Omega
  \eas
  and thus establishes (\ref{14.1}) upon integration.\abs
  ii) \ Likewise,
  \bas
    -\ln\varphi(x)
    = - \ln\varphi(x_0)
    - \int_{x_0}^x \frac{\varphi_x(y)}{\varphi(y)} dy
    \le - \ln \bigg\{ \io \varphi\bigg\}
    + \ln |\Omega|
    + |\Omega|^\frac{1}{2} \cdot \bigg\{ \io \frac{\varphi_x^2}{\varphi^2} \bigg\}^\frac{1}{2},
  \eas
  which after integration results in (\ref{14.2}).
\qed
Indeed, an application of the first part thereof enables us to create an absorptive linear summand in (\ref{8.1})
at the expense of additional expressions on its right-hand side, which however remain under control as long as
the mass functionals $\io \ue$ and $\io \ve$ remain uniformly positive:
\begin{lem}\label{lem89}
  Given $D_i>0, a_i>0, \lambda_i>0, \chi_i>0$ and $n_i\in [1,2]$ for $i\in\{1,2\}$,
  one can find $C>0$ such that if $\chi_1>0$ and $\chi_2>0$ and (\ref{ie}) holds,
  then for all $\eps\in (0,1)$
  the functions $\F$ and $\D$ from (\ref{F}) and (\ref{D}) have the property that
  \bea{89.1}
    & & \hspace*{-20mm}
    \frac{d}{dt} \F(t) + \frac{1}{C} \F(t) + \frac{1}{2} \D(t) \nn\\
    &\le& C \cdot \Big(1+\frac{\chi_1}{\chi_2}\Big) \cdot \Bigg\{
    1 + \bigg\{ \io \ue(\cdot,t) \bigg\}^7 + \bigg\{ \io \ve(\cdot,t) \bigg\}^7 \nn\\
    & & \hspace*{28mm}
    + \, \eps \cdot \bigg\{ \io \ue(\cdot,t) \bigg\}^{-(3-n_1)}
    + \eps \cdot \bigg\{ \io \ve(\cdot,t) \bigg\}^{-(3-n_2)} \Bigg\}
  \eea
  for all $t>0$.
\end{lem}
\proof
  In view of Lemma \ref{lem8} and Young's inequality, it is clearly sufficient to make sure that given
  $D_1>0, D_2>0, a_1>0, a_2>0, \lambda_1>0$ and $\lambda_2>0$, one can find $c_1>0$ such that whenever
  $\chi_1>0$ and $\chi_2>0$ and (\ref{ie}) holds, then for any $\eps\in (0,1)$ we have
  \bea{89.2}
    \F(t)
    &\le& c_1 \D(t)
    + c_1 \cdot \Big(1+\frac{\chi_1}{\chi_2}\Big) \cdot
    \Bigg\{  1+ \bigg\{ \io \ue\bigg\}^6 + \bigg\{ \io \ve\bigg\}^6 \nn\\
    & & \hspace*{45mm}
    + \, \eps \cdot \bigg\{ \io \ue \bigg\}^{-(3-n_1)} + \eps\cdot \bigg\{ \io \ve\bigg\}^{-(3-n_2)} \Bigg\}
  \eea
  for all $t>0$.
  To achieve this, we note  that e.g.~by once more combining the Gagliardo-Nirenberg inequality with YoungÄs inequality
  in a straightforward manner we can find positive constants $c_2$ and $c_3$ such that for all $\varphi\in C^1(\bom)$
  fulfilling $\varphi>0$ in $\bom$ we have
  \bas
    \io \varphi\ln \varphi
    \le \io \varphi^2
    \le c_2 \|(\sqrt{\varphi})_x\|_{L^2(\Omega)} \|\sqrt{\varphi}\|_{L^2(\Omega)}^3
    + c_2 \|\sqrt{\varphi}\|_{L^2(\Omega)}^4
    \le \io \frac{\varphi_x^2}{\varphi}
    + c_3 \cdot \bigg\{ \io \varphi \bigg\}^6 + c_3.
  \eas
  Therefore,
  \bea{89.3}
    \F(t)
    &\le& \io \frac{\uex^2}{\ue}
    + c_3 \cdot \bigg\{ \io \ue \bigg\}^6 + c_3 \nn\\
    & & + \frac{\chi_1}{\chi_2} \io \frac{\vex^2}{\ve}
    + \frac{c_3\chi_1}{\chi_2} \cdot \bigg\{ \io \ve \bigg\}^6 + \frac{c_3 \chi_1}{\chi_2} \nn\\
    & & + \frac{\eps}{(3-n_1)(3-n_2)} \io \frac{1}{\ue^{3-n_1}}
    + \frac{\chi_1\eps}{(3-n_2)(4-n_2)\chi_2} \io \frac{1}{\ve^{3-n_2}}
    \qquad \mbox{for all } t>0,
  \eea
  where the two rightmost summands can be estimated by means of Lemma \ref{lem14} i), which when applied to $p=q=3-n_i$,
  $i\in \{1,2\}$, says that for all $t>0$,
  \bas
    \frac{\eps}{(3-n_1)(4-n_1)} \io \frac{1}{\ue^{3-n_1}}
    &\le& \frac{(3-n_1) |\Omega|^2}{4-n_1} \cdot \eps \io \frac{\uex^2}{\ue^{5-n_1}}
    + \frac{4|\Omega|^{4-n_1}}{(3-n_1)(4-n_1)} \cdot \eps \cdot \bigg\{ \io \ue \bigg\}^{-(3-n_1)}
  \eas
  and
  \bas
    \frac{\chi_1 \eps}{(3-n_2)(4-n_2)\chi_2} \io \frac{1}{\ve^{3-n_2}}
    \le \frac{(3-n_2)\chi_1 |\Omega|^2}{(4-n_2)\chi_2} \cdot \eps\io \frac{\vex^2}{\ve^{5-n_2}}
    + \frac{4\chi_1 |\Omega|^{4-n_2}}{(3-n_2)(4-n_2)\chi_2} \cdot \eps \cdot \bigg\{ \io \ve \bigg\}^{-(3-n_2)}.
  \eas
  In conjunction with (\ref{89.3}) and the definition (\ref{D}) of $\D$, these inequalities readily establish
  (\ref{89.1}) with some suitably large $C>0$.
\qed
As a natural next task, we are thus concerned with positive lower bounds for $\io\ue$ and $\io \ve$.
Thanks to the second estimate provided by Corollary \ref{cor17}, however, at least within time intervals of finite
length these can be obtained by simply returning to (\ref{0eps})-(\ref{0epsb}):
\begin{lem}\label{lem18}
  Let $D_i>0, a_i>0, \lambda_i>0, \chi_i>0$ and $n_i\in [1,2]$ for $i\in\{1,2\}$,
  and suppose that (\ref{ie}) holds.
  Then for all $T>0$ there exists $C(T)>0$ such that whenever $\eps\in (0,1)$,
  \be{18.1}
    \io \ue(\cdot,t) \ge C(T)
    \quad \mbox{and} \quad
    \io \ve(\cdot,t) \ge C(T)
    \qquad \mbox{for all } t\in (0,T).
  \ee
\end{lem}
\proof
  For fixed $T>0$, using Corollary \ref{cor17} and the Cauchy-Schwarz inequality we can fix $c_1(T)>0$ such that
  \be{18.2}
    \int_0^T \|\ue\|_{L^\infty(\Omega)}
    + \int_0^T \|\ve\|_{L^\infty(\Omega)} \le c_1(T)
    \qquad \mbox{for all } \eps\in (0,1).
  \ee
  Thus, if going back to (\ref{0eps})-(\ref{0epsb}) we estimate
  \bas
    \frac{d}{dt} \io \ue
    &=& \io \frac{3\ue^3}{3\ue^2+\eps} (\lambda_1-\ue+a_1\ve) \\
    &\ge& - \io \ue^2 \\
    &\ge& - \|\ue\|_{L^\infty(\Omega)} \io \ue
    \qquad \mbox{for all $t>0$ and } \eps\in (0,1),
  \eas
  then upon integration we infer that thanks to (\ref{18.2}),
  \bas
    \io \ue(\cdot,t)
    &\ge& \bigg\{ \io u_{0\eps} \bigg\} \cdot e^{-\int_0^t \|\ue(\cdot,s)\|_{L^\infty(\Omega)} ds} \\
    &\ge& c_2 e^{-c_1(T)}
    \qquad \mbox{for all $t\in (0,T)$ and } \eps\in (0,1),
  \eas
  with $c_2:=\inf_{\eps\in (0,1)} \io u_{0\eps}$ being positive due to (\ref{ie}) and the positivity of $u_0$.
  Likewise, from the inequality
  \bas
    \frac{d}{dt} \io \ve
    &=& \io \frac{3\ve^3}{3\ve^2+\eps} (\lambda_2 - \ve -a_2\ue) \\
    &\ge& - \io \ve^2 - a_2 \io \ue\ve \\
    &\ge& - \|\ve\|_{L^\infty(\Omega)} \io \ve
    - a_2\|\ue\|_{L^\infty(\Omega)} \io \ve
    \qquad \mbox{for all $t>0$ and } \eps\in (0,1),
  \eas
  through (\ref{18.2}) we obtain that
  \bas
    \io \ve(\cdot,t)
    &\ge& \bigg\{ \io v_{0\eps} \bigg\} \cdot e^{-\int_0^t \|\ve(\cdot,s)\|_{L^\infty(\Omega)} ds}
    \cdot e^{-a_2 \int_0^t \|\ue(\cdot,s)\|_{l^\infty(\Omega)} ds} \\
    &\ge& c_3 e^{-c_1(T)} \cdot e^{-a_2c_1(T)}
    \qquad \mbox{for all $t\in (0,T)$ and any } \eps\in (0,1)
  \eas
  with $c_3:=\inf_{\eps\in (0,1)} \io v_{0\eps}>0$.
\qed
Making use of the fact that the singular expressions on the right-hand side of (\ref{89.1}) contain the small factor $\eps$,
on integration thereof we infer the following data-independent eventual regularity information as the main outcome of
this section.
\begin{lem}\label{lem19}
  For $i\in\{1,2\}$, let $D_i>0, a_i>0, \lambda_i>0$ and $n_i\in [1,2]$.
  Then there exists $K>0$ with the following property:
  If $\chi_1>0$ and $\chi_2>0$ and (\ref{ie}) holds, one can find
  $T_0=T_0(\chi_1,\chi_2,u_0,v_0)>0$ such that for each $T>T_0$ there exists
  $\eps_0=\eps_0(T,\chi_1,\chi_2,u_0,v_0)\in (0,1)$ such that
  whenever $\eps\in (0,\eps_0)$,
  \be{19.1}
    \io \ue(\cdot,t) \ln (\ue(\cdot,t)+e)
    + \frac{\chi_1}{\chi_2} \io \ve(\cdot,t) \ln (\ve(\cdot,t)+e)
    \le K \cdot \Big(1+\frac{\chi_1}{\chi_2}\Big)
    \qquad \mbox{for all } t\in (T_0,T)
  \ee
  and
  \be{19.2}
    \int_t^{t+1} \io \frac{\uex^2}{\ue}
    + \frac{\chi_1}{\chi_2} \int_t^{t+1} \io \frac{\vex^2}{\ve}
    \le K \cdot \Big(1+\frac{\chi_1}{\chi_2}\Big)
    \qquad \mbox{for all } t\in (T_0,T).
  \ee
\end{lem}
\proof
  By means of Lemma \ref{lem89} and Lemma \ref{lem87}, given any $D_1>0, D_2>0, a_1>0, a_2>0,\lambda_1>0$ and $\lambda_2>0$
  we can find $c_1>0, c_2>0$ and $c_3>0$ with the property that whenever $\chi_1>0$, $\chi_2>0$ and (\ref{ie})
  holds, for $\F$ and $\D$ taken from (\ref{F}) and (\ref{D}) we have
  \bea{19.3}
    & & \hspace*{-16mm}
    \frac{d}{dt} \F(t)
    + c_1 \F(t)
    + \frac{1}{2} \D(t) \nn\\
    & & \hspace*{-8mm}
    \le c_2 \cdot \Big(1+\frac{\chi_1}{\chi_2}\Big) \cdot \Bigg\{
    1 + \bigg\{ \io \ue \bigg\}^7
    + \bigg\{ \io \ve\bigg\}^7
    + \eps \cdot \bigg\{ \io \ue \bigg\}^{-(3-n_1)}
    + \eps \cdot \bigg\{ \io \ve \bigg\}^{-(3-n_2)} \Bigg\}
  \eea
  and
  \be{19.4}
    \io \ue\ln (\ue+e)
    + \frac{\chi_1}{\chi_2} \io \ve \ln (\ve+e)
    \le \F(t) + c_3\cdot \Big(1+\frac{\chi_1}{\chi_2}\Big) \cdot
    \bigg\{ 1 + \io \ue + \io \ve \bigg\}
  \ee
  for all $t>0$ and any $\eps\in (0,1)$.
  Moreover, combining the inequalities (\ref{3.1}) and (\ref{3.2}) provided by Lemma \ref{lem3}
  we see that if (\ref{ie}) is valid, then there exists $T_1=T_1(u_0,v_0)$ such that
  \be{19.5}
    \io \ue + \io \ve
    \le c_4:=m_\infty+1
    \qquad \mbox{for all $t>T_1$ and each } \eps\in (0,1),
  \ee
  where $m_\infty>0$ is as defined in (\ref{3.2}).
  Keeping this value of $T_1$ fixed henceforth, we now assume that $D_1,D_2,a_1,a_2,\lambda_1,\lambda_2, \chi_1$
  and $\chi_2$ are given positive constants, and that (\ref{ie}) holds.
  Then an application of Lemma \ref{lem16} yields $c_5=c_5(\chi_1,\chi_2,u_0,v_0)>0$ such that
  the correspondingly defined function $\F$ in (\ref{F}) satisfies
  \be{19.6}
    \F(T_1) \le c_5
    \qquad \mbox{for all } \eps\in (0,1),
  \ee
  and we thereupon let $T_2=T_2(\chi_1,\chi_2,u_0,v_0)>T_1$ be large enough such that
  \be{19.77}
    c_5 e^{-(T_2-T_1)} \le 1.
  \ee
  Next, fixing any $T>T_2+1$ we infer from Lemma \ref{lem18} that
  \be{19.7}
    \io \ue \ge c_6
    \quad \mbox{and} \quad
    \io \ve \ge c_6
    \qquad \mbox{for all $t\in (0,T)$ and } \eps\in (0,1)
  \ee
  with some suitably small $c_6=c_6(T,\chi_1,\chi_2,u_0,v_0)>0$, whence it becomes possible to
  pick $\eps_0\equiv \eps_0(T,\chi_1,\chi_2,u_0,v_0)\in (0,1)$ small enough such that
  \bas
    c_6^{-(3-n_1)} \eps_0
    + c_6^{-(3-n_2)} \eps_0
    \le 1.
  \eas
  Then on inserting (\ref{19.5}) and (\ref{19.7}) into (\ref{19.3}) we obtain that
  \be{19.8}
    \frac{d}{dt} \F(t)
    + c_1 \F(t) + \frac{1}{2} \D(t)
    \le c_7 \cdot \Big(1+\frac{\chi_1}{\chi_2}\Big)
    \qquad \mbox{for all $t\in (T_1,T)$ and } \eps\in (0,\eps_0),
  \ee
  where we underline that besides $c_1$, also $c_7:=c_2 \cdot (2+2c_4^7)$ is independent of both $\chi_1$ and $\chi_2$
  as well as $u_0$ and $v_0$.
  After integration using (\ref{19.6}) and (\ref{19.77}), from (\ref{19.8}) we firstly infer that
  \bea{19.9}
    \F(t)
    &\le& \F(T_1) \cdot e^{-c_1(t-T_1)}
    + c_7 \cdot \Big(1+\frac{\chi_1}{\chi_2}\Big) \cdot \int_{T_1}^t e^{-c_1(t-s)} ds \nn\\
    &\le& c_5 e^{-c_1(t-T_1)}
    + c_7 \cdot \Big(1+\frac{\chi_1}{\chi_2}\Big) \cdot \frac{1-e^{-c_1(t-T_1)}}{c_1} \nn\\
    &\le& c_5 e^{-c_1(t-T_1)}
    + \frac{c_7}{c_1} \cdot \Big(1+\frac{\chi_1}{\chi_2}\Big) \nn\\
    &\le& 1 + \frac{c_7}{c_1} \cdot \Big(1+\frac{\chi_1}{\chi_2}\Big)
    \qquad \mbox{for all $t\in T_2,T)$ and } \eps\in (0,\eps_0),
  \eea
  and that thus in view of (\ref{19.4}) and (\ref{19.5}), in particular,
  \bea{19.10}
    \io \ue\ln (\ue+e)
    + \frac{\chi_1}{\chi_2} \io \ve\ln (\ve+e)
    &\le &1 + \frac{c_7}{c_1} \cdot \Big(1+\frac{\chi_1}{\chi_2}\Big)+ c_3 \cdot \Big(1+\frac{\chi_1}{\chi_2}\Big) \cdot
    (1+c_4)\nn\\
    & & \qquad \mbox{for all $t\in (T_2,T)$ and } \eps\in (0,\eps_0).
  \eea
  Thanks to (\ref{19.9}), (\ref{19.8}) secondly entails that for all $t\in (T_2,T)$ and each $\eps\in (0,\eps_0)$,
  \bas
    \frac{1}{2} \int_t^{t+1} \D(s) ds
    &\le& \F(t) - \F(t+1)
    - c_1\int_t^{t+1} \F(s) ds
    + c_7\cdot\Big(1+\frac{\chi_1}{\chi_2}\Big) \\
    &\le& 1 + \frac{c_7}{c_1} \cdot \Big(1+\frac{\chi_1}{\chi_2}\Big)
    + c_3(c_4+1) \Big(1+\frac{\chi_1}{\chi_2}\Big)
    + c_1 c_3 (c_4+1) \Big(1+\frac{\chi_1}{\chi_2}\Big)
    + c_7 \Big(1+\frac{\chi_1}{\chi_2}\Big),
  \eas
  because
  \bas
    \F(t) \ge - c_3 \Big(1+\frac{\chi_1}{\chi_2}\Big) (1+c_4)
    \qquad \mbox{for all $t>T_1$ and any } \eps\in (0,1)
  \eas
  according to (\ref{19.4}) and (\ref{19.5}).
  Due to the definition (\ref{D}) of $\D$, the latter implies (\ref{19.2}), whereas (\ref{19.1}) results from
  (\ref{19.10}).
\qed
\mysection{Eventual bounds II. The conditional quasi-entropy $\io \uex^2 + \gamma \io \vex^2$ }\label{sect_evbounds}
Now approaching the core of our asymptotic analysis, we next intend to improve the information on eventual regularity
provided by Lemma \ref{lem19}.  
To this end, in the key Lemma \ref{lem62} we shall  
study the time evolution of $\io \uex^2 + \gamma \io \vex^2$,
and more precisely we shall see there that under suitable smallness assumptions on $\chi_1$ and $\chi_2$,
for carefully chosen values of the parameter $\gamma>0$ depending on the ratio $\frac{\chi_1}{\chi_2}$ in a subtle manner,
this functional enjoys an entropy-like property {\em as long as its values remain small}.
Using Lemma \ref{lem19} to achieve the latter at least at some large initial time, we will thereby be able to indeed conclude
eventual boundedness of both solution components in $W^{1,2}(\Omega)$.\abs
In estimating the respective cross-diffusion terms appearing in the corresponding testing procedure adequately,
we will make use of the following elementary inequalities.
\begin{lem}\label{lem61}
  Let $n\in [0,\frac{7}{2}]$, $\eps>0$ and
  \be{61.1}
    h_\eps(s):=\frac{s^{5-n}}{s^{4-n}+\eps},
    \qquad s\ge 0.
  \ee
  Then
  \be{61.2}
    0 \le h_\eps(s) \le s
    \quad \mbox{and} \quad
    0 \le h_\eps'(s) \le 5-n
    \qquad \mbox{for all } s\ge 0,
  \ee
  and moreover
  \be{61.3}
    |h_\eps''(s)| \le 2^{-\frac{7-2n}{2(4-n)}} \cdot (4-n)(5-n) \eps^\frac{1}{2(4-n)} s^{-\frac{3}{2}}
    \qquad \mbox{for all } s>0.
  \ee
\end{lem}
\proof
  Computing
  \bas
    h_\eps'(s)=\frac{s^{8-2n} + (5-n) \eps s^{4-n}}{(s^{4-n}+\eps)^2},
    \qquad s\ge 0,
  \eas
  and
  \be{61.4}
    h_\eps''(s) = \frac{-(3-n)(4-n) \eps s^{7-2n} + (4-n)(5-n) \eps^2 s^{3-n}}{(s^{4-n}+\eps)^3},
    \qquad s>0,
  \ee
  we immediately see that since $1\le 5-n$,
  \bas
    0 \le h_\eps'(s)
    \le \frac{(5-n)s^{8-2n} + (5-n)\eps s^{4-n}}{(s^{4-n}+\eps)^2}
    = (5-n) \frac{s^{4-n}}{s^{4-n}+\eps} \le 5-n
    \qquad \mbox{for all } s\ge 0,
  \eas
  whence (\ref{61.2}) becomes evident.\\
  Furthermore, since $|3-n| \le 5-n$, from (\ref{61.4}) it follows that
  \be{61.5}
    |h_\eps''(s)|
    \le (4-n)(5-n) \cdot \frac{\eps s^{7-2n} + \eps^2 e^{3-n}}{(s^{4-n}+\eps)^3}
    = (4-n)(5-n) I_\eps(s)
    \qquad \mbox{for all } s>0,
  \ee
  where clearly
  \bas
    I_\eps(s):=\frac{\eps s^{3-n}}{(s^{4-n}+\eps)^2}
    \le \frac{\eps s^{3-n}}{s^{8-2n}} = \eps s^{n-5}
    \qquad \mbox{for all } s>0,
  \eas
  and where on the other hand by Young's inequality,
  \bas
    I_\eps(s) = \frac{1}{s} \cdot \frac{\eps s^{4-n}}{(s^{4-n}+\eps)^2}
    \le \frac{1}{s} \cdot \frac{\frac{1}{2} s^{8-2n}+ \frac{1}{2} \eps^2}{(s^{4-n}+\eps)^2}
    \le \frac{1}{2s}
    \qquad \mbox{for all } s>0.
  \eas
  An interpolation therefore shows that since $n\le \frac{7}{2}$,
  \bas
    I_\eps(s)
    &=& I_\eps^\frac{1}{2(4-n)}(s) \cdot I_\eps^\frac{7-2n}{8-2n}(s) \nn\\
    &\le& \Big\{ \eps^\frac{1}{2(4-n)} s^\frac{n-5}{8-2n} \Big\} \cdot
    \Big\{ (2s)^{-\frac{7-2n}{8-2n}} \Big\}
    \qquad \mbox{for all } s>0,
  \eas
  and that thus (\ref{61.3}) is a consequence of (\ref{61.5}).
\qed
Now the presumably most crucial step toward our derivation both of Theorem \ref{theo33} and of Theorem \ref{theo333} will
be accomplished in the course of the following quite delicate argument.
Here a considerable technical intricacy originates
from the ambition to unambiguously identify smallness assumptions on $\chi_1$ and $\chi_2$, independently of $u_0$ and $v_0$,
and on $\eps$, the latter possibly data-dependent, that ensure the following global absorption property considerably going
beyond that from Lemma \ref{lem3}.
\begin{lem}\label{lem62}
  Let $D_i>0, a_i>0, \lambda_i>0$ and $n_i\in [1,2]$ for $i\in\{1,2\}$.
  Then there exist $\chis>0$ and $C>0$ with the property that whenever $\chi_1\in (0,\chis)$ and $\chi_2\in (0,\chis)$
  and (\ref{ie}) holds, one can find $\ts=\ts(\chi_1,\chi_2,u_0,v_0)>0$ and
  $\epss=\epss(\chi_1,\chi_2,u_0,v_0)\in (0,1)$ such that
  \be{62.1}
    \io \uex^2(\cdot,t) + \io \vex^2(\cdot,t) \le C
    \qquad \mbox{for all $t>\ts$ and each } \eps\in (0,\epss)
  \ee
  and
  \be{62.2}
    \|\ue(\cdot,t)\|_{L^\infty(\Omega)} +
    \|\ve(\cdot,t)\|_{L^\infty(\Omega)}
    \le C
    \qquad \mbox{for all $t>\ts$ and } \eps\in (0,\epss)
  \ee
  as well as
  \be{62.02}
    \int_t^{t+1} \io \uexx^2 + \int_t^{t+1} \io \vexx^2
    \le C
    \qquad \mbox{for all $t>\ts$ and any } \eps\in (0,\epss).
  \ee
\end{lem}
\proof
  In order to prepare our selection of $\chis$, taking $m_\infty>0$ from Lemma \ref{lem3} let us set $m:=m_\infty+1$ and
  repeatedly make use of the Gagliardo-Nirenberg inequality to find positive constants $c_1, c_2, c_3, c_4$ and $c_5$
  such that for any $\varphi\in W^{2,2}(\Omega)$ fulfilling $\|\varphi\|_{L^1(\Omega)} \le m$ we have
  \be{62.3}
    \|\varphi\|_{L^\infty(\Omega)}^2 \le c_1 \cdot \Big\{ \|\varphi_{xx}\|_{L^2(\Omega)}^\frac{4}{5}+1\Big\}
  \ee
  and
  \be{62.4}
    \|\varphi\|_{L^2(\Omega)}^2 \le c_2 \cdot \Big\{ \|\varphi_{xx}\|_{L^2(\Omega)}^\frac{2}{5}+1\Big\}
  \ee
  and
  \be{62.5}
    \|\varphi\|_{L^\infty(\Omega)} \le c_3 \cdot \Big\{ \|\varphi_x\|_{L^2(\Omega)}^\frac{2}{3}+1\Big\}
  \ee
  as well as
  \be{62.6}
    \|\varphi_x\|_{L^2(\Omega)}^2 \le c_4 \cdot \Big\{ \|\varphi_{xx}\|_{L^2(\Omega)}^\frac{6}{5}+1\Big\}
  \ee
  and
  \be{62.7}
    \|\varphi_x\|_{L^4(\Omega)}^\frac{10}{9} \le c_5\cdot \Big\{ \|\varphi_{xx}\|_{L^2(\Omega)}^\frac{7}{9}+1\Big\}.
  \ee
  By the same token, we can fix $c_6>0$ and $c_7>0$ such that
  \be{62.8}
    \|\varphi\|_{L^4(\Omega)} \le c_6\|\varphi_x\|_{L^2(\Omega)}^\frac{1}{4} \|\varphi\|_{L^2(\Omega)}^\frac{3}{4}
    \qquad \mbox{for all } \varphi\in W_0^{1,2}(\Omega)
  \ee
  and
  \be{62.9}
    \|\varphi\|_{L^\infty(\Omega)}^4 \le c_7 \cdot \Big\{ \|\varphi_x\|_{L^2(\Omega)}^2 + 1 \Big\}
    \qquad \mbox{for all $\varphi\in W^{1,2}(\Omega)$ satisfying } \|\varphi\|_{L^2(\Omega)}^2 \le m.
  \ee
  We next employ Young's inequality to choose positive numbers $c_8, c_9, c_{10}$ and $c_{11}$ with the properties that
  \be{62.10}
    \Big( \frac{a_1^2}{2D_1} + \frac{a_2^2}{2D_2}\Big) \cdot c_1 c_2 \xi\eta
    \le \frac{D_1}{8\sqrt{2}} \xi^\frac{5}{2} + c_8 \eta^\frac{5}{3}
    \qquad \mbox{for all $\xi\ge 0$ and } \eta\ge 0
  \ee
  and
  \be{62.11}
    2^\frac{2}{3} c_8 \xi\eta \le \frac{D_2}{16} \xi^3 + c_9 \eta^\frac{3}{2}
    \qquad \mbox{for all $\xi\ge 0$ and } \eta\ge 0
  \ee
  as well as
  \be{62.12}
    3\lambda_1 c_4 \xi^\frac{6}{5} \le \frac{D_1}{8} \xi^2 + c_{10}
    \qquad \mbox{for all } \xi\ge 0
  \ee
  and
  \be{62.13}
    3\lambda_2 c_4 \xi^\frac{6}{5} \le \frac{D_2}{8} \xi^2 + c_{11}
    \qquad \mbox{for all } \xi\ge 0,
  \ee
  and thereupon abbreviate
  \be{62.14}
    c_{12} := \frac{D_1+D_2}{4} + c_9
  \ee
  as well as
  \be{62.15}
    c_{13} := c_{10} + 3\lambda_1 c_4
    \qquad \mbox{and} \qquad
    c_{14} := c_{11} + 3\lambda_2 c_4.
  \ee
  Furthermore introducing
  \be{62.16}
    c_{15} := \max_{i\in\{1,2\}} \Big\{ c_3 + (5-n_i) c_5 c_6^\frac{8}{9} \Big\}
    \qquad \mbox{and} \qquad
    c_{16} := \max_{i\in\{1,2\}} 2^{-\frac{7-2n_i}{2(4-n_i)}-1} \cdot (4-n_i)(5-n_i),
  \ee
  we let
  \be{62.17}
    c_{17}:=c_{12}+c_{13}+c_{14}+ \frac{D_1+D_2}{32}+2
    \qquad \mbox{and} \qquad
    c_{18} := c_6^4 c_{16}
  \ee
  as well as
  \be{62.18}
    c_{19} := \frac{\min\{D_1,D_2\}}{32\cdot (3c_4)^\frac{5}{3}}
    \qquad \mbox{and} \qquad
    c_{20}:=c_{17} + \frac{2^\frac{2}{3}\cdot\min\{D_1,D_2\}}{16},
  \ee
  and taking $K>0$ from Lemma \ref{lem19} we write
  \be{62.19}
    c_{21}:=2K \cdot \sqrt{\frac{c_7 K}{2} + c_7}, \quad
    c_{22} := \max \Big\{ 1 \, , \, 4c_{21} \, , \, 2\cdot \Big(\frac{c_{20}}{c_{19}}\Big)^\frac{3}{5} \Big\}
    \quad \mbox{and} \quad
    c_{23}:=\min \Big\{ 1 \, , \, \Big(\frac{c_{22}}{4c_{21}}\Big)^\frac{10}{13} \Big\}
  \ee
  and finally define
  \be{62.20}
    \chis := \min \bigg\{ 1 \, , \,
    \frac{\sqrt{D_1D_2} \cdot c_{23}^\frac{3}{5}}{\sqrt{2048} \cdot c_{15} c_{22}^\frac{1}{3}} \bigg\},
  \ee
  noting that through the above construction, $\chis$ indeed only depends on $D_1,D_2,a_1,a_2,\lambda_1$ and $\lambda_2$.\abs
  We now let $\chi_1>0$ and $\chi_2>0$ be such that
  \be{62.21}
    \chi_1 \le \chi_2 \le \chis,
  \ee
  and assuming (\ref{ie}) we introduce
  \be{62.22}
    \ye(t):= \io \uex^2(\cdot,t) + \gamma \io \vex^2(\cdot,t),
    \qquad t>0, \ \eps\in (0,1),
  \ee
  where
  \be{62.23}
    \gamma:=c_{23} \rho^\frac{15}{13}
    \qquad \mbox{with} \qquad
    \rho:=\frac{\chi_1}{\chi_2} \in (0,1].
  \ee
  In order to derive an appropriate upper bound for $\ye$, we first invoke Lemma \ref{lem3} to see that thanks to our
  definition of $m$ we can pick $T_1=T_1(u_0,v_0)>0$ such that
  \be{62.24}
    \io \ue + \io \ve \le m
    \qquad \mbox{for all $t>T_1$ and any } \eps\in (0,1),
  \ee
  and thereafter employ Lemma \ref{lem19} in choosing $\ts=\ts(\chi_1,\chi_2,u_0,v_0)>T_1+1$ and
  $\eps_1=\eps_1(\chi_1,\chi_2,u_0,v_0)\in (0,1)$ such that
  \be{62.25}
    \int_{\ts-1}^{\ts} \io \frac{\uex^2}{\ue}
    + \rho \int_{\ts-1}^{\ts} \io \frac{\vex^2}{\ve} \le K\cdot (1+\rho) \le 2K
  \ee
  according to (\ref{62.23}).
  Then writing $\theta:=\min_{i\in\{1,2\}} \frac{1}{2(4-n_i)}$ we fix
  $\epss=\epss(\chi_1,\chi_2,u_0,v_0) \in (0,\eps_1)$ suitably small fulfilling
  \be{62.26}
    \epss^\theta
    \le \min \bigg\{ \frac{D_1}{576 c_{16} \chi_1} \, , \,
    \frac{\gamma D_2}{32 c_6^4 c_{16} \chi_1} \, , \,
    \frac{2\gamma^\frac{39}{10}}{c_6^4 c_{16} c_{22}^3 \chi_1} \bigg\}
  \ee
  as well as
  \be{62.27}
    \epss^\theta \le \min \bigg\{
    \frac{D_2}{576 c_{16} \chi_2} \, , \,
    \frac{D_1}{32 c_6^4 c_{16} \gamma \chi_2} \, , \,
    \frac{2}{c_6^4 c_{16} c_{22}^3 \gamma^\frac{1}{10} \chi_2} \bigg\},
  \ee
  and we claim that these selections ensure that
  \be{62.28}
    \ye(t) \le c_{22} \gamma^{-\frac{3}{10}}
    \qquad \mbox{for all $t>\ts$ and any } \eps\in (0,\epss).
  \ee
  To verify this, for $\eps \in (0,1)$ we consider the set
  \be{62.288}
    S_\eps := \Big\{ t>T_1 \ \Big| \ \ye(t) < c_{22} \gamma^{-\frac{3}{10}} \Big\},
  \ee
  and first combine the inequalities provided by Lemma \ref{lem1} to see upon dropping two nonnegative summands that
  for all $\eps\in (0,1)$,
  \bea{62.29}
    & & \hspace*{-20mm}
    \frac{1}{2} \ye'(t) + \frac{D_1}{2} \io \uexx^2 + \frac{\gamma D_2}{2} \io \vexx^2 \nn\\
    &\le& \chi_1 \io \Big(h_{1,\eps}(\ue) \vex \Big)_x \uexx
    - \gamma \chi_2 \io \Big(h_{2,\eps}(\ve)\uex \Big)_x \vexx \nn\\
    & & + 3\lambda_1 \io \uex^2 + 3\gamma \lambda_2 \io \vex^2
    + \Big( \frac{a_1^2}{2D_1} + \frac{\gamma a_2^2}{2D_2} \Big) \io \ue^2 \ve^2
    \qquad \mbox{for all } t>0,
  \eea
  where $h_{i,\eps}(s):=\frac{s^{5-n_i}}{s^{4-n_i}+\eps}$ for $s\ge 0, \eps\in (0,1)$ and $i\in\{1,2\}$.
  Here since $\gamma\le 1$ by (\ref{62.23}) and (\ref{62.19}), in view of (\ref{62.24}) we may apply
  (\ref{62.3}), (\ref{62.4}), (\ref{62.10}), (\ref{62.11}) and Young's inequality to estimate
  \bea{62.30}
    & & \hspace*{-20mm}
    \Big(\frac{a_1^2}{2D_1} + \frac{\gamma a_2^2}{2D_2}\Big) \io \ue^2 \ve^2 \nn\\
    &\le& \Big(\frac{a_1^2}{2D_1} + \frac{\gamma a_2^2}{2D_2}\Big) \|\ue\|_{L^\infty(\Omega)}^2 \|\ve\|_{L^2(\Omega)}^2
        \nn\\
    &\le& \Big(\frac{a_1^2}{2D_1} + \frac{\gamma a_2^2}{2D_2}\Big) \cdot
    c_1 c_2 \cdot \Big\{ \|\uexx\|_{L^2(\Omega)}^\frac{4}{5}+1 \Big\} \cdot
    \Big\{ \|\vexx\|_{L^2(\Omega)}^\frac{2}{5} +1 \Big\} \nn\\
    &\le& \frac{D_1}{8\sqrt{2}} \cdot \Big\{ \|\uexx\|_{L^2(\Omega)}^\frac{4}{5} +1 \Big\}^\frac{5}{2}
    + c_8 \cdot \Big\{ \|\vexx\|_{L^2(\Omega)}^\frac{2}{5} +1 \Big\}^\frac{5}{3} \nn\\
    &\le& 2^\frac{3}{2} \cdot \frac{D_1}{8\sqrt{2}} \cdot \Big\{ \|\uexx\|_{L^2(\Omega)}^2 + 1 \Big\}
    + 2^\frac{2}{3} c_8 \cdot \Big\{ \|\vexx\|_{L^2(\Omega)}^\frac{2}{3}+1\Big\} \nn\\
    &=& \frac{D_1}{4} \io \uexx^2 + \frac{D_1}{4}
    + 2^\frac{2}{3} c_8 \gamma^{-\frac{1}{3}} \cdot
    \Big\{ (\sqrt{\gamma} \|\vexx\|_{L^2(\Omega)})^\frac{2}{3} + \gamma^\frac{1}{3} \Big\} \nn\\
    &\le& \frac{D_1}{4} \io \uexx^2 + \frac{D_1}{4}
    + \frac{D_2}{16} \cdot \Big\{ (\sqrt{\gamma} \|\vexx\|_{L^2(\Omega)})^\frac{2}{3} + \gamma^\frac{1}{3} \Big\}^3
    + c_9 \gamma^{-\frac{1}{2}} \nn\\
    &\le& \frac{D_1}{4} \io \uexx^2 + \frac{D_1}{4}
    + \frac{D_2}{4} \cdot \Big\{ \gamma \|\vexx\|_{L^2(\Omega)}^2 + \gamma \Big\} + c_9 \gamma^{-\frac{1}{2}} \nn\\
    &=& \frac{D_1}{4} \io \uexx^2 + \frac{D_1}{4}
    + \frac{\gamma D_2}{4} \io \vexx^2
    + \frac{\gamma D_2}{4} + c_9 \gamma^{-\frac{1}{2}} \nn\\
    &\le& \frac{D_1}{4} \io \uexx^2
    + \frac{\gamma D_2}{4} \io \vexx^2
    + c_{12} \gamma^{-\frac{1}{2}}
    \qquad \mbox{for all $t>T_1$ and } \eps\in (0,1).
  \eea
  Next, due to (\ref{62.24}), (\ref{62.6}), (\ref{62.12}) and (\ref{62.15}) we have
  \bea{62.31}
    3\lambda_1 \io \uex^2
    &\le& 3\lambda_1 c_4 \|\uexx\|_{L^2(\Omega)}^\frac{6}{5} + 3\lambda_1 c_4 \nn\\
    &\le& \frac{D_1}{8} \|\uexx\|_{L^2(\Omega)}^2
    + c_{10} + 3\lambda_1 c_4 \nn\\
    &=& \frac{D_1}{8} \io \uexx^2 + c_{13}
    \qquad \mbox{for all $t>T_1$ and } \eps\in (0,1),
  \eea
  while similarly (\ref{62.24}), (\ref{62.6}), (\ref{62.13}) and (\ref{62.15}) ensure that
  \bea{62.32}
    3\gamma \lambda_2 \io \vex^2
    &\le& 3\gamma \lambda_2 c_4\|\vexx\|_{L^2(\Omega)}^\frac{6}{5}
    + 3\gamma \lambda_2 c_4 \nn\\
    &\le& \frac{\gamma D_2}{8} \|\vexx \|_{L^2(\Omega)}^2
    + c_{11} \gamma + 3\gamma \lambda_2 c_4 \nn\\
    &\le& \frac{\gamma D_2}{8} \io \vexx^2 + c_{14}
    \qquad \mbox{for all $t>T_1$ and } \eps\in (0,1),
  \eea
  again because $\gamma\le 1$.\abs
  We next intend to estimate the cross-diffusive contributions to (\ref{62.29}), and to this end we use integration by parts
  in firstly rewriting
  \bea{62.33}
    \chi_1 \io \Big(h_{1,\eps}(\ue) \vex\Big)_x \uexx
    &=& \chi_1 \io h_{1,\eps}(\ue) \uexx \vexx
    + \frac{\chi_1}{2} \io h_{1,\eps}'(\ue) \vex (\uex^2)_x \nn\\
    &=& \chi_1 \io h_{1,\eps}(\ue) \uexx \vexx \nn\\
    & & - \frac{\chi_1}{2} \io h_{1,\eps}'(\ue) \uex^2 \vexx
    - \frac{\chi_1}{2} \io h_{1,\eps}''(\ue) \uex^3 \vex
  \eea
  for $t>0$, where by Lemma \ref{lem61}, the Cauchy-Schwarz inequality, (\ref{62.5}) and (\ref{62.24}),
  \bea{62.34}
    \chi_1 \io h_{1,\eps}(\ue) \uexx \vexx
    &\le& \chi_1 \io |\ue \uexx \vexx| \nn\\
    &\le& \chi_1 \|\ue\|_{L^\infty(\Omega)} \|\uexx\|_{L^2(\Omega)} \|\vexx\|_{L^2(\Omega)} \nn\\
    &\le& c_3\chi_1 \cdot \Big\{ \|\uex\|_{L^2(\Omega)}^\frac{2}{3} +1\Big\}
    \|\uexx\|_{L^2(\Omega)} \|\vexx\|_{L^2(\Omega)}
  \eea
  for all $t>T_1$ and $\eps\in (0,1)$.
  Moreover, thanks to Lemma \ref{lem61}  and the H\"older inequality,
  \bea{62.35}
    \hspace*{-10mm}
    - \frac{\chi_1}{2} \io h_{1,\eps}'(\ue) \uex^2 \vexx
    &\le& \frac{(5-n_1)\chi_1}{2} \io |\uex^2 \vexx| \nn\\
    &\le& \frac{(5-n_1)\chi_1}{2} \|\uex\|_{L^4(\Omega)}^2 \|\vexx\|_{L^2(\Omega)}
    \qquad \mbox{for all $t>0$ and } \eps\in (0,1),
  \eea
  and here combining (\ref{62.8}) with (\ref{62.7}), (\ref{62.24}) and Young's inequality shows that
  \bea{62.355}
    \|\uex\|_{L^4(\Omega)}^2
    &=& \|\uex\|_{L^4(\Omega)}^\frac{8}{9} \cdot \|\uex\|_{L^4(\Omega)}^\frac{10}{9} \nn\\
    &\le& c_5 c_6^\frac{8}{9} \|\uexx\|_{L^2(\Omega)}^\frac{2}{9} \|\uex\|_{L^2(\Omega)}^\frac{2}{3}
    \cdot \Big\{ \|\uexx\|_{L^2(\Omega)}^\frac{7}{9}+1 \Big\} \nn\\
    &\le& 2 c_5 c_6^\frac{8}{9} \|\uex\|_{L^2(\Omega)}^\frac{2}{3} \cdot
    \Big\{\|\uexx\|_{L^2(\Omega)} + 1 \Big\}
    \qquad \mbox{for all $t>T_1$ and } \eps\in (0,1).
  \eea
  According to (\ref{62.34}), (\ref{62.35}) and our definition (\ref{62.16}) of $c_{15}$, we therefore readily obtain
  using Young's inequality that since $\io \uex^2(\cdot,t) \le \ye(t)$ for all $t>0$ by (\ref{62.22}),
  \bea{62.356}
    & & \hspace*{-20mm}
    \chi_1 \io h_{1,\eps}(\ue) \uexx\vexx
    - \frac{\chi_1}{2} \io h_{1,\eps}'(\ue) \uex^2 \vexx \nn\\
    &\le& \Big( c_3 + (5-n_1) c_5 c_6^\frac{8}{9} \Big) \chi_1 \cdot
    \Big\{ \|\uex\|_{L^2(\Omega)}^\frac{2}{3} +1 \Big\} \cdot
    \Big\{ \|\uexx\|_{L^2(\Omega)} + 1 \Big\} \|\vexx\|_{L^2(\Omega)} \nn\\
    &\le& c_{15} \chi_1 \cdot \Big\{ \ye^\frac{1}{3}(t)+1 \Big\} \cdot
    \Big\{ \|\uexx\|_{L^2(\Omega)} +1 \Big\} \|\vexx\|_{L^2(\Omega)} \nn\\
    &\le& \frac{D_1}{64} \cdot \Big\{ \|\uexx\|_{L^2(\Omega)} +1 \Big\}^2
    + \frac{16 c_{15}^2 \chi_1^2}{D_1} \cdot \Big\{ \ye^\frac{1}{3}(t)+1\Big\}^2 \|\vexx\|_{L^2(\Omega)}^2 \nn\\
    &\le& \frac{D_1}{32} \io \uexx^2 + \frac{D_1}{32}
    + \frac{16 c_{15}^2 \chi_1^2}{D_1} \cdot \Big\{ \ye^\frac{1}{3}(t) +1 \Big\}^2 \io \vexx^2
    \qquad \mbox{for all $t>T_1$ and } \eps\in (0,1).
  \eea
  Here we note that  since $c_{22}\ge 1$ and $\gamma\le 1$ and thus
  $1 \le c_{22}^\frac{1}{3} \gamma^{-\frac{1}{10}}$ by (\ref{62.19}), in view of (\ref{62.23}) and (\ref{62.20})
  inside the set $S_\eps$ in (\ref{62.288}) we may estimate
  \bas
    \frac{\frac{16 c_{15}^2 \chi_1^2}{D_1} \cdot \Big\{ \ye^\frac{1}{3}(t)+1\Big\}^2}{\frac{\gamma D_2}{32}}
    &=& \frac{512 c_{15}^2 \chi_1^2}{\gamma D_1 D_2} \cdot \Big\{ \ye^\frac{1}{3}(t) +1 \Big\}^2 \nn\\
    &\le& \frac{512 c_{15}^2 \chi_1^2}{\gamma D_1 D_2} \cdot (c_{22}^\frac{1}{3} \gamma^{-\frac{1}{10}} +1)^2 \nn\\
    &\le& \frac{2048 c_{15}^2 c_{22}^\frac{2}{3}}{D_1 D_2} \cdot \chi_1^2 \gamma^{-\frac{6}{5}} \nn\\
    &=& \frac{2048 c_{15}^2 c_{22}^\frac{2}{3} c_{23}^{-\frac{6}{5}}}{D_1 D_2} \cdot \chi_1^\frac{8}{13} \chi_2^\frac{18}{13} \nn\\
    &\le& \frac{2048 c_{15}^2 c_{22}^\frac{2}{3} c_{23}^{-\frac{6}{5}}}{D_1 D_2} \cdot \chis^2 \nn\\[2mm]
    &\le& 1
    \qquad \mbox{for all $t\in S_\eps$ and each } \eps \in (0,1),
  \eas
  whence (\ref{62.356}) implies that actually
  \bea{62.36}
    & & \hspace*{-20mm}
    \chi_1 \io h_{1,\eps}(\ue) \uexx\vexx
    - \frac{\chi_1}{2} \io h_{1,\eps}'(\ue) \uex^2 \vexx \nn\\
    &\le& \frac{D_1}{32} \io \uexx^2
    + \frac{\gamma D_2}{32} \io \vexx^2
    + \frac{D_1}{32}
    \qquad \mbox{for all $t\in S_\eps$ and any } \eps \in (0,1).
  \eea
  In the rightmost summand in (\ref{62.33}), once more relying on Lemma \ref{lem61} enables us to combine Young's inequality
  with Lemma \ref{lem63} to see that due to (\ref{62.16}) and our definition of $\theta$, and again thanks to (\ref{62.8}),
  for all $t>0$ and $\eps\in (0,1)$ we have
  \bea{62.366}
    - \frac{\chi_1}{2} \io h_{1,\eps}''(\ue) \uex^3 \vex
    &\le& c_{16} \chi_1 \eps^\frac{1}{2(4-n_1)} \io \Big| \frac{\uex^3 \vex}{\ue^\frac{3}{2}} \Big| \nn\\
    &\le& c_{16} \chi_1 \eps^\theta \io \frac{\uex^4}{\ue^2}
    + c_{16} \chi_1 \eps^\theta \io \vex^4 \nn\\
    &\le& 9c_{16} \chi_1 \eps^\theta \io \uexx^2
    + c_6^4 c_{16} \chi_1 \eps^\theta \|\vexx\|_{L^2(\Omega)} \|\vex\|_{L^2(\Omega)}^3 \nn\\
    &\le& 9c_{16} \chi_1 \eps^\theta \io \uexx^2 \nn\\
    & & + \frac{1}{2} c_6^4 c_{16} \chi_1 \eps^\theta \io \vexx^2
    + \frac{1}{2} c_6^4 c_{16} \chi_1 \eps^\theta \gamma^{-3} \ye^3(t),
  \eea
  because $\io \vex^2(\cdot,t)\le \gamma^{-1} \ye(t)$ for $t>0$ by (\ref{62.22}).
  Again restricted to times belonging to $S_\eps$, due to
  (\ref{62.288}) and (\ref{62.26}) this implies that
  \bea{62.37}
    \hspace*{-10mm}
    - \frac{\chi_1}{2} \io h_{1,\eps}''(\ue) \uex^3 \vex
    &\le& 9c_{16} \chi_1 \eps^\theta \io \uexx^2 \nn\\
    & & + \frac{1}{2} c_6^4 c_{16} \chi_1 \eps^\theta \io \vexx^2
    + \frac{1}{2} c_6^4 c_{16} c_{22}^3 \chi_1 \gamma^{-\frac{39}{10}} \eps^\theta \nn\\
    &\le& \frac{D_1}{64} \io \uexx^2
    + \frac{\gamma D_2}{64} \io \vexx^2
    + 1
    \qquad \mbox{for all $t\in S_\eps$ and } \eps\in (0,\epss).
  \eea
  Now by arguments quite similar to those used in (\ref{62.33})-(\ref{62.37}), on splitting
  \bea{62.38}
    - \gamma\chi_2 \io \Big(h_{2,\eps}(\ve) \uex \Big)_x \vexx
    &=& - \gamma\chi_2 \io h_{2,\eps}(\ve) \uexx \vexx
    + \frac{\gamma\chi_2}{2} \io h_{2,\eps}'(\ve) \vex^2 \uexx \nn\\
    & & + \frac{\gamma\chi_2}{2} \io h_{2,\eps}''(\ve) \uex \vex^3,
    \qquad t>0,
  \eea
  we may first essentially copy the reasoning from (\ref{62.34}), (\ref{62.35}), (\ref{62.355}), (\ref{62.356})
  and (\ref{62.36}): Indeed, again using that $\io \vex^2(\cdot,t)\le \gamma^{-1} \ye(t)$ for $t>0$
  we thereby obtain that
  \bea{62.39}
    & & \hspace*{-20mm}
    - \gamma\chi_2 \io h_{2,\eps}(\ve) \uexx \vexx
    + \frac{\gamma\chi_2}{2} \io h_{2,\eps}'(\ve) \vex^2 \uexx \nn\\
    &\le& c_{15} \gamma \chi_2 \cdot \Big\{ \|\vex\|_{L^2(\Omega)}^\frac{2}{3} +1 \Big\} \|\uexx\|_{L^2(\Omega)}
    \Big\{ \|\vexx\|_{L^2(\Omega)}+1\Big\} \nn\\
    &\le& c_{15} \gamma \chi_2 \cdot \Big\{ \gamma^{-\frac{1}{3}} \ye^\frac{1}{3}(t)+1\Big\}
    \|\uexx\|_{L^2(\Omega)} \Big\{ \|\vexx\|_{L^2(\Omega)}+1\Big\} \nn\\
    &\le& \frac{D_1}{32} \|\uexx\|_{L^2(\Omega)}^2
    +\frac{8 c_{15}^2 \gamma^2 \chi_2^2}{D_1}
    \cdot \Big\{ \gamma^{-\frac{1}{3}}\ye^\frac{1}{3}(t)+1\Big\}^2 \cdot
    \Big\{ \|\vexx\|_{L^2(\Omega)} +1 \Big\}^2 \nn\\
    &\le& \frac{D_1}{32} \io \uexx^2
    + \frac{16 c_{15}^2 \gamma^2 \chi_2^2}{D_1} \cdot
    \Big\{ \gamma^{-\frac{1}{3}} \ye^\frac{1}{3}(t)+1\Big\}^2 \cdot \bigg\{ \io \vexx^2 + 1 \bigg\} \nn\\
    &\le& \frac{D_1}{32} \io \uexx^2
    + \frac{\gamma D_2}{32} \io \vexx^2
    + \frac{\gamma D_2}{32}
    \qquad \mbox{for all $t\in S_\eps$ and } \eps\in (0,1),
  \eea
  because according to (\ref{62.288}) and the inequalities $\gamma\le 1$ and
  $c_{22}^\frac{1}{3} \gamma^{-\frac{13}{30}} \ge 1$ asserted by (\ref{62.19}) and (\ref{62.23}) we have
  \bas
    \frac{\frac{16 c_{15}^2 \gamma^2 \chi_2^2}{D_1} \cdot \Big\{ \gamma^{-\frac{1}{3}} \ye^\frac{1}{3}(t)+1\Big\}^2}
    {\frac{\gamma D_2}{32}}
    &=& \frac{512 c_{15}^2 \gamma \chi_2^2}{D_1 D_2} \cdot \Big\{ \gamma^{-\frac{1}{3}} \ye^\frac{1}{3}(t)+1\Big\}^2 \\
    &<& \frac{512 c_{15}^2 \gamma \chi_2^2}{D_1 D_2} \cdot (c_{22}^\frac{1}{3} \gamma^{-\frac{13}{30}} +1)^2 \\
    &\le& \frac{2048 c_{15}^2 c_{22}^\frac{2}{3} \gamma^\frac{2}{15} \chi_2^2}{D_1 D_2} \\
    &\le& \frac{2048 c_{15}^2 c_{22}^\frac{2}{3} \chis^2}{D_1 D_2} \nn\\[2mm]
    &\le& 1
    \qquad \mbox{for all $t\in S_\eps$ and } \eps\in (0,1)
  \eas
  due to (\ref{62.20}) and our restriction that $1 \ge c_{23}^\frac{3}{5}$.\abs
  Likewise, proceeding as in (\ref{62.366}) and (\ref{62.37}) we find that since $\io \uex^2(\cdot,t) \le \ye(t)
  < c_{22} \gamma^{-\frac{3}{10}}$ for all $t\in S_\eps$,
  \bea{62.40}
    \hspace*{-10mm}
    \frac{\gamma \chi_2}{2} \io h_{2,\eps}''(\ve) \uex \vex^3
    &\le& 9c_{16} \gamma \chi_2 \eps^\theta \io \vexx^2
    + c_6^4 c_{16} \gamma \chi_2 \eps^\theta \|\uexx\|_{L^2(\Omega)} \|\uex\|_{L^2(\Omega)}^3 \nn\\
    &\le& 9c_{16} \gamma \chi_2 \eps^\theta \io \vexx^2
    + \frac{1}{2} c_6^4 c_{16} \gamma \chi_2 \eps^\theta \io \uexx^2
    + \frac{1}{2} c_6^4 c_{16} c_{22}^3 \gamma^\frac{1}{10} \chi_2 \eps^\theta \nn\\
    &\le& \frac{D_1}{64} \io \uexx^2
    + \frac{\gamma D_2}{64} \io \vexx^2
    + 1
    \qquad \mbox{for all $t\in S_\eps$ and } \eps\in (0,\epss)
  \eea
  because of (\ref{62.27}).\abs
  When inserted into (\ref{62.29}), in view of (\ref{62.33}) and (\ref{62.38}) the estimates (\ref{62.30}),
  (\ref{62.31}), (\ref{62.32}), (\ref{62.36}), (\ref{62.37}), (\ref{62.39}) and (\ref{62.40}) in summary show that
  as $\gamma\le 1$,
  \bea{62.41}
    \frac{1}{2} \ye'(t)
    + \frac{D_1}{32} \io \uexx^2
    + \frac{\gamma D_2}{32} \io \vexx^2
    &\le& c_{12} \gamma^{-\frac{1}{2}} + c_{13} + c_{14} + \frac{D_1}{32} + 1 + \frac{\gamma D_2}{32} + 1 \nn\\
    &\le& c_{17} \gamma^{-\frac{1}{2}}
    \qquad \mbox{for all $t\in S_\eps$ and } \eps\in (0,\epss)
  \eea
  according to the definition of $c_{17}$ in (\ref{62.17}).
  In order to turn the two remaining integrals herein into an appropriate superlinear absorptive term,
  we once more make use of (\ref{62.6}), (\ref{62.24}) and Young's inequality to see that again since $\gamma\le 1$,
  \bas
    \ye^\frac{5}{3}(t)
    &=& \bigg\{ \io \uex^2 + \gamma \io \vex^2 \bigg\}^\frac{5}{3} \\
    &\le& c_4^\frac{5}{3} \cdot \bigg\{ \|\uexx\|_{L^2(\Omega)}^\frac{6}{5}+1
    + \gamma \cdot \Big( \|\vexx\|_{L^2(\Omega)}^\frac{6}{5}+1\Big) \bigg\}^\frac{5}{3} \\
    &\le& c_4^\frac{5}{3} \cdot \Big\{ \|\uexx\|_{L^2(\Omega)}^\frac{6}{5}
    + \gamma^\frac{3}{5} \|\vexx\|_{L^2(\Omega)}^\frac{6}{5} + 2 \Big\}^\frac{5}{3} \\
    &\le& (3c_4)^\frac{5}{3} \cdot \Big\{ \|\uexx\|_{L^2(\Omega)}^2 + \gamma \|\vexx\|_{L^2(\Omega)}^2 + 2^\frac{5}{3}
    \Big\}
    \qquad \mbox{for all $t>T_1$ and all } \eps\in (0,1),
  \eas
  so that for all $t>T_1$ and any $\eps\in (0,1)$,
  \bas
    \frac{D_1}{32} \io \uexx^2
    + \frac{\gamma D_2}{32} \io \vexx^2
    &\ge& \frac{\min\{D_1,D_2\}}{32} \cdot \bigg \{ \io \uexx^2 + \gamma \io \vexx^2 \bigg\} \\
    &\ge& \frac{\min\{D_1,D_2\}}{32} \cdot (3c_4)^{-\frac{5}{3}} \cdot \ye^\frac{5}{3}(t)
    - \frac{2^\frac{2}{3} \cdot \min\{D_1,D_2\}}{16}.
  \eas
  In light of (\ref{62.18}), (\ref{62.41}) therefore entails the autonomous ODI
  \bea{62.42}
    \frac{1}{2} \ye'(t) + c_{19} \ye^\frac{5}{3}(t)
    \le c_{20} \gamma^{-\frac{1}{2}}
    \qquad \mbox{for all $t\in S_\eps$ and each } \eps\in (0,\epss).
  \eea
  To deduce (\ref{62.28}) from this, we now recall (\ref{62.25}) to infer that for each $\eps\in (0,1)$
  we can fix $t_\eps\in (\ts-1,\ts)$ such that
  \be{62.43}
    \io \frac{\uex^2(\cdot,t_\eps)}{\ue(\cdot,t_\eps)}
    + \rho \io \frac{\vex^2(\cdot,t_\eps)}{\ve(\cdot,t_\eps)}
    \le 2K.
  \ee
  Through (\ref{62.9}) and again (\ref{62.24}), this firstly entails that
  \bas
    \|\ue(\cdot,t_\eps)\|_{L^\infty(\Omega)}^2
    &=& \|\sqrt{\ue(\cdot,t_\eps)}\|_{L^\infty(\Omega)}^4 \nn\\
    &\le& c_7 \cdot \Big\{ \|(\sqrt{\ue(\cdot,t_\eps)})_x\|_{L^2(\Omega)}^2 +1 \Big\} \\
    &=& \frac{c_7}{4} \io \frac{\uex^2(\cdot,t_\eps)}{\ue(\cdot,t_\eps)} + c_7 \\
    &\le& \frac{c_7 K}{2} + c_7
    \qquad \mbox{for all } \eps\in (0,1)
  \eas
  and, similarly,
  \bas
    \|\ve(\cdot,t_\eps)\|_{L^\infty(\Omega)}^2
    \le \frac{c_7}{4} \io \frac{\vex^2(\cdot,t_\eps)}{\ve(\cdot,t_\eps)} + c_7
    \le \frac{c_7 K}{2\rho} + c_7
    \qquad \mbox{for all } \eps\in (0,1).
  \eas
  Therefore, (\ref{62.43}) secondly implies that
  \bea{62.44}
    \ye(t_\eps)
    &\le& \|\ue(\cdot,t_\eps)\|_{L^\infty(\Omega)} \io \frac{\uex^2(\cdot,t_\eps)}{\ue(\cdot,t_\eps)}
    + \gamma \|\ve(\cdot,t_\eps)\|_{L^\infty(\Omega)} \io \frac{\vex^2(\cdot,t_\eps)}{\ve(\cdot,t_\eps)} \nn\\
    &\le& \sqrt{\frac{c_7 K}{2} + c_7} \cdot 2K
    + \gamma \cdot \sqrt{\frac{c_7 K}{2\rho} + c_7} \cdot \frac{2K}{\rho} \nn\\
    &\le& c_{21} + c_{21} \gamma \rho^{-\frac{3}{2}}
    \qquad \mbox{for all } \eps\in (0,1)
  \eea
  by (\ref{62.19}), for we are yet assuming that $\rho\le 1$. \abs
  But now the precise link between $\gamma$ and $\rho$ in (\ref{62.23}) ensures that
  \bas
    \frac{c_{21} \gamma \rho^{-\frac{3}{2}}}{\frac{1}{4} c_{22} \gamma^{-\frac{3}{10}}}
    = \frac{4 c_{21}}{c_{22}} \gamma^\frac{13}{10} \rho^{-\frac{3}{2}}
    = \frac{ 4c_{21} c_{23}^\frac{13}{10}}{c_{22}} \le 1
  \eas
  due to (\ref{62.19}),
  so that since (\ref{62.19}) moreover entails that
  \bas
    \frac{c_{21}}{\frac{1}{4} c_{22} \gamma^{-\frac{3}{10}}} \le \frac{4 c_{21}}{c_{22}} \le 1,
  \eas
  from (\ref{62.44}) we infer that
  \be{62.444}
    \ye(t_\eps) \le \frac{1}{2} c_{22} \gamma^{-\frac{3}{10}}
    \qquad \mbox{for all } \eps\in (0,1).
  \ee
  This especially guarantees that $t_\eps$ belongs to $S_\eps$ and that hence
  \bas
    T_\eps:=\sup \Big\{ \widehat{T}>t_\eps \ \Big| \ [t_\eps,\widehat{T}] \subset S_\eps \Big\}
  \eas
  is a well-defined element of $(t_\eps,\infty]$ for all $\eps\in (0,1)$.
  According to (\ref{62.42}), however, for each $\eps\in (0,\epss)$ we actually must have $T_\eps=\infty$:
  Otherwise, namely, the definition of $S_\eps$ would entail that $(t_\eps,T_\eps) \subset S_\eps$ but
  \be{62.45}
    \ye(T_\eps)= c_{22} \gamma^{-\frac{3}{10}},
  \ee
  while since $\oy(t):=\frac{1}{2} c_{22} \gamma^{-\frac{3}{10}}$, $t\ge t_\eps$, satisfies
  \bas
    \frac{1}{2} \oy'(t) + c_{19} \oy^\frac{5}{3}(t) - c_{20} \gamma^{-\frac{1}{2}}
    = (2^{-\frac{5}{3}} c_{19} c_{22}^\frac{5}{3} - c_{20}) \cdot \gamma^{-\frac{1}{2}} \ge 0
    \qquad \mbox{for all } t>t_\eps
  \eas
  as a further consequence of (\ref{62.19}), a comparison argument based on (\ref{62.444}) would show that
  $\ye(t) \le \oy(t)$ for all $t\in [t_\eps,T_\eps]$ and thereby contradict (\ref{62.45}).
  This confirms that indeed $[t_\eps,\infty) \subset S_\eps$ for all $\eps\in (0,\epss)$ and that thus
  (\ref{62.28}) results from the fact that $t_\eps < \ts$ for any such $\eps$.
  In view of (\ref{62.24}) and e.g.~(\ref{62.5}), we therefore conclude that both (\ref{62.1}) and (\ref{62.2})
  are valid whenever $0<\chi_1\le \chi_2 < \chis$, whereupon going back to (\ref{62.41}) we infer by integration
  that also (\ref{62.02}) holds for all such $\chi_1$ and $\chi_2$.\abs
  According to an evident symmetry property of the presented
  reasoning, however, in the case when conversely $\chi_2\le \chi_1$
  the above result can be derived in much the same manner, essentially by exchanging the roles of $\ue$ and $\ve$.
\qed
By once more recalling the convergence properties obtained in Section \ref{sect_conv}, from Lemma \ref{lem62}
we readily derive the following conclusion, inter alia asserting eventual continuity and $H^1$-boundedness of the limit couple
from Lemma \ref{lem_conv}.
\begin{cor}\label{cor24}
  Let $D_i>0, a_i>0, \lambda_i>0$ and $n_i\in [1,2]$ for $i\in\{1,2\}$,
  and let $\chi_1 \in (0,\chis)$ and $\chi_2\in (0,\chis)$
  with $\chis>0$ as given by Lemma \ref{lem62}.
  Then whenever (\ref{ie}) holds, there exists $\ts=\ts(\chi_1,\chi_2,u_0,v_0)>0$
  such that with $(\eps_j)_{j\in\N}\subset (0,1)$ and $N\subset (0,\infty)$ taken from Lemma \ref{lem_conv}
  we have
  \be{24.1}
    \ue(\cdot,t) \to u(\cdot,t)
    \quad \mbox{and} \quad
    \ve(\cdot,t) \to v(\cdot,t)
    \quad \mbox{in } L^\infty(\Omega)
    \qquad \mbox{for all } t\in (0,\infty)\setminus N
  \ee
  as well as
  \be{24.2}
    \int_{t_1}^{t_2} \io \uex^2 \to \int_{t_1}^{t_2} \io u_x^2
    \quad \mbox{and} \quad
    \int_{t_1}^{t_2} \io \vex^2  \to \int_{t_1}^{t_2} \io v_x^2
    \qquad \mbox{for all $t_1>\ts$ and any } t_2>t_1
  \ee
  as $\eps=\eps_j\searrow 0$.
  Moreover, both $u$ and $v$ belong to $C^0(\bom\times (\ts,\infty))$ and satisfy $u(\cdot,t)\in W^{1,2}(\Omega)$
  and $v(\cdot,t) \in W^{1,2}(\Omega)$ for all $t>\ts$, and one can find $C>0$ such that
  \be{24.3}
    \|u(\cdot,t)\|_{W^{1,2}(\Omega)} + \|v(\cdot,t)\|_{W^{1,2}(\Omega)} \le C
    \qquad \mbox{for all } t>\ts.
  \ee
\end{cor}
\proof
  According to Lemma \ref{lem62}, for fixed $u_0$ and $v_0$ we can find $\ts>0$, $c_1>0$, $c_2>0$ and $\epss\in (0,1)$
  such that for any $\eps\in (0,\epss)$,
  \be{24.4}
    \|\ue(\cdot,t)\|_{W^{1,2}(\Omega)}
    + \|\ve(\cdot,t)\|_{W^{1,2}(\Omega)} \le c_1
    \qquad \mbox{for all } t>\ts
  \ee
  and
  \be{24.5}
    \int_{\ts}^T \|\ue(\cdot,t)\|_{W^{2,2}(\Omega)}^2 dt
    + \int_{\ts}^T \|\ve(\cdot,t)\|_{W^{2,2}(\Omega)}^2 dt
    \le c_2 \cdot (T+1)
    \qquad \mbox{for all } T>\ts.
  \ee
  As (\ref{24.4}) warrants relative compactness of $(\ue(\cdot,t))_{\eps\in (0,\epss)}$ and of
  $(\ve(\cdot,t))_{\eps\in (0,\epss)}$ with respect to the weak topology in $W^{1,2}(\Omega)$ for all $t>\ts$,
  it is thus clear from Lemma \ref{lem_conv} that as $\eps=\eps_j\searrow 0$,
  \be{24.6}
    \ue(\cdot,t) \wto u(\cdot,t)
    \quad \mbox{and} \quad
    \ve(\cdot,t) \wto v(\cdot,t)
    \quad \mbox{in } W^{1,2}(\Omega)
    \qquad \mbox{for all } t\in (\ts,\infty)\setminus N.
  \ee
  and thereby especially entails (\ref{24.1}) by compactness of the embedding $W^{1,2}(\Omega) \hra L^\infty(\Omega)$.
  Since thus $(u(\cdot,t))_{t\in (\ts,\infty)\setminus N}$ and $(v(\cdot,t))_{t\in (\ts,\infty)\setminus N}$
  are bounded in $W^{1,2}(\Omega)$ due to (\ref{24.6}) and (\ref{24.4}), and since from Lemma \ref{lem_conv} we already
  know that both $u$ and $v$ belong to $C^0_w([0,\infty);L^1(\Omega))$, it can therefore easily be verified that
  actually $u$ and $v$ are continuous throughout $(\ts,\infty)$ as $W^{1,2}(\Omega)$-valued functions with respect
  to the weak topology therein.
  Again through compactness of the embedding $W^{1,2}(\Omega) \hra L^\infty(\Omega)$, this entails
  that indeed $u$ and $v$ are contained in $C^0((\ts,\infty);C^0(\bom))$ and satisfy (\ref{24.3}).\\
  Apart from that, for all $T>\ts$ (\ref{24.5}) ensures boundedness of the families $(\ue)_{\eps\in (0,\epss)}$ and
  $(\ve)_{\eps\in (0,\epss)}$ in $L^2((\ts,T);W^{2,2}(\Omega))$, so that recalling the time regularity properties
  from Lemma \ref{lem52} we may once more invoke the Aubin-Lions lemma to conclude that
  $\ue \to u$ and $\ve \to v$ in $L^2((\ts,t);W^{1,2}(\Omega))$, and that thus in particular (\ref{24.2}) holds.
\qed
\mysection{Toward verifying consistency with a genuine entropy structure}
Now having at hand the eventual regularity properties gained above, for establishing our main results on stabilization
in the flavor of Theorem \ref{theo33} and Theorem \ref{theo333} it will be sufficient to make sure that
the genuine gradient-like structures related to the functionals in (\ref{entropy2}) and (\ref{entropy3}) find some
appropriate rigorous counterpart in the approximate problems (\ref{0eps})-(\ref{0epsb}).
This short section provides the fundament therefor by recording results of associated testing procedures,
which will in particular foreshadow our final selections of the key parameters $n_1$ and $n_2$.
\begin{lem}\label{lem22}
  Let $D_i>0, a_i>0, \lambda_i>0$ and $n_i\in [1,2]$ for $i\in\{1,2\}$,
  and assume (\ref{ie}).\abs
  i) \ If $n_1=2$, then for all $\eps\in (0,1)$,
  \bea{22.1}
    \hspace*{-10mm}
    \frac{d}{dt} \bigg\{ - \io \ln\ue + \frac{\eps}{6} \io \frac{1}{\ue^2} \bigg\}
    &+& \eps\io \uexx^2
    + D_1 \io \frac{\uex^2}{\ue^2}
    + D_1 \eps \io \frac{\uex^2}{\ue^4} \nn\\
    &+& \eps^\frac{\alpha}{2} \io \ue^{-\alpha-2} \uex^2
    +  \eps^\frac{\alpha+2}{2} \io \ue^{-\alpha-4} \uex^2 \nn\\
    &=& \chi_1 \io \frac{1}{\ue} \uex\vex
    - \lambda_1 |\Omega|
    + \io \ue - a_1 \io \ve
    \qquad \mbox{for all } t>0.
  \eea
  ii) \ If $n_2=2$, then for all $\eps\in (0,1)$,
  \bea{22.2}
    \hspace*{-10mm}
    \frac{d}{dt} \bigg\{ - \io \ln\ve + \frac{\eps}{6} \io \frac{1}{\ve^2} \bigg\}
    &+& \eps\io \vexx^2
    + D_2 \io \frac{\vex^2}{\ve^2}
    + D_2 \eps \io \frac{\vex^2}{\ve^4} \nn\\
    &+& \eps^\frac{\alpha}{2} \io \ve^{-\alpha-2} \vex^2
    +  \eps^\frac{\alpha+2}{2} \io \ve^{-\alpha-4} \vex^2 \nn\\
    &=& - \chi_2 \io \frac{1}{\ve} \uex\vex
    - \lambda_2 |\Omega|
    + \io \ve + a_2 \io \ue
    \qquad \mbox{for all } t>0.
  \eea
  iii) \ In the case $n_2=1$, for all $\eps\in (0,1)$ we have
  \bea{22.3}
    \hspace*{-10mm}
    \frac{d}{dt} \bigg\{ \frac{1}{2} \io \ve^2 + \frac{\eps}{2} \io \frac{1}{\ve} \bigg\}
    &+& \eps\io \ve \vexx^2
    + D_2 \io \vex^2
    + D_2 \eps \io \frac{\vex^2}{\ve^3} \nn\\
    &+& \eps^\frac{\alpha}{2} \io \ve^{-\alpha} \vex^2
    +  \eps^\frac{\alpha+2}{2} \io \ve^{-\alpha-3} \vex^2 \nn\\
    &=& - \chi_2 \io \ve \uex\vex
    + \lambda_2 \io \ve^2
    - \io \ve^3 - a_2 \io \ue\ve^2 \nn\\
    & & - \eps \io  \frac{2\ve^2+3\ve}{6\ve^2 +2\eps}
      (\lambda_2-\ve-a_2\ue)
    \qquad \mbox{for all } t>0.
  \eea
\end{lem}
\proof
  i) \ Writing $\ell(s):=-\ln s + \frac{\eps}{6} \cdot \frac{1}{s^2}$ for $s>0$ and fixed $\eps\in (0,1)$, we have
  $\ell'(s)=-\frac{1}{s} - \frac{\eps}{3s^3}= - \frac{3s^2+\eps}{3s^3}$
  and $\ell''(s)=\frac{1}{s^2}+\frac{\eps}{s^4}=\frac{s^2+\eps}{s^4}$ for all $s>0$, so that using the first equation
  in (\ref{0eps}) we obtain
  \bas
    \frac{d}{dt} \io \ell(\ue)
    &=& - \eps \io \ell'(\ue) \cdot \Big( \frac{\ue^4}{\ue^2+\eps} \uexxx\Big)_x
    + \eps^\frac{\alpha}{2} \io \ell'(\ue) \cdot (\ue^{-\alpha} \uex)_x \nn\\
    & & + D_1 \io \ell'(\ue) \uexx
    - \chi_1 \io \ell'(\ue) \cdot \Big( \frac{\ue^3}{\ue^2+\eps} \vex \Big)_x \nn\\
    & & + \io \ell'(\ue) \cdot \frac{3\ue^3}{3\ue^2+\eps} \cdot (\lambda_1-\ue+a_1\ve) \nn\\
    &=& \eps \io \ell''(\ue) \cdot \frac{\ue^4}{\ue^2+\eps} \uex\uexxx
    - \eps^\frac{\alpha}{2} \io \ell''(\ue) \ue^{-\alpha} \uex^2 \nn\\
    & & - D_1 \io \ell''(\ue) \uex^2
    + \chi_1 \io \ell''(\ue) \cdot \frac{\ue^3}{\ue^2+\eps} \uex\vex \nn\\
    & & + \io \ell'(\ue) \cdot \frac{3\ue^3}{3\ue^2+\eps} \cdot (\lambda_1-\ue+a_1\ve) \nn\\
    &=& \eps \io \uex \uexxx
    - \eps^\frac{\alpha}{2} \io (\ue^{-\alpha-2} + \eps \ue^{-\alpha-4})\uex^2 \nn\\
    & & - D_1 \io \Big( \frac{1}{\ue^2} + \frac{\eps}{\ue^4}\Big) \uex^2
    + \chi_1 \io \frac{1}{\ue} \uex\vex \nn\\
    & & - \io (\lambda_1 - \ue +a_1\ve)
    \qquad \mbox{for all } t>0.
  \eas
  As clearly $\eps\io \uex\uexxx=-\eps\io \uexx^2$ for all $t>0$, this is equivalent to (\ref{22.1}).\abs
  ii) \ This part analogously follows from the second equation in (\ref{0eps}).\abs
  iii) \ We now rather let $\ell(s):=\frac{1}{2} s^2 + \frac{\eps}{2s}$ for $s>0$ and $\eps\in (0,1)$,
  and observing that then $\ell'(s)=s-\frac{\eps}{2s^2}=\frac{2s^3-\eps}{2s^2}$ and
  $\ell''(s)=1+\frac{\eps}{s^3}=\frac{s^3+\eps}{s^3}$ for all $s>0$, we use (\ref{0eps})-(\ref{0epsb}) to compute
  \bas
    \frac{d}{dt} \io \ell(\ve)
    &=& \eps \io \ell''(\ve) \cdot \frac{\ve^4}{\ve^3+\eps} \vex \vexxx
    - \eps^\frac{\alpha}{2} \io \ell''(\ve) \ve^{-\alpha} \vex^2 \nn\\
    & & - D_2 \io \ell''(\ve) \vex^2
    - \chi_2 \io \ell''(\ve) \cdot \frac{\ve^4}{\ve^3+\eps} \uex\vex \nn\\
    & & + \io \ell'(\ve) \cdot \frac{3\ve^3}{3\ve^2+\eps} \cdot (\lambda_2 - \ve -a_2\ue) \nn\\
    &=& \eps \io \ve \vex \vexxx
    - \eps^\frac{\alpha}{2} \io (\ve^{-\alpha} + \eps \ve^{-\alpha-3})\vex^2 \nn\\
    & & - D_2 \io \Big(1+\frac{\eps}{\ve^3}\Big) \vex^2
    - \chi_2 \io \ve\uex\vex \nn\\
    & & + \io \frac{3\ve (2\ve^3-\eps)}{2\cdot (3\ve^2+\eps)} \cdot (\lambda_2-\ve-a_1\ue)
    \qquad \mbox{for all } t>0.
  \eas
  Since
  \bas
    \eps\io \ve\vex\vexxx
    = - \eps \io \ve\vexx^2 - \eps \io \vex^2 \vexx
    = - \eps \io \ve \vexx^2
    \qquad \mbox{for all } t>0,
  \eas
  and since
  \bas
    \frac{3s(2s^3-\eps)}{2\cdot (3s^2+\eps)}
    = s^2-\eps \cdot \frac{2s^2+3s}{6s^2+2\eps}
    \qquad \mbox{for all } s>0,
  \eas
  this yields (\ref{22.3}).
\qed
\mysection{The case $\lambda_2>a_2\lambda_1$. Proof of Theorem \ref{theo33}}\label{sect_asy1}
In our convergence argument concentrating on the situation of Theorem \ref{theo33} first, in view of (\ref{entropy2})
our design of an approximate variant of the functional therein will be based on parts i) and ii) of Lemma \ref{lem22},
thus suggesting to fix $n_1=n_2=2$ henceforth in this case.
For our construction here and also for later reference,
let us separately state some elementary properties of the nonlinear map arising in both integrals in (\ref{entropy2}).
\begin{lem}\label{lem29}
  Let $\xis>0$ and
  \be{phi}
    \phi_{\xis}(\xi):=\xi-\xis-\xis \ln \frac{\xi}{\xis},
    \qquad \xi>0.
  \ee
  Then $\phi_{\xis}$ is positive on $(0,\infty) \setminus \{\xis\}$ with $\phi_{\xis}(\xis)=0$, and furthermore
  \be{29.2}
    \phi_{\xis}(\xi)
    \le \frac{2}{\xis} \cdot (\xi-\xis)^2
    \qquad \mbox{for all } \xi\ge \frac{\xis}{2}.
  \ee
\end{lem}
\proof
  Since $\phi_{\xis}'(\xi)=1-\frac{\xis}{\xi}$ and $\phi_{\xis}''(\xi)=\frac{\xis}{\xi^2}$ for all $\xi>0$,
  it is clear that the zero $\xis$ of $\phi_{\xis}$ is a strict and global minimum point of $\phi_{\xis}$, and that
  moreover
  \bas
    \phi_{\xis}(\xi)
    \le \sup_{\eta\ge\frac{\xis}{2}} |\phi_{\xis}''(\eta)| \cdot \frac{(\xi-\xis)^2}{2}
    = \frac{2}{\xis} \cdot (\xi-\xis)^2
  \eas
  for all $\xi\ge \frac{\xis}{2}$.
\qed
Then combining lemma \ref{lem22} with (\ref{entropy2}) suggests to generalize the entropy structure in question
as follows.
\begin{lem}\label{lem23}
  Let $n_1=n_2=2$, and suppose that
  $D_i>0, a_i>0, \lambda_i>0$ and $\chi_i>0$ for $i\in \{1,2\}$, with
  \be{23.1}
    \lambda_2>a_2 \lambda_1.
  \ee
  Then writing
  \be{A}
    A:=\frac{a_1}{a_2}
  \ee
  and taking $\us>0$ and $\vs>0$ from (\ref{usvs}), for
  \be{Eo}
    \Eo(t) := \io \phi_{\us}(\ue(\cdot,t))
    + \frac{\us \eps}{6} \io \frac{1}{\ue^2(\cdot,t)}
    + A \io \phi_{\vs}(\ve(\cdot,t))
    + \frac{A \vs \eps}{6} \io \frac{1}{\ve^2(\cdot,t)},
    \qquad t\ge 0, \ \eps\in (0,1),
  \ee
  with $\phi_{\xis}$ taken from (\ref{phi}) for $\xi>0$, we have
  \bea{23.2}
    \frac{d}{dt} \Eo(t)
    &+& \bigg\{ \frac{D_1 \us}{2} - \frac{A \chi_2^2 \vs}{2D_2} \|\ue(\cdot,t)\|_{L^\infty(\Omega)}^2 \bigg\}
    \cdot \io \frac{\uex^2}{\ue^2} \nn\\
    &+& \bigg\{ \frac{A D_2 \vs}{2} - \frac{\chi_1^2 \us}{2D_1} \|\ve(\cdot,t)\|_{L^\infty(\Omega)}^2 \bigg\}
    \cdot \io \frac{\vex^2}{\ve^2} \nn\\
    &+& \io (\ue(\cdot,t)-\us)^2
    + A \io (\ve(\cdot,t)-\vs)^2 \nn\\
    &+& \us \eps^\frac{\alpha+2}{2} \io \ue^{-\alpha-4} \uex^2
    + A\vs \eps^\frac{\alpha+2}{2} \io \ve^{-\alpha-4} \vex^2 \nn\\
    &\le& \frac{1+a_1}{2\sqrt{3}} \cdot \sqrt{\eps} \io \ue
    + \frac{A}{2\sqrt{3}} \cdot \sqrt{\eps} \io \ve
    \qquad \mbox{for all $t>0$ and } \eps\in (0,1).
  \eea
\end{lem}
\proof
  We take an appropriate linear combination of the inequalities and identities provided by Lemma \ref{lem5} and
  Lemma \ref{lem22} i) and ii) to see that for all $\eps\in (0,1)$,
  \bea{23.3}
    \hspace*{-10mm}
    \frac{d}{dt} \Eo(t)
    &+& \us \eps \io \uexx^2
    + \us \eps^\frac{\alpha}{2} \io \ue^{-\alpha-2} \uex^2
    + \us \eps^\frac{\alpha+2}{2} \io \ue^{-\alpha-4} \uex^2 \nn\\
    &+& D_1 \us \io \frac{\uex^2}{\ue^2}
    + D_1 \us \eps \io \frac{\uex^2}{\ue^4} \nn\\
    &+& A\vs  \eps \io \vexx^2
    + A\vs \eps^\frac{\alpha}{2} \io \ve^{-\alpha-2} \vex^2
    + A\vs \eps^\frac{\alpha+2}{2} \io \ve^{-\alpha-4} \vex^2 \nn\\
    &+& AD_2 \vs \io \frac{\vex^2}{\ve^2}
    + AD_2 \vs \eps \io \frac{\vex^2}{\ve^4} \nn\\
    &\le& \Big( \lambda_1 + \frac{\sqrt{\eps}}{2\sqrt{3}} \Big) \io \ue
    - \io \ue^2
    + a_1 \io \ue\ve \nn\\
    & & + A\cdot \Big( \lambda_2 + \frac{\sqrt{\eps}}{2\sqrt{3}}\Big) \io \ve
    - A \io \ve^2
    - Aa_2 \io \ue\ve
    + \frac{A a_2 \sqrt{\eps}}{2\sqrt{3}} \io \ue \nn\\
    & & + \chi_1 \us \io \frac{1}{\ue} \uex\vex
    - \lambda_1 \us |\Omega|
    + \us \io \ue
    - a_1 \us \io \ve \nn\\
    & & - A\chi_2 \vs \io \frac{1}{\ve} \uex\vex
    - A\lambda_2 \vs |\Omega|
    + A\vs \io \ve
    + A a_2 \vs \io \ue
    \qquad \mbox{for all } t>0.
  \eea
  Here using Young's inequality, we obtain that
  \bas
    \chi_1 \us \io \frac{1}{\ue} \uex\vex
    &\le& \frac{D_1\us}{2} \io \frac{\uex^2}{\ue^2}
    + \frac{\chi_1^2 \us}{2D_1} \io \vex^2 \\
    &\le& \frac{D_1\us}{2} \io \frac{\uex^2}{\ue^2}
    + \frac{\chi_1^2 \us}{2D_1} \|\ve\|_{L^\infty(\Omega)}^2 \io \frac{\vex^2}{\ve^2}
  \eas
  and, similarly,
  \bas
    -A\chi_2\vs \io \frac{1}{\ve} \uex\vex
    \le \frac{AD_2 \vs}{2} \io \frac{\vex^2}{\ve^2}
    + \frac{A\chi_2^2\ve}{2D_2} \|\ue\|_{L^\infty(\Omega)}^2 \io \frac{\uex^2}{\ue^2}
  \eas
  for all $t>0$.
  Therefore, after neglecting some signed summands and rearranging we infer from (\ref{23.3}) that
  \bea{23.4}
    \frac{d}{dt} \Eo(t)
    &+& \bigg\{ \frac{D_1 \us}{2} - \frac{A\chi_2^2 \vs}{2D_2} \|\ue\|_{L^\infty(\Omega)}^2 \bigg\}
    \cdot \io \frac{\uex^2}{\ue^2} \nn\\
    &+& \bigg\{ \frac{AD_2 \vs}{2} - \frac{\chi_1^2 \us}{2D_1} \|\ve\|_{L^\infty(\Omega)}^2 \bigg\}
    \cdot \io \frac{\vex^2}{\ve^2} \nn\\
    &+& \us \eps^\frac{\alpha+2}{2} \io \ue^{-\alpha-4}\uex^2
        +A\vs \eps^\frac{\alpha+2}{2} \io \ve^{-\alpha-4}\vex^2
        \nn\\
    &\le& I_\eps(t):=
    - \io \ue^2
    + (\lambda_1 + \us + Aa_2 \vs) \io \ue
    - \lambda_1 \us |\Omega| \nn\\
    & & - A\io \ve^2
    + (A\lambda_2 - a_1\us + A\vs) \io \ve
    - A \lambda_2 \vs |\Omega| \nn\\
    & & + (a_1 -Aa_2) \io \ue\ve \nn\\
    & & + \frac{1+Aa_2}{2\sqrt{3}} \cdot \sqrt{\eps} \io \ue
    + \frac{A}{2\sqrt{3}} \cdot \sqrt{\eps} \io \ve
    \qquad \mbox{for all } t>0.
  \eea
  We now make use of our definition (\ref{A}) of $A$, which along with (\ref{usvs}) namely ensures the precise identities
  \bas
    \lambda_1 + \us + Aa_2\vs=\lambda_1+\us+a_1\vs=2\us
    \quad \mbox{and} \quad
    A\lambda_2 - a_1\us + A\vs=2A\vs
  \eas
  as well as
  \bas
    \us^2 + A\vs^2 -\lambda_1 \us -A\lambda_2\vs
    = \frac{a_2\us(\us-\lambda_1) + a_1\vs (\vs-\lambda_2)}{a_2}
    = \frac{a_2 \us \cdot a_1\vs + a_1\vs \cdot (-a_2\us)}{a_2}
    =0
  \eas
  and clearly also $a_1-Aa_2=0$. Accordingly, in (\ref{23.4}) we may simplify
  \bas
    I_\eps(t)
    = - \io (\ue-\us)^2 - A\io (\ve-\vs)^2
    + \frac{1+a_1}{2\sqrt{3}} \cdot \sqrt{\eps} \io \ue
    + \frac{A}{2\sqrt{3}} \cdot \sqrt{\eps} \io \ve
  \eas
  for all $t>0$ and $\eps\in (0,1)$, and thus end up with (\ref{23.2}).
\qed
In view of the eventual bounds derived in Section \ref{sect_evbounds}, for small values of $\chi_1$ and $\chi_2$
but arbitrarily large $u_0$ and $v_0$ this readily entails an inequality for $\Eo$ which in the formal limit
$\eps\searrow 0$ indeed predicts eventual decrease thereof.
\begin{lem}\label{lem25}
  Let $n_1=n_2=2$, and for $i\in \{1,2\}$ let
  $D_i>0, a_i>0, \lambda_i>0$ and $\chi_i>0$ be such that (\ref{23.1}) holds.
  Then with $\chis>0$ taken from Lemma \ref{lem62}, one can find $\chiss\in (0,\chis)$ and $C>0$ such that
  if $\chi_1\in (0,\chiss)$ and $\chi_2\in (0,\chiss)$, and if (\ref{ie}) holds,
  there exist $T_0=T_0(\chi_1,\chi_2,u_0,v_0)>0$ and $\eps_0=\eps_0(\chi_1,\chi_2,u_0,v_0)\in (0,1)$
  such that for all 
  $\eps\in (0,\eps_0)$, the function $\Eo$ defined in (\ref{Eo}) satisfies
  \be{25.1}
    \frac{d}{dt} \Eo(t) + \frac{1}{C} \Do(t) \le C\sqrt{\eps}
    \qquad \mbox{for all } t>T_0,
  \ee
  where for $\eps\in (0,1)$,
  \bea{Do}
    \Do(t)
    &:=&
    \io \frac{\uex^2(\cdot,t)}{\ue^2(\cdot,t)}
    + \io \frac{\vex^2(\cdot,t)}{\ve^2(\cdot,t)}
    + \io (\ue(\cdot,t)-\us)^2
    + \io (\ve(\cdot,t)-\vs)^2 \nn\\
    &+& \eps^\frac{\alpha+2}{2} \io \ue^{-\alpha-4}(\cdot,t) \uex^2(\cdot,t)
    + \eps^\frac{\alpha+2}{2} \io \ve^{-\alpha-4}(\cdot,t) \vex^2(\cdot,t),
    \qquad t>0,
  \eea
  and where     
  $\us>0$ and $\vs>0$ are taken from 
  (\ref{usvs}).
\end{lem}
\proof
  We first invoke Lemma \ref{lem62} to find $c_1>0$ such that if $\chi_1\in (0,\chis)$ and $\chi_2\in (0,\chis)$, then
  whenever (\ref{ie}) holds, we can find $T_0=T_0(\chi_1,\chi_2,u_0,v_0)>0$ and
  $\eps_0=\eps_0(\chi_1,\chi_2,u_0,v_0)\in (0,1)$ such that
  \be{25.2}
    \|\ue(\cdot,t)\|_{L^\infty(\Omega)}
    + \|\ve(\cdot,t)\|_{L^\infty(\Omega)} \le c_1
    \qquad \mbox{for all $t>T_0$ and any } \eps\in (0,\eps_0).
  \ee
  Choosing $\chiss\in (0,\chis)$ small enough such that with $A>0$ given by (\ref{A}) we have
  \bas
    \frac{A\chiss^2 \vs}{2D_2} \cdot c_1^2 \le \frac{D_1 \us}{4}
    \qquad \mbox{and} \qquad
    \frac{\chiss^2 \us}{2D_1} \cdot c_1^2
    \le \frac{AD_2 \vs}{4},
  \eas
  from Lemma \ref{lem23} we infer that if $\chi_1\in (0,\chiss)$, $\chi_2\in (0,\chiss)$ and (\ref{ie})
  holds, then with $T_0$ and $\eps_0$ as above and for all $\eps\in (0,\eps_0)$,
  \bas
    \frac{d}{dt} \Eo(t)
    &+& \frac{D_1\us}{4} \io \frac{\uex^2}{\ue^2}
    + \frac{AD_2\vs}{4} \io \frac{\vex^2}{\ve^2}
    + \io (\ue-\us)^2 + A\io (\ve-\vs)^2 \nn\\
    &+& \us\eps^\frac{\alpha+2}{2} \io \ue^{-\alpha-4} \uex^2
    + A\vs \eps^\frac{\alpha+2}{2} \io \ve^{-\alpha-4} \vex^2 \nn\\
    &\le& \Big\{ \frac{(1+a_1) c_1 |\Omega|}{2\sqrt{3}} + \frac{Ac_1 |\Omega|}{2\sqrt{3}} \Big\} \cdot \sqrt{\eps}
    \qquad \mbox{for all } t>T_0,
  \eas
  which directly yields (\ref{25.1}).
\qed
Now once more thanks to the interpolation properties from Lemma \ref{lem14}, the dissipation rate in (\ref{25.1}) dominates
a certain superlinear power of $\Eo$ at least within bounded time intervals:
\begin{lem}\label{lem27}
  Let $n_1=n_2=2$, and let
  $D_i>0, a_i>0, \lambda_i>0$ and $\chi_i>0$ for $i\in \{1,2\}$, assuming (\ref{23.1}).
  Then whenever (\ref{ie}) holds, for all $T>0$ there exists $C(T)>0$ such that for all $\eps\in (0,1)$, with $\Eo$ and $\Do$
  taken from (\ref{Eo}) and (\ref{Do}) we have
  \be{27.1}
    \Eo^\frac{\alpha+2}{2}(t)
    \le C(T) \Do(t) + C(T)
    \qquad \mbox{for all $t\in (0,T)$.}
  \ee
\end{lem}
\proof
  According to Lemma \ref{lem18}, given $T>0$ we can find $c_1(T) \in (0,1)$ such that for all $\eps\in (0,1)$,
  \be{27.2}
    \io \ue(\cdot,t) \ge c_1(T)
    \quad \mbox{and} \quad
    \io \ve(\cdot,t) \ge c_1(T)
    \qquad \mbox{for all } t\in (0,T),
  \ee
  whereas Lemma \ref{lem3} provides $c_2>0$ such that
  \be{27.3}
    \io \ue(\cdot,t) + \io \ve(\cdot,t) \le c_2
    \qquad \mbox{for all } t>0
  \ee
  whenever $\eps\in (0,1)$.
  We now recall the definition (\ref{Eo}) of $\Eo$ to firstly estimate
  \bas
    \Eo(t)
    &\le& \io \ue - \us \io \ln \ue + \us \ln \us \cdot |\Omega| + \frac{\us \eps}{6} \io \frac{1}{\ue^2} \\
    & & + A \io \ve - A\vs \io \ln \ve + A\vs\ln\vs \cdot |\Omega| + \frac{A\vs \eps}{6} \io \frac{1}{\ve^2}
    \qquad \mbox{for all } t>0,
  \eas
  from which by nonnegativity of $\Eo$, as asserted by Lemma \ref{lem29}, for all $t>0$ it follows that due to (\ref{27.3}),
  \bea{27.4}
    \hspace*{-10mm}
    \Eo^\frac{\alpha+2}{2}(t)
    &\le& 8^\frac{\alpha+2}{2} \cdot \Bigg\{
    \us^\frac{\alpha+2}{2} \cdot \bigg\{ -\io \ln\ue \bigg\}_+^\frac{\alpha+2}{2}
    + \Big(\frac{\us}{6}\Big)^\frac{\alpha+2}{2} \eps^\frac{\alpha+2}{2} \cdot
        \bigg\{ \io \frac{1}{\ue^2} \bigg\}^\frac{\alpha+2}{2} \nn\\
    & & \hspace*{20mm}
    + (A\vs)^\frac{\alpha+2}{2} \cdot \bigg\{ -\io \ln\ve \bigg\}_+^\frac{\alpha+2}{2}
    + \Big(\frac{A\vs}{6}\Big)^\frac{\alpha+2}{2} \eps^\frac{\alpha+2}{2} \cdot
        \bigg\{ \io \frac{1}{\ve^2} \bigg\}^\frac{\alpha+2}{2}
    + c_3 \Bigg\}
  \eea
  with $c_3:=c_2^\frac{\alpha+2}{2} + |\us\ln\us \cdot |\Omega||^\frac{\alpha+2}{2}
  +(Ac_2)^\frac{\alpha+2}{2} + |A\vs\ln \vs \cdot |\Omega||^\frac{\alpha+2}{2}$.
  Here we use Lemma \ref{lem14} ii) along with (\ref{27.2}) and Young's inequality to see that writing
  $c_4(T):=|\Omega| \cdot \ln \frac{1}{c_1(T)} + |\Omega| \cdot |\ln |\Omega||$ we have
  \bea{27.5}
    \hspace*{-10mm}
    \bigg\{ - \io \ln\ue \bigg\}_+^\frac{\alpha+2}{2}
    &\le& \Bigg\{ |\Omega|^\frac{3}{2} \cdot \bigg\{ \io \frac{\uex^2}{\ue^2} \bigg\}^\frac{1}{2}
    - |\Omega| \ln \bigg\{ \io \ue \bigg\} + |\Omega| \ln |\Omega| \Bigg\}^\frac{\alpha+2}{2} \nn\\
    &\le& \Bigg\{ |\Omega|^\frac{3}{2} \cdot \bigg\{ \io \frac{\uex^2}{\ue^2} \bigg\}^\frac{1}{2} + c_4
        \Bigg\}^\frac{\alpha+2}{2} \nn\\
    &\le& \Big( 2|\Omega|^\frac{3}{2}\Big)^\frac{\alpha+2}{2} \cdot
    \bigg\{ \io \frac{\uex^2}{\ue^2} \bigg\}^\frac{\alpha+2}{4}
    + (2c_4)^\frac{\alpha+2}{2} \nn\\
    &\le& \Big( 2|\Omega|^\frac{3}{2}\Big)^\frac{\alpha+2}{2} \cdot \io \frac{\uex^2}{\ue^2}
    + \Big( 2|\Omega|^\frac{3}{2}\Big)^\frac{\alpha+2}{2}
    + (2c_4)^\frac{\alpha+2}{2}
    \qquad \mbox{for all } t\in (0,T).
  \eea
  Moreover, applying Lemma \ref{lem14} i) to $p=2$ and $q=\alpha+2$ shows that again thanks to (\ref{27.2}) and
  Young's inequality,
  \bas
    \eps^\frac{\alpha+2}{2} \cdot \bigg\{ \io \frac{1}{\ue^2} \bigg\}^\frac{\alpha+2}{2}
    &\le& \eps^\frac{\alpha+2}{2} \cdot \Bigg\{
    (\alpha+2)^\frac{4}{\alpha+2} |\Omega|^\frac{\alpha+4}{\alpha+2} \cdot
    \bigg\{ \io \ue^{-\alpha-4} \uex^2 \bigg\}^\frac{2}{\alpha+2} \nn\\
    & & \hspace*{20mm}
    + 2^\frac{4}{\alpha+2} |\Omega|^3 \cdot \bigg\{ \io \ue \bigg\}^{-2} \Bigg\}^\frac{\alpha+2}{2} \nn\\
    &\le& 2^\frac{\alpha}{2} (\alpha+2)^2 |\Omega|^\frac{\alpha+4}{2} \cdot \eps^\frac{\alpha+2}{2}
    \io \ue^{-\alpha-4} \uex^2
    + 2^\frac{\alpha+4}{2} |\Omega|^\frac{3(\alpha+2)}{2} c_1^{-\alpha-2}
  \eas
  for all $t\in (0,T)$.
  Together with a similar argument for the respective expressions containing the second solution component,
  this enables us to infer from (\ref{27.4}) the existence of $c_5>0$ and $c_6(T)>0$ such that for all $\eps\in (0,1)$,
  \bas
    \Eo^\frac{\alpha+2}{2}(t)
    \le c_5 \cdot \bigg\{ \io \frac{\uex^2}{\ue^2} + \io \frac{\vex^2}{\ve^2}
    + \eps^\frac{\alpha+2}{2} \io \ue^{-\alpha-4} \uex^2 + \eps^\frac{\alpha+2}{2} \io \ve^{-\alpha-4} \vex^2 \bigg\}
    + c_6(T)
    \quad \mbox{for all } t\in (0,T).
  \eas
  In view of the definition (\ref{Do}) of $\Do$, this implies (\ref{27.1}).
\qed
For drawing conclusions from the resulting autonomous ODI for $\Eo$, let us a quantitative consequence of a simple
comparison argument.
\begin{lem}\label{lem28}
  Let $t_0\in\R, T>t_0, a>0, b>0$ and $\beta>1$, and suppose that $y\in C^0([t_0,T)) \cap C^1((t_0,T))$
  is nonnegative and such that
  \be{28.1}
    y'(t) + ay^\beta(t) \le b
    \qquad \mbox{for all } t\in (t_0,T).
  \ee
  Then
  \be{28.2}
    y(t) \le \Big\{ (\beta-1) a (t-t_0) \Big\}^{-\frac{1}{\beta-1}} + \Big(\frac{b}{a}\Big)^\frac{1}{\beta}
    \qquad \mbox{for all } t\in (t_0,T).
  \ee
\end{lem}
\proof
  Without loss of generality we assume that $t_0=0$, and observe that
  $\oy(t):=((\beta-1) at)^{-\frac{1}{\beta-1}} + (\frac{b}{a})^\frac{1}{\beta}$, $t>0$, satisfies
  \bas
    \oy'(t) + a\oy^\beta(t) - b
    &=& (\beta-1)^{-\frac{1}{\beta-1}} a^{-\frac{1}{\beta-1}} \cdot \Big(-\frac{1}{\beta-1} t^{-\frac{1}{\beta-1}-1} \Big)
    + a\cdot \Big\{ \Big( (\beta-1) at\Big)^{-\frac{1}{\beta-1}} + \Big(\frac{b}{a}\Big)^\frac{1}{\beta} \Big\}^\beta
    - b \\
    &\ge& - (\beta-1)^{-\frac{\beta}{\beta-1}} a^{-\frac{1}{\beta-1}} t^{-\frac{\beta}{\beta-1}}
    + a\cdot \Big((\beta-1) at \Big)^{-\frac{\beta}{\beta-1}}
    + a\cdot \frac{b}{a} - b =0
  \eas
  for all $t>0$
  due to the fact that $(\xi+\eta)^\beta \ge \xi^\beta+\eta^\beta$ for all $\xi\ge 0$ and $\eta\ge 0$.
  Since $\oy(t) \nearrow + \infty$ as $t\searrow 0$, by continuity of $y$ at $t=0$ this readily implies the claimed
  inequality by means of a comparison argument.
\qed
We can thus reap the fruit of Lemma \ref{lem25} and Lemma \ref{lem27}
and thereby obtain, through the dissipation mechanism expressed in (\ref{25.1}),
the following preliminary decay information for our solution to the original problem.
\begin{lem}\label{lem26}
  Let $n_1=n_2=2$, let
  $D_i>0, a_i>0, \lambda_i>0$ and $\chi_i>0$ for $i\in \{1,2\}$ be such that (\ref{23.1}) holds,
  and let $\chi_1\in (0,\chiss)$ and $\chi_2\in (0,\chiss)$ with $\chiss>0$ as in Lemma \ref{lem25}.
  Then assuming (\ref{ie}), one can find $T_0=T_0(\chi_1,\chi_2,u_0,v_0)>0$ such that
  for the limit functions $u$ and $v$ obtained in Lemma \ref{lem_conv} we have
  \be{26.1}
    \int_{T_0}^\infty \io u_x^2
    + \int_{T_0}^\infty \io v_x^2
    + \int_{T_0}^\infty \io (u-\us)^2
    + \int_{T_0}^\infty \io (v-\vs)^2
    <\infty
  \ee
  with $\us>0$ and $\vs>0$ taken from (\ref{usvs}).
\end{lem}
\proof
  According to Lemma \ref{lem25} and Lemma \ref{lem62}, supposing (\ref{ie}) to be valid we
  can pick $\eps_0=\eps_0(\chi_1,\chi_2,u_0,v_0)\in (0,1)$, $T_1=T_1(\chi_1,\chi_2,u_0,v_0)>0$
  and $c_i=c_i(\chi_1,\chi_2,u_0,v_0)>0$, $i\in \{1,2,3\}$, such that
  \be{26.2}
    \frac{d}{dt} \Eo(t) + c_1 \Do(t) \le c_2 \sqrt{\eps}
    \qquad \mbox{for all } t>T_1
  \ee
  as well as
  \be{26.22}
    \|\ue(\cdot,t)\|_{L^\infty(\Omega)}
    +\|\ve(\cdot,t)\|_{L^\infty(\Omega)}
    \le c_3
    \qquad \mbox{for all $t>T_1$ and any } \eps\in (0,\eps_0),
  \ee
  where $\Eo$ and $\Do$ are as defined in (\ref{Eo}) and (\ref{Do}), respectively.
  Keeping this value of $T_1$ fixed, we thereafter invoke Lemma \ref{lem27} to pick $c_4>0$
  such that for any choice of $\eps\in (0,1)$,
  \bas
    \Eo^\frac{\alpha+2}{2}(t) \le c_4 \Do(t) + c_4
    \qquad \mbox{for all } t\in (0,T_1+2),
  \eas
  which when combined with (\ref{26.2}) shows that
  \bas
    \frac{d}{dt} \Eo(t)
    + \frac{c_1}{c_4} \Eo^\frac{\alpha+2}{2}(t)
    \le c_1 + c_2\sqrt{\eps}
    \le c_1+c_2
    \qquad \mbox{for all $t\in (T_1,T_1+2)$ and all } \eps\in (0,\eps_0).
  \eas
  Through Lemma \ref{lem28}, the latter implies that
  \bas
    \Eo(t) \le \Big\{ \frac{\alpha}{2} \cdot \frac{c_1}{c_4} \cdot (t-T_1) \Big\}^{-\frac{2}{\alpha}}
    + \bigg( \frac{ c_1+c_2}{\frac{c_1}{c_4}} \bigg)^\frac{2}{\alpha+2}
    \qquad \mbox{for all $t\in (T_1,T_1+2)$ and } \eps\in (0,\eps_0)
  \eas
  and that hence, in particular,
  \bas
    \Eo(T_1+1)
    \le c_5:=\Big( \frac{2c_4}{c_1 \alpha} \Big)^\frac{2}{\alpha}
    + \Big( \frac{(c_1+c_2)c_4}{c_1} \Big)^\frac{2}{\alpha+2}
    \qquad \mbox{for all } \eps\in (0,\eps_0).
  \eas
  Using this as $\eps$-independent information at the initial time $T_1+1$,
  we now return to (\ref{26.2}) to infer upon an integration therein that
  \bea{26.4}
    \hspace*{-10mm}
    c_1 \int_{T_1+1}^T \Do(t) dt
    &\le& \Eo(T_1+1) - \Eo(T) + c_2\sqrt{\eps} \cdot (T-T_1-1) \nn\\
    &\le& c_5 + c_2\sqrt{\eps} \cdot (T-T_1-1)
    \qquad \mbox{for all $T>T_1+1$ and each } \eps\in (0,\eps_0),
  \eea
  again because $\Eo$ is nonnegative.
  Since (\ref{26.22}) ensures that according to (\ref{Do}) we can estimate
  \bas
    \Do(t)
    \ge \frac{1}{c_3^2} \io \uex^2 + \frac{1}{c_3^2} \io \vex^2
    + \io (\ue-\us)^2 + \io (\ve-\vs)^2
    \qquad \mbox{for all $t>T_1$ and } \eps\in (0,\eps_0),
  \eas
  on taking $\eps=\eps_j\searrow 0$ with $(\eps_j)_{j\in\N}\subset (0,1)$ as in Lemma \ref{lem_conv},
  from (\ref{26.4}) and an argument based on lower semicontinuity of norms with respect to weak convergence in Hilbert
  spaces we conclude that
  \bas
    \frac{1}{c_3^2} \int_{T_1+1}^T \io u_x^2
    + \frac{1}{c_3^2} \int_{T_1+1}^T \io v_x^2
    + \int_{T_1+1}^T \io (u-\us)^2
    + \int_{T_1+1}^T \io (v-\vs)^2
    \le \frac{c_5}{c_1}
    \qquad \mbox{for all } t>T_1+1,
  \eas
  and that thus (\ref{26.1}) holds with $T_0:=T_1+1$.
\qed
We next make use of favorable smallness properties, as implied by (\ref{26.1}) for certain arbitrarily large times
which we use as new starting instants, to see on going back to (\ref{25.1}), but this time simply neglecting the positive summand $\Do$ therein, that $u$ and $v$ in fact approach their expected limits in a sense much stronger
than indicated in Lemma \ref{lem26}, albeit not yet identified as spatially uniform but rather in a topology
associated with the entropy functional in (\ref{entropy2}).
\begin{lem}\label{lem31}
  Let $n_1=n_2=2$, let
  $D_i>0, a_i>0, \lambda_i>0$ and $\chi_i>0$, $i\in \{1,2\}$, satisfy (\ref{23.1}),
  and let $\chi_1\in (0,\chiss)$ and $\chi_2\in (0,\chiss)$ with $\chiss>0$ taken from Lemma \ref{lem25}.
  Then assuming (\ref{ie}) and letting $u,v$ and $N$ be as given by Lemma \ref{lem_conv},
  with $\phi_{\us}$ and $\phi_{\vs}$ taken from (\ref{phi}) and
  (\ref{usvs}) we have
  \be{31.1}
    \io \phi_{\us}(u(\cdot,t))\to 0
    \quad \mbox{and} \quad
    \io \phi_{\vs}(v(\cdot,t)) \to 0
    \qquad \mbox{as } (0,\infty)\setminus N \ni t \to\infty.
  \ee
\end{lem}
\proof
  In order to prepare our convergence argument, let us first once more resort to Lemma \ref{lem25} in choosing
  $c_1>0$, $T_0=T_0(u_0,v_0)>0$ and $\eps_\star=\eps_\star(u_0,v_0)\in (0,1)$ such that for all $\eps\in (0,\eps_\star)$,
  \be{31.2}
    \frac{d}{dt} \Eo(t) \le c_1 \sqrt{\eps}
    \qquad \mbox{for all } t>T_0.
  \ee
  Then abbreviating $c_2:=|\Omega|^\frac{1}{2} + |\Omega|^{-\frac{1}{2}}$, given $\eta>0$ we fix $\delta>0$ conveniently
  small such that
  \be{31.3}
    2c_2\sqrt{\delta} \le \frac{\us}{2}
    \qquad \mbox{and} \qquad
    2c_2\sqrt{\delta} \le \frac{\vs}{2}
  \ee
  as well as
  \be{31.4}
    \delta \le \frac{\us \eta}{64 c_2^2 |\Omega|}
    \qquad \mbox{and} \qquad
    A \delta \le \frac{\vs \eta}{64 c_2^2 |\Omega|}.
  \ee
  According to Lemma \ref{lem26}, we can thereafter pick some $T_1>T_0+1$ suitably large such that
  \be{31.5}
    \int_{T_1-1}^\infty \io u_x^2
    + \int_{T_1-1}^\infty \io v_x^2
    + \int_{T_1-1}^\infty \io (u-\us)^2
    + \int_{T_1-1}^\infty \io (v-\vs)^2
    \le \delta,
  \ee
  and we claim that then
  \be{31.6}
    \io \phi_{\us}(u(\cdot,t)) + A \io \phi_{\vs}(v(\cdot,t))
    \le \eta
    \qquad \mbox{for all $t>T_1$ such that } t\not\in N.
  \ee
  To verify this, given any such $t$ we first observe that due to (\ref{31.5}),
  \bas
    \int_{t-1}^t \bigg\{ \io u_x^2 + \io v_x^2 + \io (u-\us)^2 + \io (v-\vs)^2 \bigg\} \le \delta,
  \eas
  whence it is possible to find $t_\star=t_\star(t) \in (t-1,t)\setminus N$ such that
  \be{31.7}
    \io u_x^2(\cdot,t_\star) + \io v_x^2(\cdot,t_\star)
    + \io (u(\cdot,t_\star)-\us)^2
    + \io (v(\cdot,t_\star)-\vs)^2
    \le \delta.
  \ee
  In particular, this entails the existence of $x_\star\in \bom$ such that $(u(x_\star,t_\star)-\us)^2
  \le \frac{\delta}{|\Omega|}$ and that hence
  \bas
    |u(x,t_\star)-\us|
    &\le& |u(x,t_\star)-u(x_\star,t_\star)|
    + |u(x_\star,t_\star) -\us| \\
    &=& \bigg| \int_{x_\star}^x u_x(y,t_\star) dy \bigg|
    + |u(x_\star,t_\star) -\us| \\
    &\le& |\Omega|^\frac{1}{2} \cdot \bigg\{ \io u_x^2(y,t_\star) dy \bigg\}^\frac{1}{2}
    + |u(x_\star,t_\star) -\us| \\
    &\le& |\Omega|^\frac{1}{2} \delta^\frac{1}{2}
    + \Big( \frac{\delta}{|\Omega|} \Big)^\frac{1}{2}
    \qquad \mbox{for all } x\in\Omega,
  \eas
  which along with an identical argument for $v$ shows that
  \be{31.8}
    \|u(\cdot,t_\star)-\us\|_{L^\infty(\Omega)} \le c_2 \sqrt{\delta}
    \qquad \mbox{and} \qquad
    \|v(\cdot,t_\star)-\vs\|_{L^\infty(\Omega)} \le c_2 \sqrt{\delta}.
  \ee
  Now since $t_\star \in (0,\infty)\setminus N$, Corollary \ref{cor24} applies so as to warrant that with some
  $\eps_{\star\star}=\eps_{\star\star}(t_\star)\in (0,\eps_\star)$ and $(\eps_j)_{j\in\N} \subset (0,1)$ as
  in Lemma \ref{lem_conv} we have
  \bas
    \|\ue(\cdot,t_\star)-u(\cdot,t_\star)\|_{L^\infty(\Omega)}
    + \|\ve(\cdot,t_\star)-v(\cdot,t_\star)\|_{L^\infty(\Omega)}
    \le c_2\sqrt{\delta}
    \qquad \mbox{for all } \eps\in (\eps_j)_{j\in\N} \cap (0,\eps_{\star\star}),
  \eas
  together with (\ref{31.8}) ensuring that
  \be{31.9}
    \|\ue(\cdot,t_\star)-\us\|_{L^\infty(\Omega)} \le 2c_2\sqrt{\delta}
    \quad \mbox{and} \quad
    \|\ve(\cdot,t_\star)-\vs\|_{L^\infty(\Omega)} \le 2c_2\sqrt{\delta}
    \qquad \mbox{for all } \eps\in (\eps_j)_{j\in\N} \cap (0,\eps_{\star\star}).
  \ee
  According to our requirements on $\delta$ in (\ref{31.3}), these estimates especially entail the inequalities
  \be{31.99}
    \ue(\cdot,t_\star) \ge \us - 2c_2\sqrt{\delta}
    \ge \frac{\us}{2}
    \quad \mbox{and} \quad
    \ve(\cdot,t_\star) \ge \vs - 2c_2\sqrt{\delta}
    \ge \frac{\vs}{2}
    \quad \mbox{in $\Omega$ \qquad for all } \eps\in (\eps_j)_{j\in\N} \cap (0,\eps_{\star\star}),
  \ee
  which firstly enables us to conclude from (\ref{29.2}) in conjunction with (\ref{31.9}) and (\ref{31.4}) that
  for all $\eps\in (\eps_j)_{j\in\N} \cap (0,\eps_{\star\star})$,
  \bea{31.10}
    \io \phi_{\us}(\ue(\cdot,t_\star))
    \le \frac{2}{\us} \io \Big(\ue(\cdot,t_\star)-\us\Big)^2
    \le \frac{2|\Omega|}{\us} \cdot \|\ue(\cdot,t_\star)-\us\|_{L^\infty(\Omega)}^2
    \le \frac{2|\Omega|}{\us} \cdot (2c_2\sqrt{\delta})^2
    \le \frac{\eta}{8}  \eea
  and similarly
  \bea{31.11}
   A \io \phi_{\vs}(\ve(\cdot,t_\star)) \le \frac{\eta}{8}.
  \eea
  Secondly, (\ref{31.99}) guarantees that if we pick $\eps_{\star\star\star}=\eps_{\star\star\star}(t_\star)
  \in (0,\eps_{\star\star})$ small enough such that
  \bas
    \frac{2|\Omega|\eps_{\star\star\star}}{3\us} \le \frac{\eta}{8}
    \qquad \mbox{and} \qquad
    \frac{2A|\Omega|\eps_{\star\star\star}}{3\vs} \le \frac{\eta}{8},
  \eas
  then in the contributions to $\Eo$ containing the factor $\eps$ we can estimate
  \bas
    \frac{\us\eps}{6} \io \frac{1}{\ue^2(\cdot,t_\star)}
    \le \frac{\us\eps}{6} \cdot \Big(\frac{2}{\us}\Big)^2 |\Omega|
    \le \frac{\eta}{8}
    \qquad \mbox{for all } \eps\in (\eps_j)_{j\in\N} \cap (0,\eps_{\star\star\star})
  \eas
  and
  \bas
    \frac{A\vs\eps}{6} \io \frac{1}{\ve^2(\cdot,t_\star)}
    \le \frac{A\vs\eps}{6} \cdot \Big(\frac{2}{\vs}\Big)^2 |\Omega|
    \le \frac{\eta}{8}
    \qquad \mbox{for all } \eps\in (\eps_j)_{j\in\N} \cap (0,\eps_{\star\star\star}).
  \eas
  When combined with (\ref{31.10}) and (\ref{31.11}), in view of (\ref{Eo}) these inequalities show that
  \bas
    \Eo(t_\star) \le 4\cdot \frac{\eta}{8}=\frac{\eta}{2}
    \qquad \mbox{for all } \eps\in (\eps_j)_{j\in\N} \cap (0,\eps_{\star\star\star}),
  \eas
  so that letting $\eps_{\star\star\star\star}=\eps_{\star\star\star\star}(t_\star) \in (0,\eps_{\star\star\star})$
  be such that $c_1 \sqrt{\eps_{\star\star\star\star}} \le \frac{\eta}{2}$,
  on integrating (\ref{31.2}) we infer that at the time in question we have
  \bea{31.12}
    \Eo(t)
    \le \Eo(t_\star) + c_1\sqrt{\eps} \cdot (t-t_\star)
    \le \frac{\eta}{2} + c_1\sqrt{\eps}
    \le \frac{\eta}{2} + \frac{\eta}{2}=\eta
    \qquad \mbox{for all } \eps\in (\eps_j)_{j\in\N} \cap (0,\eps_{\star\star\star\star}),
  \eea
  because $t_\star\ge t-1$.
  Since from Lemma \ref{lem_conv} and our assumption that $t\not\in N$
  we know that $\ue(\cdot,t) \to u(\cdot,t)$ and $\ve(\cdot,t) \to v(\cdot,t)$ a.e.~in $\Omega$ as $\eps=\eps_j\searrow 0$,
  upon an application of Fatou's lemma we readily obtain (\ref{31.6}), and thus the statement of the lemma,
  as a consequence of (\ref{31.12}).
\qed
Now thanks to the eventual precompactness features implied by Corollary \ref{cor24}, due to the positivity
of $\phi_{\xis}$ outside the point $\xis>0$ the latter readily implies the desired convergence statement.
\begin{lem}\label{lem32}
  Let $n_1=n_2=2$, let
  $D_i>0, a_i>0, \lambda_i>0$ and $\chi_i>0$ for $i\in \{1,2\}$ be such that (\ref{23.1}) is fulfilled,
  and let $\chi_1\in (0,\chiss)$ and $\chi_2\in (0,\chiss)$ with $\chiss>0$ as in Lemma \ref{lem25}.
  Then whenever (\ref{ie}) holds, as $t\to\infty$, the limit functions $u$ and $v$ obtained in Lemma \ref{lem_conv} satisfy
  \be{32.1}
    u(\cdot,t) \to \us
    \quad \mbox{in } L^\infty(\Omega)
    \qquad \mbox{and} \qquad
    v(\cdot,t) \to \vs
    \quad \mbox{in } L^\infty(\Omega),
  \ee
  where $\us>0$ and $\vs>0$ are as in (\ref{usvs}).
\end{lem}
\proof
  From Corollary \ref{cor24} we know that there exists $T_0>0$ such that
  $(u(\cdot,t))_{t>T_0}$ is bounded in $W^{1,2}(\Omega)$ and hence relatively compact in $C^0(\bom)$,
  whence if (\ref{32.1}) was false, the we could find $(t_k)_{k\in\N} \subset (T_0,\infty)$
  and $u_\infty \in C^0(\bom)$ such that $u_\infty \not\equiv \us$ and that $t_k\to\infty$ as well as
  $u(\cdot,t_k) \to u_\infty$ in $L^\infty(\Omega)$ as $k\to\infty$.
  As in view of Corollary \ref{cor24} we may assume $u$ to be continuous in $\bom\times (T_0,\infty)$,
  by density of $[t_k,t_k+1] \setminus N$ in $[t_k,t_k+1]$ we can pick $\widehat{t}_k \in [t_k, t_k+1]$
  such that $\|u(\cdot,\widehat{t}_k)-u(\cdot,t_k)\|_{L^\infty(\Omega)} \le \frac{1}{k}$, meaning that also
  $u(\cdot,\widehat{t}_k) \to u_\infty$ in $L^\infty(\Omega)$ as $k\to\infty$.
  As the function $\phi_{\us}$ from (\ref{phi}) is positive in $(0,\infty)\setminus \{\us\}$ by Lemma \ref{lem29},
  however, the hypothesis $u_\infty \not\equiv \us$ implies that $\io \phi_{\us}(u(\cdot,\widehat{t}_k)) \not\to 0$
  as $k\to\infty$, in contradiction to Lemma \ref{lem31}. Along with a similar argument for $v$, this establishes
  the claim.
\qed
The proof of our main result on kinetics-driven stabilization has thereby already been accomplished:\abs
\proofc of Theorem \ref{theo33}.\quad
  We only need to combine Corollary \ref{cor24} with Lemma \ref{lem32}.
\qed

\mysection{The case $\lambda_2 \le a_2\lambda_1$. Proof of Theorem \ref{theo333}}\label{sect_asy2}
In the context addressed in Theorem \ref{theo333}, in view of Lemma \ref{lem22} the form of the functional
in (\ref{entropy3}) suggests to choose still $n_1=2$ but now rather $n_2=1$.
Our analysis will then quite closely follow the lines presented in the previous section, so that here it will
be sufficient to concentrate on the main modifications only.\abs
The fundament for convergence is constituted by a natural counterpart of Lemma \ref{lem23}.
\begin{lem}\label{lem233}
  Let $n_1=2$ and $n_2=1$, and suppose that
  $D_i>0, a_i>0, \lambda_i>0$ and $\chi_i>0$ for $i\in\{1,2\}$, and that
  \be{233.1}
    \lambda_2 \le a_2 \lambda_1.
  \ee
  Then with $A=\frac{a_1}{a_2}$ as before and $\phi_{\lambda_1}$ as determined by (\ref{phi}), for
  \bea{Et}
    \Et(t) := \io \phi_{\lambda_1}(\ue(\cdot,t))
    + \frac{\lambda_1 \eps}{6} \io \frac{1}{\ue^2(\cdot,t)}
     +A \io \ve(\cdot, t)
    &+& \frac{A}{2\lambda_2} \io \ve^2(\cdot,t)
    + \frac{A \eps}{2\lambda_2} \io
    \frac{1}{\ve(\cdot,t)},\nn\\[1mm]
    & &  \qquad t\ge 0, \ \eps\in (0,1),
  \eea
  we have
  \bea{233.2}
    \hspace*{-10mm}
    \frac{d}{dt} \Et(t)
    &+& \bigg\{ \frac{D_1 \lambda_1}{2} - \frac{A \chi_2^2}{2D_2\lambda_2} \|\ue(\cdot,t)\|_{L^\infty(\Omega)}^2
        \|\ve(\cdot,t)\|_{L^\infty(\Omega)}^2 \bigg\}
    \cdot \io \frac{\uex^2}{\ue^2} \nn\\
    &+& \bigg\{ \frac{A D_2}{2\lambda_2} - \frac{\lambda_1^2 \chi_1^2}{2D_1} \bigg\}
    \cdot \io \frac{\vex^2}{\ve^2} \nn\\
    &+& \io (\ue(\cdot,t)-\lambda_1)^2
    + \frac{A}{\lambda_2} \io \ve^3 \nn\\
    &+& \lambda_1 \eps^\frac{\alpha+2}{2} \io \ue^{-\alpha-4} \uex^2
    + \frac{A}{\lambda_2} \eps^\frac{\alpha+2}{2} \io \ve^{-\alpha-3} \vex^2 \nn\\
    &\le& \frac{1+Aa_2+ 2 Aa_2\lambda_2^{-1}}{2\sqrt{3}} \cdot \sqrt{\eps} \io \ue
    + \frac{A+2 A\lambda_2^{-1}}{2\sqrt{3}} \cdot \sqrt{\eps} \io \ve
      \eea
    for all $t>0$ and $\eps\in (0,1)$.
\end{lem}
\proof
  Combining Lemma \ref{lem5} with Lemma \ref{lem22} i) and iii), on dropping several nonnegative summands we see that
  \bea{233.3}
    & & \hspace*{-20mm}
    \frac{d}{dt} \Et(t)
    + D_1 \lambda_1 \io \frac{\uex^2}{\ue^2}
    + \lambda_1 \eps^\frac{\alpha+2}{2} \io \ue^{-\alpha-4} \uex^2 \nn\\
    & & + \frac{AD_2}{\lambda_2} \io \vex^2
    + \frac{A}{\lambda_2}\eps^\frac{\alpha+2}{2} \io \ve^{-\alpha-3} \vex^2 \nn\\
    &\le& \Big( \lambda_1 + \frac{\sqrt{\eps}}{2\sqrt{3}}\Big) \io \ue
    - \io \ue^2
    + a_1 \io \ue\ve \nn\\
    & & + \lambda_1\chi_1 \io \frac{\uex}{\ue} \vex
    - \lambda_1^2 |\Omega|
    + \lambda_1 \io \ue
    - a_1 \lambda_1 \io \ve \nn\\
    & & +A \Big(\lambda_2 +\frac{\sqrt{\eps}}{2\sqrt{3}}\Big) \io \ve -A
       \io \ve^2 -A a_2\io \ue\ve
       +\frac{Aa_2\sqrt{\eps}}{2\sqrt{3}}\io \ue \nn\\
    & & - \frac{A\chi_2}{\lambda_2} \io \ve\uex\vex
    + A \io \ve^2
    - \frac{A}{\lambda_2} \io \ve^3
    - \frac{Aa_2}{\lambda_2} \io \ue\ve^2 \nn\\
    & & - \frac{A\eps}{\lambda_2} \io \frac{2\ve^2+3\ve}{6\ve^2+2\eps} (\lambda_2-\ve-a_2\ue) \nn\\
    &=& - \io (\ue-\lambda_1)^2
    + (A\lambda_2 - a_1\lambda_1) \io \ve  -\frac{A}{\lambda_2}\io \ve^3 -\frac{Aa_2}{\lambda_2}\io \ue\ve^2\nn\\
    & & + \lambda_1\chi_1 \io \frac{\uex}{\ue} \vex
    - \frac{A\chi_2}{\lambda_2} \io \ve\uex\vex \nn\\
    & & + \frac{(1+Aa_2)\sqrt{\eps}}{2\sqrt{3}} \io \ue + \frac{A\sqrt{\eps}}{2\sqrt{3}} \io \ve\nn\\
    & & - \frac{A\eps}{\lambda_2} \io \frac{2\ve^2+3\ve}{6\ve^2+2\eps} (\lambda_2-\ve-a_2\ue)
    \qquad \mbox{for all $t>0$ and $\eps\in (0,1)$.}
  \eea
  Here we note that according to (\ref{233.1}) we have
  $A\lambda_2 - a_1\lambda_1=\frac{a_1\lambda_2}{a_2}-a_1\lambda_1 \le 0$,
  and that by Young's inequality,
  \bas
    \lambda_1\chi_1 \io \frac{\uex}{\ue}\vex
    &\le& \frac{D_1 \lambda_1}{2} \io \frac{\uex^2}{\ue^2}
    + \frac{\lambda_1^2\chi_1^2}{2D_1} \io \vex^2
  \eas
  and
  \bas
    -\frac{A\chi_2}{\lambda_2} \io \ve\uex\vex
    &\le& \frac{AD_2}{2\lambda_2} \io \vex^2
    + \frac{A\chi_2^2}{2D_2\lambda_2} \io \ve^2 \uex^2 \\
    &\le& \frac{AD_2}{2\lambda_2} \io \vex^2
    + \frac{A\chi_2^2}{2D_2\lambda_2} \|\ue\|_{L^\infty(\Omega)}^2 \|\ve\|_{L^\infty(\Omega)}^2 \io \frac{\uex^2}{\ue^2}
  \eas
  for all $t>0$ and $\eps\in (0,1)$.
  As moreover maximizing $\varphi(s):=\frac{2s^2+3s}{6s^2+2\eps}$, $s\ge 0$, shows that
  \bas
    -\frac{A\eps}{\lambda_2} \io \frac{2\ve^2+3\ve}{6\ve^2+2\eps} (\lambda_2-\ve-a_2\ue)
    &\le& \frac{A\eps}{\lambda_2} \|\varphi\|_{L^\infty((0,\infty))} \io (\ve+a_2\ue) \\
    &=& \frac{A\eps}{\lambda_2} \cdot \frac{1}{3} \cdot \bigg\{ \io \ve + a_2 \io \ue
        \bigg\}\\
    &<& \frac{A\lambda_2^{-1}}{\sqrt{3}} \cdot \sqrt{\eps} \cdot \bigg\{ \io \ve + a_2 \io \ue
        \bigg\}
  \eas
for all $t>0$ and $\eps\in (0,1)$,  from (\ref{233.3}) we directly
obtain (\ref{233.2}). \qed
This implies an inequality of the form in Lemma \ref{lem25}.
\begin{lem}\label{lem255}
  Let $n_1=2$ and $n_2=1$, let
  $D_i>0, a_i>0, \lambda_i>0$ and $\chi_i>0$ for $i\in\{1,2\}$, and suppose that (\ref{233.1}) holds.
  Then with $\chis>0$ as in Lemma \ref{lem62},
  there exists $\chiss\in (0,\chis)$ and $C>0$ such that if $\chi_1\in (0,\chiss), \chi_2\in (0,\chiss)$ and
  (\ref{ie}) is valid, then
  there exist $T_0=T_0(\chi_1,\chi_2,u_0,v_0)>0$ and $\eps_0=\eps_0(\chi_1,\chi_2,u_0,v_0)\in (0,1)$
  such that whenever
  $\eps\in (0,\eps_0)$, for $\Et$ as in (\ref{Et}) we have
  \be{255.1}
    \frac{d}{dt} \Et(t) + \frac{1}{C} \Dt(t) \le C\sqrt{\eps}
    \qquad \mbox{for all } t>T_0,
  \ee
  where
  \bea{Dt}
    \Dt(t)
    &:=&
    \io \frac{\uex^2(\cdot,t)}{\ue^2(\cdot,t)}
    + \io \vex^2(\cdot,t)
    + \io (\ue(\cdot,t)-\lambda)^2
    + \io \ve^3(\cdot,t) \nn\\
    &+& \eps^\frac{\alpha+2}{2} \io \ue^{-\alpha-4}(\cdot,t) \uex^2(\cdot,t)
    + \eps^\frac{\alpha+2}{2} \io \ve^{-\alpha-3}(\cdot,t) \vex^2(\cdot,t),
    \qquad t>0,
  \eea
  and where     
  $\us>0$ and $\vs>0$ are taken from 
  (\ref{usvs}).
\end{lem}
\proof
  Again on the basis of Lemma \ref{lem62},
  this can be derived from Lemma \ref{lem233} in much the same manner as Lemma \ref{lem25} was deduced from
  Lemma \ref{lem23}.
\qed
Up to modifications in technical details, the strategy in the proof of Lemma \ref{lem27} finds its analogue
in the following.
\begin{lem}\label{lem277}
  Let $n_1=2$ and $n_2=1$, let
  $D_i>0, a_i>0, \lambda_i>0$ and $\chi_i>0$ for $i\in\{1,2\}$, and assume (\ref{233.1}) as well as (\ref{ie}).
  Then for all $T>0$ there exists $C(T)>0$ such that for all $\eps\in (0,1)$,
  the functions $\Et$ and $\Dt$ defined in (\ref{Et}) and (\ref{Dt}) satisfy
  \be{277.1}
    \Et^\frac{\alpha+2}{2}(t)
    \le C(T) \Dt(t) + C(T)
    \qquad \mbox{for all $t\in (0,T)$.}
  \ee
\end{lem}
\proof
  We proceed in a way similar to that in Lemma \ref{lem27} and first note that
  \be{277.11}
    \inf_{\eps\in (0,1)} \inf_{t\in (0,T)} \io \ue>0,
    \quad
    \inf_{\eps\in (0,1)} \inf_{t\in (0,T)} \io \ve>0
    \quad \mbox{and} \quad
    \sup_{\eps\in (0,1)} \sup_{t\in (0,T)} \bigg\{ \io \ue + \io \ve \bigg\}<\infty,
  \ee
  to derive the existence of $c_1(T)>0$ such that for all $t\in (0,T)$ and each $\eps\in (0,1)$,
  \bea{277.2}
    \Et^\frac{\alpha+2}{2}(t)
    &\le& c_1(T) \cdot \bigg\{ - \io \ln\ue \bigg\}_+^\frac{\alpha+2}{2}
    + c_1 (T) \eps^\frac{\alpha+2}{2} \cdot \bigg\{\io \frac{1}{\ue^2} \bigg\}^\frac{\alpha+2}{2} \nn\\
    & & + c_1(T) \cdot \bigg\{ \io \ve^2 \bigg\}^\frac{\alpha+2}{2}
    + c_1(T) \eps^\frac{\alpha+2}{2} \cdot \bigg\{ \io \frac{1}{\ve} \bigg\}^\frac{\alpha+2}{2}
    + c_1(T),
  \eea
  where again combining Lemma \ref{lem14} with (\ref{277.11}) and Young's inequality provides $c_2(T)>0$ fulfilling
  \bea{277.3}
    & & \hspace*{-20mm}
    c_1(T) \cdot \bigg\{ - \io \ln\ue \bigg\}_+^\frac{\alpha+2}{2}
    + c_1(T) \eps^\frac{\alpha+2}{2} \cdot \bigg\{ \io \frac{1}{\ue^2} \bigg\}^\frac{\alpha+2}{2} \nn\\
    &\le& c_2(T) \io \frac{\uex^2}{\ue^2}
    + c_2(T) \eps^\frac{\alpha+2}{2} \io \ue^{-\alpha-4} \uex^2
    + c_2(T)
    \qquad \mbox{for all $t\in (0,T)$ and } \eps\in (0,1).
  \eea
  Another application of Lemma \ref{lem14} i), now to $p=1$ and $q=\alpha+1$, reveals that again due to (\ref{277.11})
  and Young's inequality we can moreover find $c_3(T)>0$ such that
  \bea{277.4}
    \hspace*{-10mm}
    c_1(T) \eps^\frac{\alpha+2}{2} \cdot \bigg\{ \io \frac{1}{\ve} \bigg\}^\frac{\alpha+2}{2}
    &\le& c_1(T) \eps^\frac{\alpha+2}{2} \cdot \Bigg\{
    (\alpha+1)^\frac{2}{\alpha+1} |\Omega|^\frac{\alpha+2}{\alpha+1} \cdot
    \bigg\{ \io \ve^{-\alpha-3} \vex^2 \bigg\}^\frac{1}{\alpha+1} \nn\\
    & & \hspace*{24mm}
    + 2^\frac{2}{\alpha+1} |\Omega|^2 \cdot
    \bigg\{ \io \ve \bigg\}^{-1} \Bigg\}^\frac{\alpha+2}{2} \nn\\
    &\le& c_3(T) \eps^\frac{\alpha+2}{2} \cdot \bigg\{ \io \ve^{-\alpha-3}\vex^2 \bigg\}^\frac{\alpha+2}{2(\alpha+1)}
    + c_3(T) \eps^\frac{\alpha+2}{2} \nn\\
    &\le& c_3(T) \eps^\frac{\alpha+2}{2} \io \ve^{-\alpha-3} \vex^2
    + 2c_3(T)
    \qquad \mbox{for all $t\in (0,T)$ and } \eps\in (0,1),
  \eea
  because $\frac{\alpha+2}{2(\alpha+1)} \le 1$.
  Finally, the Gagliardo-Nirenberg inequality along with (\ref{277.11}) and Young's inequality ensures the existence
  of $c_4>0$ and $c_5(T)>0$ such that for all $t\in (0,T)$ and any $\eps\in (0,1)$,
  \bas
    c_1(T) \cdot \bigg\{ \io \ve^2 \bigg\}^\frac{\alpha+2}{2}
    &\le& c_1(T) \cdot \bigg\{ c_4 \|\vex\|_{L^2(\Omega)}^\frac{\alpha+2}{3} \|\ve\|_{L^1(\Omega)}^\frac{2(\alpha+2)}{3}
    + c_4 \|\ve\|_{L^1(\Omega)}^{\alpha+2} \bigg\} \\
    &\le& c_5(T) \|\vex\|_{L^2(\Omega)}^\frac{\alpha+2}{3}
    + c_5(T) \nn\\
    &\le& c_5(T) \io \vex^2 + 2c_5(T),
  \eas
  for $\frac{\alpha+2}{3}\le 2$.
  Together with (\ref{277.3}), (\ref{277.4}) and (\ref{277.2}), this immediately leads to (\ref{277.1}).
\qed
Therefore, Lemma \ref{lem255} implies decay in a yet very weak form similar to that in Lemma \ref{lem26}.
\begin{lem}\label{lem266}
  Let $n_1=2$ and $n_2=1$, let
  $D_i>0, a_i>0, \lambda_i>0$ and $\chi_i>0$ for $i\in\{1,2\}$, assume (\ref{233.1}),
  and let $\chi_1\in (0,\chiss)$ and $\chi_2\in (0,\chiss)$ with $\chiss>0$ as in Lemma \ref{lem255}.
  Then whenever (\ref{ie}) holds, there exists $T_0=T_0(\chi_1,\chi_2,u_0,v_0)>0$ such that
  $u$ and $v$ from Lemma \ref{lem_conv} satisfy
  \be{266.1}
    \int_{T_0}^\infty \io u_x^2
    + \int_{T_0}^\infty \io v_x^2
    + \int_{T_0}^\infty \io (u-\lambda_1)^2
    + \int_{T_0}^\infty \io v^3
    <\infty.
  \ee
\end{lem}
\proof
  Once more thanks to Lemma \ref{lem62}, this can be seen by exploiting Lemma \ref{lem255} together with Lemma \ref{lem277}
  in essentially the same way as Lemma \ref{lem25} and Lemma \ref{lem27} have been used in the derivation of
  Lemma \ref{lem26}.
\qed
With this information returning to the inequality from Lemma \ref{lem255} yields stabilization in a sense paralleling
that of Lemma \ref{lem31}.
\begin{lem}\label{lem311}
  Let $n_1=2$ and $n_2=1$, let
  $D_i>0, a_i>0, \lambda_i>0$ and $\chi_i>0$ for $i\in\{1,2\}$ be such that (\ref{233.1}) is valid,
  and let $\chi_1\in (0,\chiss)$ and $\chi_2\in (0,\chiss)$ with $\chiss>0$ as in Lemma \ref{lem25}.
  Then assuming (\ref{ie}) and letting $u,v$ and $N$ be as given by Lemma \ref{lem_conv}
  and $\phi_{\lambda_1}$ be as defined through (\ref{phi}), we have
  \be{311.1}
    \io \phi_{\lambda_1}(u(\cdot,t))\to 0
    \quad \mbox{and} \quad
    \io v^3(\cdot,t) \to 0
    \qquad \mbox{as } (0,\infty)\setminus N \ni t \to\infty.
  \ee
\end{lem}
\proof
  A verification of this can be achieved by adapting the proof of Lemma \ref{lem31} in an obvious manner.
\qed
Finally, by compactness the latter can be turned into uniform convergence as in Lemma \ref{lem32}.
\begin{lem}\label{lem322}
  Let $n_1=2$ and $n_2=1$, let
  $D_i>0, a_i>0, \lambda_i>0$ and $\chi_i>0$ for $i\in\{1,2\}$ be such that (\ref{233.1}) holds,
  and let $\chi_1\in (0,\chiss)$ and $\chi_2\in (0,\chiss)$ with $\chiss>0$ as in Lemma \ref{lem255}.
  Then whenever (\ref{ie}) is satisfied, the limit functions $u$ and $v$ in Lemma \ref{lem_conv} have the properties that
  \bas
    u(\cdot,t) \to \lambda_1
    \quad \mbox{in } L^\infty(\Omega)
    \qquad \mbox{and} \qquad
    v(\cdot,t) \to 0
    \quad \mbox{in } L^\infty(\Omega)
  \eas
  as $t\to\infty$.
\end{lem}
\proof
  Again relying on Corollary \ref{cor24} and Lemma \ref{lem29}, one can readily obtain this as a consequence of
  Lemma \ref{lem311} by means of an argument in the flavor of that presented in the proof of Lemma \ref{lem32}.
\qed
We have thereby established our main result on asymptotic dominance of the predator population when
$\lambda_2 \le a_2\lambda_1$.\abs
\proofc of Theorem \ref{theo333}. \quad
  All statements have been verified in Corollary \ref{cor24} and Lemma \ref{lem322}.
\qed

\bigskip

{\bf Acknowledgement.} \quad
  Youshan Tao was supported by the {\em National Natural Science Foundation of China
  (No. 11861131003)}. The second author acknowledges support of the {\em Deutsche Forschungsgemeinschaft}
  in the context of the project
  {\em Emergence of structures and advantages in cross-diffusion systems} (No.~411007140, GZ: WI 3707/5-1).

\end{document}